\documentclass[11pt,reqno]{amsart}
\usepackage{amsmath,amssymb,amsthm}
\usepackage[a4paper,margin=27mm]{geometry}
\usepackage{cite}
\usepackage{url}
\usepackage{microtype}
\usepackage[colorlinks=true,linkcolor=blue,citecolor=red,urlcolor=blue]{hyperref}

\theoremstyle{plain}
\newtheorem{Thm}{Theorem}[section]
\newtheorem{Lem}{Lemma}[section]
\newtheorem{Prop}{Proposition}[section]
\newtheorem{Cor}{Corollary}[section]
\theoremstyle{remark}
\newtheorem{Rem}{Remark}[section]

\numberwithin{equation}{section}
\allowdisplaybreaks[2]
\newcommand{\dd}{\,\mathrm{d}}

\title[Stability of spiky patterns in a half-cylinder]{Nonlinear Stability of Spiky Patterns for a Three-Dimensional Chemotaxis--Consumption Model with Signal-Dependent Motilities}
\author[Lin Chang]{Lin Chang}
\address{Handan University, Hebei, China}
\email{changlin23@163.com}
\author[Song Xu]{Song Xu}
\address{School of Mathematical Sciences, South China Normal University, Guangzhou, China}
\email{songxu@m.scnu.edu.cn}
\author[Yuling Xu]{Yuling Xu}
\address{Handan Vocational College of Science and Technology, Hebei, China}
\email{m15122197516@163.com}
\thanks{Corresponding author: Yuling Xu.}
\date{}

\begin{document}
\begin{abstract}
We prove the nonlinear asymptotic stability of one-dimensional spiky steady states for a three-dimensional chemotaxis--consumption system with signal-dependent motility $w^p$ ($p>0$, $\alpha\ge2$) in the half-cylinder $\Omega=\mathbb{R}_+\times\mathbb{T}^2$. The cell density satisfies a no-flux boundary condition at $x_1=0$, while the chemical concentration satisfies $w(0,x',t)=b>0$ and vanishes as $x_1\to+\infty$. For any prescribed steady-state cell mass $\lambda>0$, we show that any sufficiently small compatible initial perturbation with zero mass generates a unique global strong solution that converges to the spiky steady state in $L^\infty(\Omega)$, with the chemical concentration converging uniformly relative to the steady profile, namely $\|w(\cdot,t)/W-1\|_{L^\infty(\Omega)}\to0$ as $t\to\infty$.

The stability analysis for an arbitrarily large steady mass $\lambda>0$ in three dimensions involves four main difficulties: the algebraic degeneracy of the diffusion coefficient $W^p(x_1)\sim (1+x_1)^{-2/(\alpha-1)}$ as $x_1\to\infty$, the $O(\lambda)$ linear coupling between the density perturbation and the chemical logarithmic gradient on the transversal nonzero modes, the nonlinear boundary conditions at $x_1=0$ produced by the Hopf--Cole transformation, and the different spatial weights required across derivative orders. We resolve these difficulties by introducing the algebraic weights $(\omega_1,\omega_2,\omega_3)=(W^{-p}U^{-1},U^{-1},W^pU^{-1})$, coupling a weighted Neumann inverse on the transversal nonzero modes with the scalar logarithmic potential $\Psi+q$ to cancel the $O(\lambda)$ linear coupling terms exactly, performing tangential and time-differentiated weighted energy estimates whose boundary traces vanish identically at $x_1=0$, and closing the time-derivative energy before recovering the full weighted Hessian dissipation from the parabolic system.
\end{abstract}
\maketitle

\section{Introduction}
We consider the chemotaxis--consumption system with signal-dependent motilities in the half-cylinder $\Omega=\mathbb{R}_+\times\mathbb{T}^2$:
\begin{equation}\label{eq:model}
\left\{
\begin{aligned}
u_t &= \nabla\cdot\left(w^p\nabla u-(\alpha-1)puw^p\nabla\log w\right),\\
w_t &= D\Delta w-wu,
\end{aligned}
\right.
\end{equation}
for $x=(x_1,x')\in\Omega$ and $t>0$, where $x'=(x_2,x_3)\in\mathbb{T}^2$. Here $u(x,t)\ge0$ is the cell density, $w(x,t)>0$ is the chemical concentration, $D>0$ is the chemical diffusion coefficient, and $w^p$ is the cell motility. Throughout this paper, we assume $p>0$ and $\alpha\ge2$. We prescribe the initial data
\begin{equation}\label{eq:initial-data}
(u,w)(x,0)=(u_0,w_0)(x),\qquad x\in\Omega,
\end{equation}
periodic boundary conditions in $x'\in\mathbb{T}^2$, and the normal boundary and far-field conditions
\begin{equation}\label{eq:original-boundary}
\left\{
\begin{aligned}
&\left(u_{x_1}-(\alpha-1)pu(\log w)_{x_1}\right)(0,x',t)=0,\\
&w(0,x',t)=b,\\
&(u,w)(x_1,x',t)\to(0,0),\qquad x_1\to+\infty,
\end{aligned}
\right.
\end{equation}
where $b>0$ is a given constant. The no-flux boundary condition at $x_1=0$ prevents cells from crossing the boundary, while the Dirichlet condition $w(0,x',t)=b$ supplies the chemical signal consumed in the interior.

Chemotaxis--consumption systems describe the directed migration of cells toward higher concentrations of a chemical signal depleted by the cells; see the survey by Lankeit and Winkler~\cite{LankeitWinkler2023} for results on global existence and long-time behavior. In the half-cylinder $\Omega$ under \eqref{eq:original-boundary}, the balance between chemotactic aggregation and degenerate diffusion determines a family of nonconstant steady states $(U(x_1),W(x_1))$ parametrized by the cell mass
\[
\lambda = \int_0^\infty U(x_1)\dd x_1 > 0.
\]
As recalled in Section~\ref{sec:preliminaries}, both $U(x_1)$ and $W(x_1)$ decay algebraically as $x_1\to\infty$ and concentrate near $x_1=0$.

The stability of boundary spikes in one space dimension has been studied in recent works. For a half-line model with singular motility, Song and Li~\cite{SongLi2022} proved the nonlinear stability of spiky steady states via the Hopf--Cole transformation and weighted energy estimates. For the power-law motility $w^p$ in \eqref{eq:model}, Song, Li, and Zhang~\cite{XLZ2024} established the one-dimensional nonlinear stability of $(U,W)$ on $\mathbb{R}_+$ by combining antiderivative weighted energy estimates with one-dimensional elliptic estimates; see also~\cite{GuoHouQiuXu2024,LiWangSongXu2026,LiuQiuWangXu2025,LiuWangXu2026JMPA} for recent developments on antiderivative, weighted energy, and Cole--Hopf methods for one-dimensional wave patterns. In the three-dimensional half-cylinder $\Omega=\mathbb{R}_+\times\mathbb{T}^2$, perturbations depend on the transversal variables $x'\in\mathbb{T}^2$, so the one-dimensional antiderivative method does not apply directly to the full perturbation. Related multi-dimensional stability problems have been extensively studied for planar waves and shear flows in hyperbolic--parabolic conservation laws and kinetic equations. Recent works by Xu and collaborators have established the nonlinear stability and time-decay rates of planar shock waves for the three-dimensional compressible Navier--Stokes and Navier--Stokes--Fourier equations, the three-dimensional Boltzmann equation, and multi-dimensional scalar viscous conservation laws~\cite{ChangLiuXu2025,DengXu2026,HuangXu2026,LiuWangXu2026}, planar rarefaction waves and their vanishing viscosity limits~\cite{DuanHuangLiXu2026Rarefaction,HouLiuWangXu2024,HuangXu2025,HuangXuYuan2022}, planar contact and entropy waves in the compressible Navier--Stokes and Landau equations~\cite{DuanHuangLiXu2025,DuanHuangLiXu2026Landau,HuangLiWuXu2026}, compressible vortex sheets with viscosity~\cite{HuangXinXuYuan2023}, multi-dimensional composite waves~\cite{HouXu2023}, and three-dimensional compressible Couette flows~\cite{LiWangXuZhang2026}. In periodic channels such as $\mathbb{R}\times\mathbb{T}^2$ or $\mathbb{R}_+\times\mathbb{T}^2$, decomposing a three-dimensional perturbation into its transversal zero mode and nonzero mode over $\mathbb{T}^2$ separates the one-dimensional conservation law governing the transversal average from the nonzero-mode component controlled via the transversal Poincar\'e inequality; see~\cite{ChangLiuXu2025,DengXu2026,DuanHuangLiXu2025,DuanHuangLiXu2026Landau,HouLiuWangXu2024,HouXu2023,HuangXuYuan2022,LiWangXuZhang2026,LiuWangXu2026}. For the chemotaxis--consumption system \eqref{eq:model}, however, the degenerate diffusion $w^p$, the logarithmic sensitivity, and the large-mass chemotactic drift require a different weighted structure. In this paper, we establish the global nonlinear stability of the spiky steady state $(U,W)$ in $\mathbb{R}_+\times\mathbb{T}^2$ for any prescribed mass $\lambda>0$.

\subsection{Main difficulties and proof strategy}
Under the Hopf--Cole transformation $v=-\nabla\log w$, the steady state $(U,W)$ becomes $(U,Ve_1)$ with $V=-(\log W)_{x_1}$, and we study the perturbation $(\phi,\psi)=(u-U,v-Ve_1)$. Proving global stability in three dimensions without smallness on the mass $\lambda>0$ requires addressing four main analytical difficulties.

\medskip
\noindent\emph{1. Degenerate diffusion at infinity and the algebraic weight structure.}
Since the steady profile decays algebraically as $x_1\to\infty$ with
\[
W(x_1)\sim (1+x_1)^{-\frac{2}{(\alpha-1)p}}, \qquad U(x_1)\sim (1+x_1)^{-2}, \qquad V(x_1)\sim (1+x_1)^{-1},
\]
the diffusion coefficient $W^p(x_1)\sim (1+x_1)^{-2/(\alpha-1)}$ in the density equation degenerates at infinity, and unweighted $L^2$ energy estimates do not close. Because spatial derivatives and steady-state coefficients carry different algebraic decay rates, we introduce three weights:
\[
\omega_1=\frac{1}{W^pU},\qquad \omega_2=\frac{1}{U},\qquad \omega_3=\frac{W^p}{U}.
\]
For $\alpha\ge2$, all three weights have positive lower bounds and satisfy
\[
\omega_2^2=\omega_1\omega_3,\qquad \|W^p\nabla^2\phi\|_{\omega_1}=\|\nabla^2\phi\|_{\omega_3}.
\]
In our energy norm, $\phi$ and its tangential derivatives $\partial_{x'}^\beta\phi$ are measured in $L^2_{\omega_1}$, $\nabla\phi$ and $\psi,\nabla\psi,\nabla^2\psi$ are measured in $L^2_{\omega_2}$, and $\nabla^2\phi$ is measured in $L^2_{\omega_3}$.

\medskip
\noindent\emph{2. Large-mass linear coupling on transversal nonzero modes.}
Motivated by the transversal mode decomposition and weighted elliptic multiplier methods developed for multi-dimensional planar waves~\cite{ChangLiuXu2025,DuanHuangLiXu2025,DuanHuangLiXu2026Landau,DuanHuangLiXu2026Rarefaction,HuangLiWuXu2026,HuangXu2026,HuangXuYuan2022}, we split the perturbation into its transversal average $P_0$ and transversal nonzero mode $P_{\neq}=I-P_0$. The zero mode $(P_0\phi,P_0\psi_1)$ satisfies a one-dimensional conservation law on $\mathbb{R}_+$; under the zero-mass conditions \eqref{eq:zero-mass-compatibility}, we introduce the antiderivatives $(\Phi,\Psi)$ and close their weighted estimates via the Hardy inequality. On the nonzero mode $P_{\neq}\phi$, however, no one-dimensional antiderivative exists. In a direct $L^2_{\omega_1}$ estimate of $\phi$, the linearized chemotactic drift $\nabla\cdot(W^pU\psi)$ produces an $O(\lambda)$ linear coupling term between $\nabla\phi$ and $\psi$, which cannot be absorbed as a small error when $\lambda>0$ is large.

To overcome this difficulty for arbitrary $\lambda>0$, we use the curl-free relation $\nabla\times\psi=0$ to write $\psi$ as the spatial gradient of a scalar logarithmic potential:
\[
\psi=\nabla(\Psi+q),\qquad P_{\neq}(\Psi+q)=q,\qquad \log\frac{w}{W}=-(\Psi+q).
\]
We then define an auxiliary weighted Neumann potential $z$ on the nonzero-mode space $P_0z=0$ by solving the degenerate elliptic problem
\[
-(\alpha-1)p\,\nabla\cdot(W^pU\nabla z)=P_{\neq}\phi\quad\text{in }\Omega,\qquad z_{x_1}(0,x',t)=0.
\]
Testing the density equation against $z$ yields two key properties: the linear diffusion term produces the weighted density dissipation $\|P_{\neq}\phi\|_{\omega_2}^2$, while the linear chemotactic drift term becomes
\[
-\int_\Omega P_{\neq}\phi\,q\dd x.
\]
When added to the $L^2(\Omega)$ energy identity of the scalar potential equation for $q$, this linear coupling term cancels identically, yielding the closed nonzero-mode dissipation
\[
\|P_{\neq}\phi\|_{\omega_2}^2+\|P_{\neq}\psi\|_{L^2(\Omega)}^2
\]
for any $\lambda>0$.

\medskip
\noindent\emph{3. Nonlinear boundary conditions at $x_1=0$.}
Under $v=-\nabla\log w$, the boundary conditions \eqref{eq:original-boundary} become homogeneous tangential conditions $\psi_2=\psi_3=0$ and two nonlinear normal boundary conditions at $x_1=0$:
\[
\phi_{x_1}+(\alpha-1)p(V\phi+U\psi_1+\phi\psi_1)=0,\qquad D\psi_{1x_1}=2DV\psi_1+D\psi_1^2-\phi.
\]
Because of these nonlinear normal boundary conditions, normal derivatives $\partial_{x_1}\phi$ and $\partial_{x_1}\psi$ cannot be used as energy multipliers without generating higher-order boundary terms at $x_1=0$. Instead, the boundary conditions are preserved under tangential derivatives $\partial_{x'}^\beta$ ($1\le|\beta|\le2$) and time differentiation $\partial_t$. When testing the tangentially or time-differentiated density equation against $\partial_{x'}^\beta\phi/(W^pU)$ or $\phi_t/(W^pU)$, the boundary term at $x_1=0$ is the differentiated normal flux and vanishes identically. Similarly, when testing the differentiated vector equation against $\partial_{x'}^\beta\psi/U$ or $\psi_t/U$, the normal diffusion boundary term $-D\psi_{1x_1}$ combines with the drift, nonlinear, and density boundary terms to vanish at $x_1=0$ by the second boundary condition and $\psi_2=\psi_3=0$.

\medskip
\noindent\emph{4. Independent closure of time derivatives and recovery of the weighted Hessian.}
At the second derivative level, estimating the second-order tangential derivatives $\partial_{x'}^\beta\phi$ ($|\beta|=2$) involves trilinear commutator terms such as
\[
\frac{(\partial_{x_j}\psi_k)\nabla\phi}{\sqrt{U}}, \qquad \frac{(\partial_{x_j}\phi)\partial_{x_k}\psi}{\sqrt{U}},
\]
which require the full space--time Hessian integral
\[
\int_0^t\left(\|\nabla^2\phi\|_{\omega_3}^2+\|\nabla^2\psi\|_{\omega_2}^2\right)\dd\tau.
\]
To close the energy estimates without circularity, we estimate the time derivatives before the second-order tangential derivatives. After closing the lower-order estimate and the first-order tangential estimate ($|\beta|=1$), we combine the energy identity for $(\phi_t,\psi_t)$ with the identity obtained by testing the perturbation system against $(\phi_t/(W^pU),\psi_t/U)$. Controlling $\|(\Psi+q)_t/\sqrt{U}\|_{L^6(\Omega)}$ via the scalar potential equation closes
\[
\|\phi_t(t)\|_{\omega_1}^2+\|\psi_t(t)\|_{\omega_2}^2+\int_0^t\left(\|\phi_\tau\|_{\omega_1}^2+\|\psi_\tau\|_{\omega_2}^2+\|\nabla\phi_\tau\|_{\omega_2}^2+\|\nabla\psi_\tau\|_{\omega_2}^2\right)\dd\tau
\]
independently of the second-order spatial derivatives. We then solve the parabolic equations for $(\phi_{x_1x_1},\partial_{x_1}^2\psi)$ and combine them with the first-order tangential dissipation to bound the full space--time Hessian integral $\int_0^t(\|\nabla^2\phi\|_{\omega_3}^2+\|\nabla^2\psi\|_{\omega_2}^2)\dd\tau$. This bound closes the second-order tangential estimate ($|\beta|=2$), and current-time elliptic and interpolation estimates complete the closure of $N(t)$.

\medskip
The paper is organized as follows. Section~\ref{sec:preliminaries} introduces the steady profile, the perturbation variables, the weighted spaces, and Theorem~\ref{thm:main}. Section~\ref{sec:reformulation} gives the transformed system, the zero-mode mass conservation, the potential representation, and the zero-mode and scalar potential equations. Section~\ref{sec:apriori} establishes the a priori weighted energy estimates. Section~\ref{sec:proof-main} states the local continuation result and proves Theorem~\ref{thm:main}. Appendix~\ref{app:weighted} collects the weighted functional inequalities, and Appendix~\ref{app:neumann} proves the existence and flux bounds for the weighted nonzero-mode Neumann problem.

\medskip
\noindent\textbf{Notation.} We normalize the measure on $\mathbb{T}^2$ by $\int_{\mathbb{T}^2}1\dd x'=1$ and write $e_1=(1,0,0)^\top$ and $\nabla_{x'}=(\partial_{x_2},\partial_{x_3})$. For a multi-index $\beta=(\beta_2,\beta_3)\in\mathbb{N}_0^2$, $|\beta|=\beta_2+\beta_3$ and $\partial_{x'}^\beta=\partial_{x_2}^{\beta_2}\partial_{x_3}^{\beta_3}$. The symbol $\|\cdot\|_{L^r}$ (or $\|\cdot\|_r$ for $1\le r\le\infty$) denotes the unweighted $L^r(\Omega)$ or $L^r(\mathbb{R}_+)$ norm. We write $A\sim B$ when $cB\le A\le CB$ for positive constants $c,C$ depending only on $(p,\alpha,D,b,\lambda)$. Throughout the paper, $C$ and $c$ denote positive constants depending only on $(p,\alpha,D,b,\lambda)$ unless otherwise stated, and all constants in the a priori estimates are independent of $T$.

\section{Preliminaries and main theorem}\label{sec:preliminaries}\label{sec:main-result}
\subsection{The spiky steady state}
For any prescribed cell mass $\lambda>0$, the one-dimensional steady state $(U(x_1),W(x_1))$ solves the boundary value problem
\begin{equation}\label{eq:steady-problem}
\left\{
\begin{aligned}
&\left[W^p\left(U_{x_1}-(\alpha-1)pU(\log W)_{x_1}\right)\right]_{x_1}=0,\\
&DW_{x_1x_1}-WU=0,\\
&\int_0^\infty U(x_1)\dd x_1=\lambda,
\end{aligned}
\right.
\end{equation}
subject to the boundary and far-field conditions $W(0)=b$, $(U,W)(+\infty)=(0,0)$, and
\[
\left(U_{x_1}-(\alpha-1)pU(\log W)_{x_1}\right)(0)=0.
\]
Integrating \eqref{eq:steady-problem} yields the explicit algebraic spiky profile (see \cite{XLZ2024}):
\begin{equation}\label{eq:explicit-steady-profile}
\begin{aligned}
U(x_1) &= \frac{\lambda^2(\alpha-1)^2p^2}{2D[(\alpha-1)p+2]} \left(1+\frac{\lambda(\alpha-1)^2p^2}{2D[(\alpha-1)p+2]}x_1\right)^{-2},\\
W(x_1) &= b\left(1+\frac{\lambda(\alpha-1)^2p^2}{2D[(\alpha-1)p+2]}x_1\right)^{-\frac{2}{(\alpha-1)p}}.
\end{aligned}
\end{equation}
Under the Hopf--Cole transformation, the steady state corresponds to $(U,Ve_1)$, where
\begin{equation}\label{eq:steady-V}
V(x_1) = -(\log W)_{x_1} = \frac{\lambda(\alpha-1)p}{D[(\alpha-1)p+2]} \left(1+\frac{\lambda(\alpha-1)^2p^2}{2D[(\alpha-1)p+2]}x_1\right)^{-1}.
\end{equation}
Direct computation shows that $(U,V)$ satisfies the structural identities
\begin{equation}\label{eq:steady-profile-identities}
U_{x_1}+(\alpha-1)pUV=0, \qquad V^2=\frac{2U}{D[(\alpha-1)p+2]}, \qquad DV_{x_1}-DV^2+U=0.
\end{equation}
In particular, $-V_{x_1}=\frac{(\alpha-1)p}{D[(\alpha-1)p+2]}U>0$.

\subsection{Perturbation variables and transversal modes}
For any positive chemical concentration $w>0$, we define the Hopf--Cole vector field
\begin{equation}\label{eq:hopf-cole}
v=-\nabla\log w,
\end{equation}
and the perturbations around $(U,Ve_1)$ by
\begin{equation}\label{eq:perturbations}
\phi=u-U, \qquad \psi=v-(V,0,0)^\top.
\end{equation}
In terms of the relative chemical concentration $w/W$, we have
\begin{equation}\label{eq:psi-gradient-log}
\psi=-\nabla\log\frac{w}{W},
\end{equation}
which immediately implies the curl-free constraint
\begin{equation}\label{eq:curl-free}
\nabla\times\psi=0.
\end{equation}
For any integrable scalar or vector function $g$ on $\Omega$, we define the transversal zero-mode projection $P_0$ and nonzero-mode projection $P_{\neq}$ by
\begin{equation}\label{eq:mode-projections}
 P_0 g(x_1,t)=\int_{\mathbb{T}^2}g(x_1,x',t)\dd x', \qquad  P_{\neq} g=g- P_0 g.
\end{equation}
Since $\psi_j=-\partial_{x_j}\log(w/W)$ for $j=2,3$, periodicity in $x'\in\mathbb{T}^2$ implies that the tangential components of $\psi$ have zero transversal mean:
\begin{equation}\label{eq:psi23-zero-average}
 P_0\psi_j=-\int_{\mathbb{T}^2}\partial_{x_j}\log\frac{w}{W}\dd x'=0,\qquad j=2,3.
\end{equation}
Thus $\psi_2= P_{\neq}\psi_2$ and $\psi_3= P_{\neq}\psi_3$ are purely nonzero-mode quantities and obey the transversal Poincar\'e inequality. For the longitudinal zero modes $P_0\phi$ and $P_0\psi_1$, we introduce the antiderivatives
\begin{equation}\label{eq:Phi-Psi-definition}
(\Phi,\Psi)(x_1,t) = -\int_{x_1}^{+\infty} \left(P_0\phi, P_0\psi_1\right)(y_1,t)\dd y_1,
\end{equation}
so that
\begin{equation}\label{eq:Phi-Psi-derivatives}
\Phi_{x_1}= P_0\phi, \qquad \Psi_{x_1}= P_0\psi_1.
\end{equation}
For the nonzero mode $P_{\neq}\psi_1$, we define the potential
\begin{equation}\label{eq:tpsi-definition}
 q(x_1,x',t) = \int_0^{x_1} P_{\neq}\psi_1(y_1,x',t)\dd y_1,
\end{equation}
which satisfies $q(0,x',t)=0$, $P_0 q=0$, and
\begin{equation}\label{eq:tpsi-normal}
 q_{x_1}= P_{\neq}\psi_1.
\end{equation}
Evaluation at $t=0$ is denoted by a subscript $0$, such as $\phi_0$, $\psi_0=(\psi_{10},\psi_{20},\psi_{30})^\top$, $\Phi_0$, $\Psi_0$, and $q_0$.

We impose the zero-mass compatibility conditions on the initial perturbation:
\begin{equation}\label{eq:zero-mass-compatibility}
\int_0^\infty  P_0\phi_0(y_1)\dd y_1=0, \qquad \int_0^\infty  P_0\psi_{10}(y_1)\dd y_1=0,
\end{equation}
together with the initial boundary compatibility conditions at $x_1=0$:
\[
\left\{
\begin{aligned}
&w_0(0,x')=b, \qquad \psi_{20}(0,x')=\psi_{30}(0,x')=0,\\
&\phi_{0x_1}(0,x')+(\alpha-1)p\left(V(0)\phi_0(0,x')+U(0)\psi_{10}(0,x')+\phi_0(0,x')\psi_{10}(0,x')\right)=0,\\
&D\psi_{10x_1}(0,x')=2DV(0)\psi_{10}(0,x')+D\psi_{10}^2(0,x')-\phi_0(0,x').
\end{aligned}
\right.
\]
As shown in Lemmas~\ref{lem:zero-mass-preservation} and \ref{lem:potential-representation}, the zero-mass conditions \eqref{eq:zero-mass-compatibility} are preserved for all $t\ge0$, ensuring the homogeneous Dirichlet boundary conditions $\Phi(0,t)=\Psi(0,t)=0$ and the exact potential reconstruction $w=W\exp(-(\Psi+q))$.

\subsection{Weighted spaces}
We define the three algebraic weights
\begin{equation}\label{eq:weights}
\omega_1(x_1)=\frac{1}{W^p(x_1)U(x_1)},\qquad \omega_2(x_1)=\frac{1}{U(x_1)},\qquad \omega_3(x_1)=\frac{W^p(x_1)}{U(x_1)}.
\end{equation}
For $O\in\{\mathbb{R}_+,\Omega\}$, a positive weight $\omega(x_1)$, and an integer $k\ge0$, we define the weighted $L^2$ and Sobolev norms by
\[
\|f\|_{\omega,O}^2=\int_{O}\omega(x_1)|f|^2\dd x, \qquad \|f\|_{k,\omega,O}^2=\sum_{|\gamma|\le k}\|\partial^\gamma f\|_{\omega,O}^2,
\]
and denote the corresponding weighted Sobolev space by $H^k_\omega(O)$. The domain subscript $O$ is omitted when clear from context. Since $\alpha\ge2$ and $p>0$, the explicit formulas \eqref{eq:explicit-steady-profile}--\eqref{eq:steady-V} imply that $(1+x_1)^2\sim \omega_2(x_1)$ and
\begin{equation}\label{eq:weight-lower-bounds}
\omega_j(x_1)\ge c>0,\qquad \omega_3(x_1)\le C\omega_2(x_1)\le C\omega_1(x_1), \qquad j=1,2,3,
\end{equation}
as well as the uniform profile bounds
\begin{equation}\label{eq:profile-pointwise-bounds}
\left|\frac{U_{x_1}}{U}\right|\le C,\qquad \frac{V^2}{U}\le C,\qquad |V_{x_1}|\le CU,\qquad W^p\le C.
\end{equation}
Moreover, the weights satisfy the exact scaling relations
\[
\omega_2^2=\omega_1\omega_3,\qquad \|W^p\nabla^2\phi\|_{\omega_1}^2=\|\nabla^2\phi\|_{\omega_3}^2.
\]
To track the perturbation in our bootstrap argument, we define the weighted energy functional
\begin{equation}\label{eq:bootstrap-norm}
\begin{aligned}
N^2(t) = \sup_{0\le\tau\le t}\Bigg[
&\|\Phi(\tau)\|_{1,\omega_1}^2 + \|\Psi(\tau)\|_{L^2(\mathbb{R}_+)}^2 + \|\Psi_{x_1}(\tau)\|_{1,\omega_2}^2 + \|\phi(\tau)\|_{\omega_1}^2 + \|\nabla\phi(\tau)\|_{\omega_2}^2 \\
&+ \|\nabla^2\phi(\tau)\|_{\omega_3}^2 + \sum_{1\le|\beta|\le2} \|\partial_{x'}^\beta\phi(\tau)\|_{\omega_1}^2 + \|\psi(\tau)\|_{\omega_2}^2 + \|\nabla\psi(\tau)\|_{1,\omega_2}^2
\Bigg].
\end{aligned}
\end{equation}
For $T>0$, we define the solution space $X(0,T)$ to be the class of perturbations satisfying the boundary, potential, and zero-mass conditions above, together with the regularity
\begin{equation}\label{eq:solution-space-zero-final}
\left\{
\begin{aligned}
&\Phi\in C\big([0,T];H^2(\mathbb{R}_+)\cap H^1_{\omega_1}(\mathbb{R}_+)\big),\qquad \Phi_{x_1}\in L^2\big(0,T;H^1_{\omega_2}(\mathbb{R}_+)\big),\\
&\Psi\in C\big([0,T];H^2(\mathbb{R}_+)\big),\qquad \Psi_{x_1}\in C\big([0,T];H^1_{\omega_2}(\mathbb{R}_+)\big) \cap L^2\big(0,T;H^1_{\omega_2}(\mathbb{R}_+)\big),
\end{aligned}
\right.
\end{equation}
and
\begin{equation}\label{eq:solution-space-physical-final}
\left\{
\begin{aligned}
&(\phi,\psi)\in C\big([0,T];H^2(\Omega)\big),\qquad \phi\in C\big([0,T];L^2_{\omega_1}(\Omega)\big),\\
&\nabla\phi\in C\big([0,T];L^2_{\omega_2}(\Omega)\big) \cap L^2\big(0,T;L^2_{\omega_2}(\Omega)\big),\\
&\nabla^2\phi\in C\big([0,T];L^2_{\omega_3}(\Omega)\big) \cap L^2\big(0,T;L^2_{\omega_3}(\Omega)\big),\\
&\partial_{x'}^\beta\phi\in C\big([0,T];L^2_{\omega_1}(\Omega)\big), \qquad 1\le|\beta|\le2,\\
&\psi\in C\big([0,T];L^2_{\omega_2}(\Omega)\big),\qquad \nabla\psi\in C\big([0,T];H^1_{\omega_2}(\Omega)\big) \cap L^2\big(0,T;H^1_{\omega_2}(\Omega)\big).
\end{aligned}
\right.
\end{equation}
Note that $\Phi_{x_1x_1}= P_0(\phi_{x_1})$ and $\Psi_{x_1x_1}= P_0(\psi_{1x_1})$ are controlled in their respective weighted spaces by Jensen's inequality.

\subsection{Main theorem}
Our main stability result reads as follows.
\begin{Thm}\label{thm:main}
Let $p>0$, $\alpha\ge2$, $D>0$, $b>0$, and $\lambda>0$ be fixed, and let $(U,W)$ be the spiky steady state defined in \eqref{eq:explicit-steady-profile}. Assume that $u_0\ge0$, $w_0>0$ at every finite spatial point, and $w_0(0,x')=b$. Suppose that the initial perturbation satisfies $N(0)<\infty$ in \eqref{eq:bootstrap-norm} along with the initial boundary compatibility conditions, the logarithmic-gradient relation \eqref{eq:psi-gradient-log} at $t=0$, and the zero-mass conditions \eqref{eq:zero-mass-compatibility}. Then there exists a constant $\varepsilon_0>0$ such that if
\begin{equation}\label{eq:main-smallness}
N(0)\le \varepsilon_0,
\end{equation}
the initial-boundary value problem \eqref{eq:model}--\eqref{eq:original-boundary} admits a unique global solution $(u,w)$ whose perturbation belongs to $X(0,T)$ for every $T>0$ and satisfies
\begin{equation}\label{eq:main-uniform-bound}
\sup_{t\ge0}N(t)\le C N(0).
\end{equation}
Moreover, the solution converges uniformly to the steady state as $t\to\infty$:
\begin{equation}\label{eq:main-asymptotic}
\lim_{t\to\infty}\left(
\|u(\cdot,t)-U\|_{L^\infty(\Omega)}
+\left\|\frac{w(\cdot,t)}{W}-1\right\|_{L^\infty(\Omega)}
+\|-\nabla\log w(\cdot,t)-Ve_1\|_{L^\infty(\Omega)}
\right)=0.
\end{equation}
\end{Thm}

\begin{Rem}
The condition $\alpha\ge2$ ensures that $\omega_3(x_1)=W^p(x_1)/U(x_1)\sim (1+x_1)^{2-2/(\alpha-1)}\ge c>0$, so that the weighted second-order norm $\|\nabla^2\phi\|_{\omega_3}$ controls the unweighted norm $\|\nabla^2\phi\|_{L^2(\Omega)}$. The steady mass $\lambda>0$ can be arbitrarily large; the threshold $\varepsilon_0$ and constant $C$ depend on $\lambda$, and the complete space--time dissipation bound is stated in Corollary~\ref{cor:complete-closed-apriori}.
\end{Rem}

\section{Reformulation of the original problem}\label{sec:reformulation}
\subsection{The transformed initial-boundary value problem}
Differentiating the Hopf--Cole relation $v=-\nabla\log w$ with respect to $t$ and substituting the second equation of \eqref{eq:model} yields
\[
(\log w)_t = D\left(|v|^2-\nabla\cdot v\right)-u.
\]
Taking the spatial gradient, we rewrite the system \eqref{eq:model} in terms of $(u,v)$ as
\begin{equation}\label{eq:transformed-system}
\left\{
\begin{aligned}
u_t &= \nabla\cdot\left(w^p\big(\nabla u+(\alpha-1)puv\big)\right),\\
v_t &= -\nabla\left(D\big(|v|^2-\nabla\cdot v\big)-u\right),\\
w(x_1,x',t) &= b\exp\left(-\int_0^{x_1}v_1(y_1,x',t)\dd y_1\right),
\end{aligned}
\right.
\end{equation}
with initial data $(u,v)|_{t=0}=(u_0,-\nabla\log w_0)$. At the boundary $x_1=0$, the Dirichlet condition $w(0,x',t)=b$ implies $v_2(0,x',t)=v_3(0,x',t)=0$ and $(\log w)_t(0,x',t)=0$. Combining this with $\nabla\times v=0$ and the no-flux boundary condition in \eqref{eq:original-boundary} gives the transformed boundary conditions at $x_1=0$:
\begin{equation}\label{eq:transformed-boundary}
\left\{
\begin{aligned}
&u_{x_1}+(\alpha-1)puv_1=0,\\
&Dv_{1x_1}-Dv_1^2+u=0,\\
&v_2=v_3=0.
\end{aligned}
\right.
\end{equation}

\subsection{The physical perturbation and its boundary conditions}
Subtracting the steady-state equations \eqref{eq:steady-profile-identities} from \eqref{eq:transformed-boundary}, we obtain the boundary conditions for the perturbation $(\phi,\psi)$ at $x_1=0$:
\begin{equation}\label{eq:perturbation-bc}
\left\{
\begin{aligned}
&\phi_{x_1}+(\alpha-1)p\left(V\phi+U\psi_1+\phi\psi_1\right)=0,\\
&D\psi_{1x_1}=2DV\psi_1+D\psi_1^2-\phi,\\
&\psi_2=\psi_3=0,
\end{aligned}
\right.
\qquad x_1=0.
\end{equation}
In particular, the tangential components vanish identically on the boundary:
\begin{equation}\label{eq:tangential-boundary}
\psi_2(0,x',t)=\psi_3(0,x',t)=0.
\end{equation}
Using the curl-free identity $\nabla(\nabla\cdot\psi)=\Delta\psi$ from \eqref{eq:curl-free} and the representation $w=W\exp(-(\Psi+q))$ justified in Lemma~\ref{lem:potential-representation} below, we obtain the interior perturbation system
\begin{equation}\label{eq:physical-system}
\left\{
\begin{aligned}
\phi_t &= \nabla\cdot\Big\{W^p e^{-p(\Psi+q)} \big[\nabla\phi+(\alpha-1)pV\phi e_1 +(\alpha-1)p(U+\phi)\psi\big]\Big\},\\
\psi_t &= D\Delta\psi-2D\nabla(V\psi_1)-D\nabla|\psi|^2+\nabla\phi,
\end{aligned}
\right.
\end{equation}
with initial data $(\phi,\psi)|_{t=0}=(\phi_0,\psi_0)$. To avoid any circularity, the proof of zero-mode mass conservation below uses the original coefficient $w^p$ prior to invoking the exponential potential representation.

\begin{Lem}\label{lem:zero-mass-preservation}
Let $T>0$ and suppose $(\phi,\psi)\in X(0,T)$ satisfies the boundary conditions \eqref{eq:perturbation-bc}. If the initial zero-mass compatibility conditions \eqref{eq:zero-mass-compatibility} hold, then for every $t\in[0,T]$,
\begin{equation}\label{eq:zero-mass-preserved}
\int_0^\infty  P_0\phi(y_1,t)\dd y_1=0, \qquad \int_0^\infty  P_0\psi_1(y_1,t)\dd y_1=0,
\end{equation}
and consequently
\begin{equation}\label{eq:zero-mode-boundary}
\Phi(0,t)=\Psi(0,t)=0.
\end{equation}
\end{Lem}
\begin{proof}
By \eqref{eq:psi-gradient-log} and $w(0,x',\tau)=W(0)=b$, we have
\[
\log\frac{w(x_1,x',\tau)}{W(x_1)} = -\int_0^{x_1}\psi_1(y_1,x',\tau)\dd y_1.
\]
Using the Sobolev embedding $H^2(\mathbb{T}^2)\hookrightarrow L^\infty(\mathbb{T}^2)$ and the Cauchy--Schwarz inequality in $x_1$, we find
\[
\left\|\log\frac{w}{W}\right\|_{L^\infty(\Omega)}
\le C\left(\int_0^\infty U(x_1)\dd x_1\right)^{1/2}
\left(\sum_{|\beta|\le2}\int_\Omega \frac{|\partial_{x'}^\beta\psi_1|^2}{U}\dd x\right)^{1/2} < \infty.
\]
Thus $w^p$ is bounded on $[0,T]\times\Omega$. Averaging the density equation and the first component of the vector equation over $\mathbb{T}^2$ gives
\[
\begin{aligned}
\partial_\tau P_0\phi &= \partial_{x_1} P_0\Big\{w^p\big[\phi_{x_1}+(\alpha-1)pV\phi+(\alpha-1)p(U+\phi)\psi_1\big]\Big\},\\
\partial_\tau P_0\psi_1 &= \partial_{x_1} P_0\left(D\psi_{1x_1}-2DV\psi_1-D|\psi|^2+\phi\right).
\end{aligned}
\]
By \eqref{eq:perturbation-bc}, both normal fluxes inside the $\partial_{x_1}$ derivatives vanish at $x_1=0$. Since $\int_0^\infty U\dd x_1<\infty$ and $\int_0^\infty W^pU\dd x_1<\infty$, the Cauchy--Schwarz inequality and the weighted bounds in $X(0,T)$ imply that both fluxes belong to $L^1((0,T)\times\Omega)$. Let $\chi_R\in C_c^\infty([0,\infty))$ satisfy $\chi_R\equiv1$ on $[0,R]$, $\chi_R\equiv0$ on $[2R,\infty)$, $0\le\chi_R\le1$, and $|\chi_R'|\le C/R$. Multiplying the two conservation laws by $\chi_R(x_1)$, integrating over $[0,t]\times\mathbb{R}_+$, and letting $R\to\infty$ by dominated convergence yields \eqref{eq:zero-mass-preserved}. The boundary values \eqref{eq:zero-mode-boundary} then follow from \eqref{eq:Phi-Psi-definition}.
\end{proof}

\begin{Lem}\label{lem:potential-representation}
Let $T>0$ and suppose $(\phi,\psi)\in X(0,T)$ satisfies the boundary conditions \eqref{eq:perturbation-bc} and the curl-free constraint \eqref{eq:curl-free}. Then for every $t\in[0,T]$, the potentials $\Psi$ and $q$ satisfy
\begin{equation}\label{eq:w-decomposition}
w = W e^{-(\Psi+q)},
\end{equation}
\begin{equation}\label{eq:potential-identity}
\nabla(\Psi+q) = \psi,
\end{equation}
\begin{equation}\label{eq:potential-components}
\partial_{x_1}(\Psi+q) = \psi_1, \qquad \partial_{x_2} q = \psi_2, \qquad \partial_{x_3} q = \psi_3,
\end{equation}
and
\begin{equation}\label{eq:zero-nonzero-potential}
 P_0\psi = (\Psi_{x_1},0,0)^\top, \qquad  P_{\neq}\psi = \nabla q.
\end{equation}
\end{Lem}
\begin{proof}
By the definitions \eqref{eq:Phi-Psi-definition} and \eqref{eq:tpsi-definition}, we have $\Psi_{x_1}= P_0\psi_1$, $q_{x_1}= P_{\neq}\psi_1$, and $P_0 q=0$. Thus
\[
\partial_{x_1}(\Psi+q) =  P_0\psi_1+ P_{\neq}\psi_1 = \psi_1.
\]
For the tangential derivatives $j=2,3$, the curl-free condition \eqref{eq:curl-free} gives $\partial_{x_j}\psi_1=\partial_{x_1}\psi_j$. Since $\partial_{x_j} P_0\psi_1=0$ and $\psi_j(0,x',t)=0$ by \eqref{eq:tangential-boundary}, we have
\[
\partial_{x_j} q(x_1,x',t) = \int_0^{x_1}\partial_{x_j} P_{\neq}\psi_1(y_1,x',t)\dd y_1
= \int_0^{x_1}\partial_{x_1}\psi_j(y_1,x',t)\dd y_1 = \psi_j(x_1,x',t).
\]
This establishes \eqref{eq:potential-identity} and \eqref{eq:potential-components}. Combining \eqref{eq:potential-identity} with $\psi=-\nabla\log(w/W)$ from \eqref{eq:psi-gradient-log} shows that
\[
\nabla\left(\log\frac{w}{W}+\Psi+q\right)=0.
\]
Evaluating at $x_1=0$ and using $w(0,x',t)=W(0)=b$, $\Psi(0,t)=0$ (by Lemma~\ref{lem:zero-mass-preservation}), and $q(0,x',t)=0$, we deduce that $\log(w/W)+\Psi+q=0$ throughout $\Omega$, which gives \eqref{eq:w-decomposition}. Finally, applying $P_0$ and $P_{\neq}$ to \eqref{eq:potential-identity} yields \eqref{eq:zero-nonzero-potential}.
\end{proof}

\begin{Rem}
The second condition in \eqref{eq:zero-mass-compatibility} normalizes the far-field behavior of $w_0$ relative to $W$. Specifically, since $-\log(w_0/W)|_{x_1=0}=0$, we have $\int_0^R P_0\psi_{10}\dd x_1=- P_0\log(w_0(R,\cdot)/W(R))$, so $\int_0^\infty P_0\psi_{10}\dd x_1=0$ is equivalent to $\lim_{R\to\infty} P_0\log(w_0(R,\cdot)/W(R))=0$. Moreover, the weighted bounds in Lemma~\ref{lem:mode-reconstruction} imply
\[
|\Psi_0(x_1)| \le \| P_0\psi_{10}\|_{\omega_2}\left(\int_{x_1}^\infty U(y_1)\dd y_1\right)^{1/2}\to0, \qquad | q_0(x)|\le CN(0)\sqrt{U(x_1)}\to0
\]
as $x_1\to\infty$, so that $w_0(x_1,x')/W(x_1)=\exp(-(\Psi_0(x_1)+q_0(x_1,x')))\to1$ uniformly in $x'\in\mathbb{T}^2$.
\end{Rem}

\subsection{The equations for the zero mode and scalar potential}
Averaging the physical perturbation system \eqref{eq:physical-system} over $\mathbb{T}^2$ and integrating with respect to $x_1$ over $[x_1,\infty)$ via \eqref{eq:Phi-Psi-definition}, we obtain the one-dimensional zero-mode system for $(\Phi,\Psi)$:
\begin{equation}\label{eq:zero-mode-system}
\left\{
\begin{aligned}
W^{-p}\Phi_t &= e^{-p\Psi}\left(\Phi_{x_1x_1}+(\alpha-1)pV\Phi_{x_1}\right)+(\alpha-1)pU\Psi_{x_1}\\
&\quad +(\alpha-1)pe^{-p\Psi}\Phi_{x_1}\Psi_{x_1}+\left(e^{-p\Psi}-1\right)(\alpha-1)pU\Psi_{x_1}+W^{-p} P_0 R_1,\\
\Psi_t &= D\Psi_{x_1x_1}-2DV\Psi_{x_1}-D\Psi_{x_1}^2+\Phi_{x_1}+ P_0 R_2,\\
\Phi(0,t) &= 0,\qquad \Psi(0,t)=0,
\end{aligned}
\right.
\end{equation}
where the nonzero-mode coupling terms $R_1$ and $R_2$ are given by
\begin{equation}\label{eq:zero-mode-remainders}
\begin{aligned}
R_1 &= W^pe^{-p\Psi}\Big\{\big(e^{-p q}-1\big)\big[\phi_{x_1}+(\alpha-1)p(V\phi+U\psi_1+\phi\psi_1)\big]
+(\alpha-1)p\, P_{\neq}\phi\, P_{\neq}\psi_1\Big\},\\
R_2 &= -D| P_{\neq}\psi|^2.
\end{aligned}
\end{equation}

Furthermore, substituting $\psi=\nabla(\Psi+q)$ into the vector equation in \eqref{eq:physical-system} shows that the spatial gradient of
\[
(\Psi+q)_t - D\Delta(\Psi+q) + 2DV(\Psi+q)_{x_1} + D|\nabla(\Psi+q)|^2 - \phi
\]
vanishes identically in $\Omega$. At $x_1=0$, we have $\Psi(0,t)= q(0,x',t)=0$, $\psi_2(0,x',t)=\psi_3(0,x',t)=0$, and $\Delta(\Psi+q)|_{x_1=0}=\psi_{1x_1}(0,x',t)$, so the second boundary condition in \eqref{eq:perturbation-bc} forces the entire expression above to vanish at $x_1=0$. Consequently, the scalar logarithmic potential $\Psi+q$ satisfies the exact equation
\begin{equation}\label{eq:nonzero-potential-exact}
(\Psi+q)_t = D\Delta(\Psi+q) - 2DV(\Psi+q)_{x_1} - D|\psi|^2 + \phi \quad\text{in }\Omega.
\end{equation}
Applying the nonzero-mode projection $P_{\neq}$ to \eqref{eq:nonzero-potential-exact} yields the evolution equation for $q= P_{\neq}(\Psi+q)$, which will be coupled with the weighted Neumann compensator in Lemma~\ref{lem:nonzero-neumann-potential-coupling}.

\section{A priori estimates}\label{sec:apriori}
Throughout this section, let $T>0$ and let $(\phi,\psi)\in X(0,T)$ be a smooth compatible solution satisfying the bootstrap assumption
\begin{equation}\label{eq:bootstrap-assumption}
N(T)\le\delta_0\le1,
\end{equation}
where $\delta_0>0$ is a sufficiently small threshold depending only on the fixed parameters $(p,\alpha,D,b,\lambda)$. All integrations by parts on the unbounded domain are justified via spatial cutoffs $\chi_R(x_1)$ as $R\to\infty$ using the weighted integrability in $X(0,T)$.

\subsection{The longitudinal zero mode}
We begin by deriving the weighted $L^2$ and $H^1$ energy estimates for the zero-mode antiderivatives $(\Phi,\Psi)$.

\begin{Lem}\label{lem:zero-first-energy}
Let $t\in[0,T]$. Assume $N(t)\le\delta_0\le1$ with $\delta_0>0$ sufficiently small. Then
\begin{equation}\label{eq:zero-first-energy}
\begin{aligned}
&\int_0^\infty\left(\frac{\Phi^2(t)}{W^pU}+\Psi^2(t)\right)\dd x_1
+\int_0^t\int_0^\infty\left(\Phi^2+U\Psi^2+\frac{\Phi_{x_1}^2}{U}+\Psi_{x_1}^2\right)\dd x_1\dd\tau\\
&\le C\int_0^\infty\left(\frac{\Phi_0^2}{W^pU}+\Psi_0^2\right)\dd x_1
+CN(t)\int_0^t\int_\Omega\frac{|\nabla\phi|^2+|\nabla\psi|^2}{U}\dd x\dd\tau.
\end{aligned}
\end{equation}
\end{Lem}
\begin{proof}
We rewrite the zero-mode system \eqref{eq:zero-mode-system} as
\begin{equation}\label{eq:zero-mode-energy-system}
\left\{
\begin{aligned}
W^{-p}\Phi_t &= e^{-p\Psi}\left(\Phi_{x_1x_1}+(\alpha-1)pV\Phi_{x_1}\right)+(\alpha-1)pU\Psi_{x_1}\\
&\quad +(\alpha-1)pe^{-p\Psi}\Phi_{x_1}\Psi_{x_1}+\left(e^{-p\Psi}-1\right)(\alpha-1)pU\Psi_{x_1}+W^{-p} P_0 R_1,\\
\Psi_t &= D\Psi_{x_1x_1}-2DV\Psi_{x_1}-D\Psi_{x_1}^2+\Phi_{x_1}+ P_0 R_2.
\end{aligned}
\right.
\end{equation}
Multiplying the first equation of \eqref{eq:zero-mode-energy-system} by $\Phi/U$ and the second by $(\alpha-1)p\Psi$, integrating over $\mathbb{R}_+$, and using $(1/U)_{x_1}=(\alpha-1)pV/U$ together with the boundary conditions $\Phi(0,t)=\Psi(0,t)=0$, we notice that the linear cross terms $(\alpha-1)p\Phi\Psi_{x_1}$ and $(\alpha-1)p\Psi\Phi_{x_1}$ integrate to zero:
\[
(\alpha-1)p\int_0^\infty\left(\Phi\Psi_{x_1}+\Psi\Phi_{x_1}\right)\dd x_1 = (\alpha-1)p\int_0^\infty(\Phi\Psi)_{x_1}\dd x_1 = 0.
\]
Hence we obtain the energy identity
\begin{equation}\label{eq:zero-mode-basic-identity}
\begin{aligned}
&\frac{1}{2}\frac{\dd}{\dd t}\int_0^\infty\left(\frac{\Phi^2}{W^pU}+(\alpha-1)p\Psi^2\right)\dd x_1
+\int_0^\infty\frac{\Phi_{x_1}^2}{U}\dd x_1
+(\alpha-1)pD\int_0^\infty\left(\Psi_{x_1}^2-V_{x_1}\Psi^2\right)\dd x_1\\
&= \alpha p\int_0^\infty\frac{e^{-p\Psi}\Phi\Phi_{x_1}\Psi_{x_1}}{U}\dd x_1
-(\alpha-1)pD\int_0^\infty\Psi\Psi_{x_1}^2\dd x_1\\
&\quad +\int_0^\infty\left(e^{-p\Psi}-1\right)\left[(\alpha-1)p\Phi\Psi_{x_1}-\frac{\Phi_{x_1}^2}{U}\right]\dd x_1\\
&\quad +\int_0^\infty\frac{\Phi}{W^pU} P_0 R_1\dd x_1
+(\alpha-1)p\int_0^\infty\Psi P_0 R_2\dd x_1\\
&= \sum_{i=1}^5 A_i.
\end{aligned}
\end{equation}
By the profile Hardy inequality \eqref{eq:profile-hardy-appendix}, we have
\begin{equation}\label{eq:Phi-hardy-coercivity}
\int_0^\infty\Phi^2\dd x_1 \le C\int_0^\infty\frac{\Phi_{x_1}^2}{U}\dd x_1.
\end{equation}
Moreover, using $(1/U)_{x_1}=(\alpha-1)pV/U$ and $V^2\le CU$, we have
\[
\begin{aligned}
\frac{\Psi_{x_1}^2(x_1,t)}{U(x_1)}
&= -\int_{x_1}^\infty\frac{2\Psi_{x_1}\Psi_{x_1x_1}}{U}\dd y_1
-\int_{x_1}^\infty\left(\frac{1}{U}\right)_{x_1}\Psi_{x_1}^2\dd y_1\\
&\le C\int_0^\infty\frac{\Psi_{x_1}^2+\Psi_{x_1x_1}^2}{U}\dd y_1,
\end{aligned}
\]
which together with \eqref{eq:bootstrap-norm} and \eqref{eq:Psi-Linfty-smallness} yields
\begin{equation}\label{eq:Psi-basic-Linfty-bounds}
\left\|\frac{\Psi_{x_1}(\cdot,t)}{\sqrt{U}}\right\|_{L^\infty(\mathbb{R}_+)} \le CN(t), \qquad \|\Psi(\cdot,t)\|_{L^\infty(\mathbb{R}_+)} \le CN(t)\le C,
\end{equation}
and hence $|e^{-p\Psi}|\le C$ and $|e^{-p\Psi}-1|\le C|\Psi|\le CN(t)$.

We now estimate $A_1,\dots,A_5$ in \eqref{eq:zero-mode-basic-identity}. By \eqref{eq:Psi-basic-Linfty-bounds} and Cauchy--Schwarz,
\[
\begin{aligned}
|A_1| &\le C\left\|\frac{\Psi_{x_1}}{\sqrt{U}}\right\|_{L^\infty(\mathbb{R}_+)}\int_0^\infty|\Phi|\frac{|\Phi_{x_1}|}{\sqrt{U}}\dd x_1
\le CN(t)\int_0^\infty\left(\Phi^2+\frac{\Phi_{x_1}^2}{U}\right)\dd x_1,\\
|A_2| &\le C\|\Psi\|_{L^\infty(\mathbb{R}_+)}\int_0^\infty\Psi_{x_1}^2\dd x_1
\le CN(t)\int_0^\infty\Psi_{x_1}^2\dd x_1,\\
|A_3| &\le C\|\Psi\|_{L^\infty(\mathbb{R}_+)}\int_0^\infty\left(|\Phi||\Psi_{x_1}|+\frac{\Phi_{x_1}^2}{U}\right)\dd x_1
\le CN(t)\int_0^\infty\left(\Phi^2+\Psi_{x_1}^2+\frac{\Phi_{x_1}^2}{U}\right)\dd x_1.
\end{aligned}
\]
For $A_4$, using \eqref{eq:zero-mode-remainders} and $|e^{-p q}-1|\le C| q|$ from \eqref{eq:exponential-bounds}, we split
\begin{equation}\label{eq:A4-decomposition}
\begin{aligned}
|A_4| &\le C\int_0^\infty\left(\int_{\mathbb{T}^2}| q\phi_{x_1}|\dd x'\right)\frac{|\Phi|}{U}\dd x_1
+C\int_0^\infty\left(\int_{\mathbb{T}^2}| q V\phi|\dd x'\right)\frac{|\Phi|}{U}\dd x_1\\
&\quad +C\int_0^\infty\left(\int_{\mathbb{T}^2}| q U\psi_1|\dd x'\right)\frac{|\Phi|}{U}\dd x_1
+C\int_0^\infty\left(\int_{\mathbb{T}^2}| q\phi\psi_1|\dd x'\right)\frac{|\Phi|}{U}\dd x_1\\
&\quad +C\int_0^\infty\left(\int_{\mathbb{T}^2}| P_{\neq}\phi\, P_{\neq}\psi_1|\dd x'\right)\frac{|\Phi|}{U}\dd x_1
= \sum_{j=1}^5 I_j.
\end{aligned}
\end{equation}
Using the modal potential bounds \eqref{eq:weighted-potential-estimates-concise}, the reconstruction bound \eqref{eq:mode-reconstruction-phi-appendix}, the transversal Poincar\'e inequality \eqref{eq:weighted-transversal-poincare}, and the Hardy inequality \eqref{eq:Phi-hardy-coercivity}, we estimate $I_1,\dots,I_5$ term by term:
\begin{equation}\label{eq:A41}
\begin{aligned}
I_1 &\le C\left\|\frac{q}{\sqrt{U}}\right\|_{L^\infty_{x_1}L^2_{x'}}\left(\int_\Omega\frac{|\phi_{x_1}|^2}{U}\dd x\right)^{1/2}\|\Phi\|_{L^2(\mathbb{R}_+)}\\
&\le CN(t)\left(\|\nabla\phi\|_{\omega_2}^2+\int_0^\infty\frac{\Phi_{x_1}^2}{U}\dd x_1\right),
\end{aligned}
\end{equation}
\begin{equation}\label{eq:A42}
\begin{aligned}
I_2 &\le C\| q\|_{L^\infty(\Omega)}\left(\int_\Omega\frac{V^2\phi^2}{U^2}\dd x\right)^{1/2}\|\Phi\|_{L^2(\mathbb{R}_+)}\\
&\le CN(t)\left(\|\phi\|_{\omega_2}^2+\|\Phi\|_{L^2(\mathbb{R}_+)}^2\right)
\le CN(t)\left(\int_0^\infty\frac{\Phi_{x_1}^2}{U}\dd x_1+\|\nabla\phi\|_{\omega_2}^2\right),
\end{aligned}
\end{equation}
\begin{equation}\label{eq:A43}
\begin{aligned}
I_3 &\le C\| q\|_{L^\infty(\Omega)}\|\psi_1\|_{L^2(\Omega)}\|\Phi\|_{L^2(\mathbb{R}_+)}\\
&\le CN(t)\left(\int_0^\infty\Psi_{x_1}^2\dd x_1+\| P_{\neq}\psi_1\|_{L^2(\Omega)}^2+\|\Phi\|_{L^2(\mathbb{R}_+)}^2\right)\\
&\le CN(t)\left(\int_0^\infty\Psi_{x_1}^2\dd x_1+\int_0^\infty\frac{\Phi_{x_1}^2}{U}\dd x_1+\|\nabla\psi\|_{\omega_2}^2\right),
\end{aligned}
\end{equation}
\begin{equation}\label{eq:A44}
\begin{aligned}
I_4 &\le C\| q\|_{L^\infty(\Omega)}\left\|\frac{\psi_1}{\sqrt{U}}\right\|_{L^\infty_{x_1}L^2_{x'}}\|\phi\|_{\omega_2}\|\Phi\|_{L^2(\mathbb{R}_+)}\\
&\le CN(t)\left(\int_0^\infty\frac{\Phi_{x_1}^2}{U}\dd x_1+\|\nabla\phi\|_{\omega_2}^2\right),
\end{aligned}
\end{equation}
and
\begin{equation}\label{eq:A45}
\begin{aligned}
I_5 &\le C\left\|\frac{P_{\neq}\psi_1}{\sqrt{U}}\right\|_{L^\infty_{x_1}L^2_{x'}}\| P_{\neq}\phi\|_{\omega_2}\|\Phi\|_{L^2(\mathbb{R}_+)}\\
&\le CN(t)\left(\int_0^\infty\frac{\Phi_{x_1}^2}{U}\dd x_1+\|\nabla\phi\|_{\omega_2}^2\right).
\end{aligned}
\end{equation}
Summing \eqref{eq:A41}--\eqref{eq:A45} gives
\begin{equation}\label{eq:first-zero-N21-final}
|A_4| \le CN(t)\left(\int_0^\infty\frac{\Phi_{x_1}^2}{U}\dd x_1+\int_0^\infty\Psi_{x_1}^2\dd x_1+\|\nabla\phi\|_{\omega_2}^2+\|\nabla\psi\|_{\omega_2}^2\right).
\end{equation}
For $A_5$, using $R_2=-D| P_{\neq}\psi|^2$, $\|\Psi\|_{L^\infty(\mathbb{R}_+)}\le CN(t)$, and the transversal Poincar\'e inequality \eqref{eq:weighted-transversal-poincare}, we obtain
\begin{equation}\label{eq:first-zero-N22-final}
\begin{aligned}
|A_5| &\le C\|\Psi\|_{L^\infty(\mathbb{R}_+)}\| P_{\neq}\psi\|_{L^2(\Omega)}^2
\le C\|\Psi\|_{L^\infty(\mathbb{R}_+)}\|\nabla_{x'}\psi\|_{L^2(\Omega)}^2
\le CN(t)\|\nabla\psi\|_{\omega_2}^2.
\end{aligned}
\end{equation}
Substituting the bounds for $A_1,\dots,A_5$ into \eqref{eq:zero-mode-basic-identity} and using $-V_{x_1}=\frac{(\alpha-1)p}{D[(\alpha-1)p+2]}U$ together with \eqref{eq:Phi-hardy-coercivity}, we arrive at
\begin{equation}\label{eq:first-zero-before-absorption}
\begin{aligned}
&\frac{\dd}{\dd t}\int_0^\infty\left(\frac{\Phi^2}{W^pU}+(\alpha-1)p\Psi^2\right)\dd x_1
+c\int_0^\infty\left(\Phi^2+U\Psi^2+\frac{\Phi_{x_1}^2}{U}+\Psi_{x_1}^2\right)\dd x_1\\
&\le CN(t)\int_0^\infty\left(\Phi^2+\frac{\Phi_{x_1}^2}{U}+\Psi_{x_1}^2\right)\dd x_1
+CN(t)\left(\|\nabla\phi\|_{\omega_2}^2+\|\nabla\psi\|_{\omega_2}^2\right).
\end{aligned}
\end{equation}
Choosing $\delta_0>0$ sufficiently small so that $CN(t)\le C\delta_0\le c/2$ and integrating over $[0,t]$ yields \eqref{eq:zero-first-energy}.
\end{proof}

\begin{Lem}\label{lem:zero-second-energy}
Let $t\in[0,T]$. Assume $N(t)\le\delta_0\le1$ with $\delta_0>0$ sufficiently small. Then
\begin{equation}\label{eq:zero-second-energy}
\begin{aligned}
&\int_0^\infty\left(\frac{\Phi_{x_1}^2(t)}{W^pU}+\frac{\Psi_{x_1}^2(t)}{U}\right)\dd x_1
+\int_0^t\int_0^\infty\frac{\Phi_{x_1x_1}^2+\Psi_{x_1x_1}^2}{U}\dd x_1\dd\tau\\
&\le C\int_0^\infty\left(\frac{\Phi_{0x_1}^2}{W^pU}+\frac{\Psi_{0x_1}^2}{U}\right)\dd x_1
+C\int_0^t\int_0^\infty\left(\frac{\Phi_{x_1}^2}{U}+\Psi_{x_1}^2\right)\dd x_1\dd\tau\\
&\quad +CN(t)\int_0^t\int_\Omega\frac{|\nabla\phi|^2+|\nabla\psi|^2}{U}\dd x\dd\tau.
\end{aligned}
\end{equation}
\end{Lem}
\begin{proof}
Multiplying the first equation of \eqref{eq:zero-mode-energy-system} by $-\Phi_{x_1x_1}/U$, integrating over $\mathbb{R}_+$, and integrating the time-derivative term by parts in $x_1$, we obtain
\begin{equation}\label{eq:second-zero-B-identity}
\begin{aligned}
&\frac{1}{2}\frac{\dd}{\dd t}\int_0^\infty\frac{\Phi_{x_1}^2}{W^pU}\dd x_1+\int_0^\infty\frac{\Phi_{x_1x_1}^2}{U}\dd x_1\\
&= \int_0^\infty\big(1-e^{-p\Psi}\big)\frac{\Phi_{x_1x_1}^2}{U}\dd x_1
-(\alpha-1)p\int_0^\infty e^{-p\Psi}\Phi_{x_1x_1}\Psi_{x_1}\dd x_1\\
&\quad -(\alpha-1)p\int_0^\infty\frac{V}{U}e^{-p\Psi}\Phi_{x_1}\Phi_{x_1x_1}\dd x_1
-(\alpha-1)p\int_0^\infty\frac{e^{-p\Psi}\Phi_{x_1}\Phi_{x_1x_1}\Psi_{x_1}}{U}\dd x_1\\
&\quad -\int_0^\infty\left(\frac{1}{W^pU}\right)_{x_1}W^pe^{-p\Psi}\Big[\Phi_{x_1x_1}+(\alpha-1)p\big(V\Phi_{x_1}+U\Psi_{x_1}+\Phi_{x_1}\Psi_{x_1}\big)\Big]\Phi_{x_1}\dd x_1\\
&\quad -\int_0^\infty\left[\frac{\Phi_{x_1x_1}}{W^pU}+\left(\frac{1}{W^pU}\right)_{x_1}\Phi_{x_1}\right] P_0 R_1\dd x_1
+\int_0^\infty\left(\frac{\Phi_t\Phi_{x_1}}{W^pU}\right)_{x_1}\dd x_1\\
&= \sum_{i=1}^7 B_i.
\end{aligned}
\end{equation}
Using $|1-e^{-p\Psi}|\le CN(t)$, $\|\Psi_{x_1}\|_{L^\infty(\mathbb{R}_+)}\le CN(t)$, $V^2/U\le C$, and the identity
\[
\left(\frac{1}{W^pU}\right)_{x_1}W^p = \frac{\alpha pV}{U},
\]
we apply Young's inequality to estimate $B_1,\dots,B_5$:
\begin{equation}\label{eq:second-zero-B1}
|B_1| \le CN(t)\int_0^\infty\frac{\Phi_{x_1x_1}^2}{U}\dd x_1,
\end{equation}
\begin{equation}\label{eq:second-zero-B2}
|B_2| \le \frac{1}{16}\int_0^\infty\frac{\Phi_{x_1x_1}^2}{U}\dd x_1+C\int_0^\infty\Psi_{x_1}^2\dd x_1,
\end{equation}
\begin{equation}\label{eq:second-zero-B3}
|B_3| \le \frac{1}{16}\int_0^\infty\frac{\Phi_{x_1x_1}^2}{U}\dd x_1+C\int_0^\infty\frac{\Phi_{x_1}^2}{U}\dd x_1,
\end{equation}
\begin{equation}\label{eq:second-zero-B4}
|B_4| \le CN(t)\int_0^\infty\left(\frac{\Phi_{x_1}^2}{U}+\frac{\Phi_{x_1x_1}^2}{U}\right)\dd x_1,
\end{equation}
and
\begin{equation}\label{eq:second-zero-B5}
\begin{aligned}
|B_5| &\le C\int_0^\infty\frac{|V|}{U}\Big(|\Phi_{x_1x_1}|+|V\Phi_{x_1}|+|U\Psi_{x_1}|+|\Phi_{x_1}\Psi_{x_1}|\Big)|\Phi_{x_1}|\dd x_1\\
&\le \frac{1}{16}\int_0^\infty\frac{\Phi_{x_1x_1}^2}{U}\dd x_1+C\int_0^\infty\left(\frac{\Phi_{x_1}^2}{U}+\Psi_{x_1}^2\right)\dd x_1+CN(t)\int_0^\infty\frac{\Phi_{x_1}^2}{U}\dd x_1.
\end{aligned}
\end{equation}
For $B_6$, since $|(1/(W^pU))_{x_1}W^p|=\alpha p|V|/U\le C/U$, we have
\[
|B_6| \le C\int_0^\infty\frac{|\Phi_{x_1x_1}|+|\Phi_{x_1}|}{U}\int_{\mathbb{T}^2}\Big(| q\phi_{x_1}|+| q V\phi|+| q U\psi_1|+| q\phi\psi_1|+| P_{\neq}\phi\, P_{\neq}\psi_1|\Big)\dd x'\dd x_1.
\]
Estimating the five integrand terms against $(|\Phi_{x_1x_1}|+|\Phi_{x_1}|)/U$ in the exact same way as $I_1,\dots,I_5$ in \eqref{eq:A41}--\eqref{eq:A45} (with $(|\Phi_{x_1x_1}|+|\Phi_{x_1}|)/\sqrt{U}$ in $L^2(\mathbb{R}_+)$ in place of $|\Phi|/\sqrt{U}$), we obtain
\begin{equation}\label{eq:second-zero-B6}
|B_6| \le CN(t)\left(\int_0^\infty\left(\frac{\Phi_{x_1}^2+\Phi_{x_1x_1}^2}{U}+\Psi_{x_1}^2\right)\dd x_1+\|\nabla\phi\|_{\omega_2}^2+\|\nabla\psi\|_{\omega_2}^2\right).
\end{equation}
For $B_7$, $\Phi(0,t)=0$ implies $\Phi_t(0,t)=0$, and the far-field contribution vanishes by the cutoff argument in Remark~\ref{rem:far-field-zero-mode}, so
\begin{equation}\label{eq:second-zero-B7}
B_7 = 0.
\end{equation}
Substituting \eqref{eq:second-zero-B1}--\eqref{eq:second-zero-B7} into \eqref{eq:second-zero-B-identity}, absorbing the $CN(t)\|\Phi_{x_1x_1}\|_{\omega_2}^2$ terms for $\delta_0$ sufficiently small, and integrating over $[0,t]$ yields
\begin{equation}\label{eq:second-zero-integrated}
\begin{aligned}
&\int_0^\infty\frac{\Phi_{x_1}^2(t)}{W^pU}\dd x_1+\int_0^t\int_0^\infty\frac{\Phi_{x_1x_1}^2}{U}\dd x_1\dd\tau\\
&\le C\int_0^\infty\frac{\Phi_{0x_1}^2}{W^pU}\dd x_1+C\int_0^t\int_0^\infty\left(\frac{\Phi_{x_1}^2}{U}+\Psi_{x_1}^2\right)\dd x_1\dd\tau
+CN(t)\int_0^t\left(\|\nabla\phi\|_{\omega_2}^2+\|\nabla\psi\|_{\omega_2}^2\right)\dd\tau.
\end{aligned}
\end{equation}

Next, multiplying the second equation of \eqref{eq:zero-mode-energy-system} by $-\Psi_{x_1x_1}/U$ and using the identity
\[
-\Psi_t\frac{\Psi_{x_1x_1}}{U} = \frac{1}{2}\left(\frac{\Psi_{x_1}^2}{U}\right)_t + \left(\frac{1}{U}\right)_{x_1}\Psi_{x_1}\Psi_t - \left(\frac{\Psi_{x_1}\Psi_t}{U}\right)_{x_1},
\]
we substitute the equation for $\Psi_t$ into the middle term and integrate by parts in the $P_0 R_2$ term to obtain
\begin{equation}\label{eq:second-zero-psi-identity}
\begin{aligned}
&\frac{1}{2}\frac{\dd}{\dd t}\int_0^\infty\frac{\Psi_{x_1}^2}{U}\dd x_1+D\int_0^\infty\frac{\Psi_{x_1x_1}^2}{U}\dd x_1\\
&= -\int_0^\infty\frac{\Phi_{x_1}\Psi_{x_1x_1}}{U}\dd x_1
+2D\int_0^\infty\frac{V}{U}\Psi_{x_1}\Psi_{x_1x_1}\dd x_1\\
&\quad -\int_0^\infty\left(\frac{1}{U}\right)_{x_1}\left(D\Psi_{x_1x_1}-2DV\Psi_{x_1}-D\Psi_{x_1}^2+\Phi_{x_1}\right)\Psi_{x_1}\dd x_1\\
&\quad +D\int_0^\infty\frac{\Psi_{x_1}^2\Psi_{x_1x_1}}{U}\dd x_1
-\int_0^\infty\left(\frac{\Psi_{x_1}}{U}\right)_{x_1} P_0 R_2\dd x_1
+\int_0^\infty\left(\frac{\Psi_{x_1}\Psi_t}{U}\right)_{x_1}\dd x_1\\
&= \sum_{i=8}^{13} B_i.
\end{aligned}
\end{equation}
Since $|(1/U)_{x_1}|=(\alpha-1)p|V|/U\le C/\sqrt{U}$ and $\|\Psi_{x_1}/\sqrt{U}\|_{L^\infty(\mathbb{R}_+)}\le CN(t)$ by \eqref{eq:Psi-x-weighted-Linfty}, Young's inequality gives
\begin{equation}\label{eq:second-zero-B8}
|B_8| \le \frac{D}{16}\int_0^\infty\frac{\Psi_{x_1x_1}^2}{U}\dd x_1+C\int_0^\infty\frac{\Phi_{x_1}^2}{U}\dd x_1,
\end{equation}
\begin{equation}\label{eq:second-zero-B9}
|B_9| \le \frac{D}{16}\int_0^\infty\frac{\Psi_{x_1x_1}^2}{U}\dd x_1+C\int_0^\infty\Psi_{x_1}^2\dd x_1,
\end{equation}
\begin{equation}\label{eq:second-zero-B10}
\begin{aligned}
|B_{10}| &\le \frac{D}{16}\int_0^\infty\frac{\Psi_{x_1x_1}^2}{U}\dd x_1+C\int_0^\infty\left(\frac{\Phi_{x_1}^2}{U}+\Psi_{x_1}^2+\frac{\Psi_{x_1}^4}{U}\right)\dd x_1\\
&\le \frac{D}{16}\int_0^\infty\frac{\Psi_{x_1x_1}^2}{U}\dd x_1+C\int_0^\infty\left(\frac{\Phi_{x_1}^2}{U}+\Psi_{x_1}^2\right)\dd x_1,
\end{aligned}
\end{equation}
and
\begin{equation}\label{eq:second-zero-B11}
|B_{11}| \le D\left\|\frac{\Psi_{x_1}}{\sqrt{U}}\right\|_{L^\infty(\mathbb{R}_+)}\int_0^\infty\frac{|\Psi_{x_1}\Psi_{x_1x_1}|}{\sqrt{U}}\dd x_1
\le CN(t)\int_0^\infty\left(\Psi_{x_1}^2+\frac{\Psi_{x_1x_1}^2}{U}\right)\dd x_1.
\end{equation}
For $B_{12}$, using $|(\Psi_{x_1}/U)_{x_1}|\le |\Psi_{x_1x_1}|/U+C|\Psi_{x_1}|/\sqrt{U}$ and $\| P_{\neq}\psi/\sqrt{U}\|_{L^\infty(\Omega)}\le CN(t)$ from \eqref{eq:weighted-vector-Linfty} and \eqref{eq:Psi-x-weighted-Linfty}, we have
\[
 P_0(| P_{\neq}\psi|^2) \le \left\|\frac{P_{\neq}\psi}{\sqrt{U}}\right\|_{L^\infty(\Omega)}\sqrt{U}\,\| P_{\neq}\psi\|_{L^2(\mathbb{T}^2)}
\le CN(t)\sqrt{U}\,\| P_{\neq}\psi\|_{L^2(\mathbb{T}^2)},
\]
and hence
\begin{equation}\label{eq:second-zero-B12}
\begin{aligned}
|B_{12}| &\le CN(t)\int_0^\infty\left(\frac{|\Psi_{x_1x_1}|}{\sqrt{U}}+|\Psi_{x_1}|\right)\| P_{\neq}\psi\|_{L^2(\mathbb{T}^2)}\dd x_1\\
&\le CN(t)\int_0^\infty\left(\frac{\Psi_{x_1x_1}^2}{U}+\Psi_{x_1}^2\right)\dd x_1+CN(t)\|\nabla\psi\|_{\omega_2}^2.
\end{aligned}
\end{equation}
Finally, $\Psi_t(0,t)=0$ and the far-field cutoff argument give
\begin{equation}\label{eq:second-zero-B13}
B_{13}=0.
\end{equation}
Substituting \eqref{eq:second-zero-B8}--\eqref{eq:second-zero-B13} into \eqref{eq:second-zero-psi-identity}, absorbing $CN(t)\|\Psi_{x_1x_1}\|_{\omega_2}^2$ for $\delta_0$ sufficiently small, and integrating over $[0,t]$ gives
\begin{equation}\label{eq:second-zero-psi-integrated}
\begin{aligned}
&\int_0^\infty\frac{\Psi_{x_1}^2(t)}{U}\dd x_1+\int_0^t\int_0^\infty\frac{\Psi_{x_1x_1}^2}{U}\dd x_1\dd\tau\\
&\le C\int_0^\infty\frac{\Psi_{0x_1}^2}{U}\dd x_1+C\int_0^t\int_0^\infty\left(\frac{\Phi_{x_1}^2}{U}+\Psi_{x_1}^2\right)\dd x_1\dd\tau
+CN(t)\int_0^t\|\nabla\psi\|_{\omega_2}^2\dd\tau.
\end{aligned}
\end{equation}
Adding \eqref{eq:second-zero-integrated} and \eqref{eq:second-zero-psi-integrated} yields \eqref{eq:zero-second-energy}.
\end{proof}

\begin{Rem}\label{rem:far-field-zero-mode}
The preceding product estimates and \eqref{eq:zero-mode-system} imply that
\[
\begin{aligned}
\|\Phi_t\|_{\omega_1}^2 &\le C\left(\|\Phi_{x_1x_1}\|_{\omega_2}^2+\|\Phi_{x_1}\|_{\omega_2}^2+\|\Psi_{x_1}\|_{L^2(\mathbb{R}_+)}^2+N^2(t)\big(\|\nabla\phi\|_{\omega_2}^2+\|\phi\|_{\omega_2}^2+\|\psi\|_{L^2(\Omega)}^2\big)\right),\\
\|\Psi_t\|_{\omega_2}^2 &\le C\left(\|\Psi_{x_1x_1}\|_{\omega_2}^2+\|\Psi_{x_1}\|_{L^2(\mathbb{R}_+)}^2+\|\Phi_{x_1}\|_{\omega_2}^2+N^2(t)\| P_{\neq}\psi\|_{L^2(\Omega)}^2\right).
\end{aligned}
\]
Hence $\omega_1\Phi_t\Phi_{x_1}$ and $\omega_2\Psi_t\Psi_{x_1}$ belong to $L^1((0,t)\times\mathbb{R}_+)$, so multiplying by $\chi_R(x_1)$ and letting $R\to\infty$ rigorously justifies $B_7=B_{13}=0$.
\end{Rem}

\begin{Cor}\label{cor:combined-zero-mode}
Assume $N(t)\le\delta_0$ with $\delta_0>0$ sufficiently small. Then
\begin{equation}\label{eq:combined-zero-mode}
\begin{aligned}
&\int_0^\infty\left(\frac{\Phi^2(t)+\Phi_{x_1}^2(t)}{W^pU}+\Psi^2(t)+\frac{\Psi_{x_1}^2(t)}{U}\right)\dd x_1\\
&\quad +\int_0^t\int_0^\infty\left(\Phi^2+U\Psi^2+\frac{\Phi_{x_1}^2+\Phi_{x_1x_1}^2+\Psi_{x_1x_1}^2}{U}+\Psi_{x_1}^2\right)\dd x_1\dd\tau\\
&\le C\int_0^\infty\left(\frac{\Phi_0^2+\Phi_{0x_1}^2}{W^pU}+\Psi_0^2+\frac{\Psi_{0x_1}^2}{U}\right)\dd x_1
+CN(t)\int_0^t\int_\Omega\frac{|\nabla\phi|^2+|\nabla\psi|^2}{U}\dd x\dd\tau.
\end{aligned}
\end{equation}
\end{Cor}
\begin{proof}
Substituting the dissipation bound for $\int_0^t\int_0^\infty(\Phi_{x_1}^2/U+\Psi_{x_1}^2)\dd x_1\dd\tau$ from Lemma~\ref{lem:zero-first-energy} into the right-hand side of Lemma~\ref{lem:zero-second-energy} and adding the two inequalities gives \eqref{eq:combined-zero-mode}.
\end{proof}

\subsection{Transversal nonzero modes}
To control the transversal nonzero modes without imposing any smallness restriction on the steady mass $\lambda>0$, we couple the weighted Neumann compensator on the nonzero modes with the scalar potential equation for $q= P_{\neq}(\Psi+q)$.

\begin{Lem}\label{lem:nonzero-neumann-potential-coupling}
Let $t\in[0,T]$. Assume $N(t)\le\delta_0\le1$ with $\delta_0>0$ sufficiently small. Then
\begin{equation}\label{eq:nonzero-neumann-potential-coupled}
\begin{aligned}
&\int_\Omega | q(t)|^2\dd x
+\int_0^t\int_\Omega\left(\frac{| P_{\neq}\phi|^2}{U}+| P_{\neq}\psi|^2+U| q|^2\right)\dd x\dd\tau\\
&\le C\int_\Omega\left(\frac{|\phi_0|^2}{W^pU}+\frac{|\psi_0|^2}{U}\right)\dd x
+CN(t)\int_0^t\int_0^\infty\left(\frac{\Phi_{x_1}^2}{U}+\Psi_{x_1}^2\right)\dd x_1\dd\tau.
\end{aligned}
\end{equation}
\end{Lem}
\begin{proof}
By Lemma~\ref{lem:potential-representation}, we have $\nabla(\Psi+q)=\psi$ and $P_{\neq}(\Psi+q)= q$. Using $U_{x_1}+(\alpha-1)pUV=0$ from \eqref{eq:steady-profile-identities}, we rewrite the density equation in \eqref{eq:physical-system} in the symmetric divergence form
\begin{equation}\label{eq:nonzero-density-flux}
\phi_\tau = \nabla\cdot\left\{W^pe^{-p(\Psi+q)}\left[U\nabla\left(\frac{\phi}{U}\right)+(\alpha-1)p(U+\phi)\nabla(\Psi+q)\right]\right\}.
\end{equation}
For each $\tau\in[0,t]$, Lemma~\ref{lem:weighted-nonzero-neumann} applied to $P_{\neq}\phi(\cdot,\tau)$ yields a unique weak solution $z(\cdot,\tau)$ of the weighted Neumann problem
\begin{equation}\label{eq:nonzero-neumann-z}
-(\alpha-1)p\,\nabla\cdot(W^pU\nabla z)= P_{\neq}\phi \quad\text{in }\Omega,
\end{equation}
subject to the boundary and zero-mean conditions
\begin{equation}\label{eq:nonzero-neumann-z-bc}
z_{x_1}(0,x',\tau)=0, \qquad  P_0 z=0.
\end{equation}
The weak formulation of \eqref{eq:nonzero-neumann-z}--\eqref{eq:nonzero-neumann-z-bc} reads
\begin{equation}\label{eq:nonzero-neumann-weak}
(\alpha-1)p\int_\Omega W^pU\nabla z\cdot\nabla h\dd x = \int_\Omega  P_{\neq}\phi\,h\dd x
\end{equation}
for every zero-transversal-mean test function $h$ with $\int_\Omega W^pU|\nabla h|^2\dd x<\infty$. Since $U=U(x_1)$ depends only on $x_1$, we have $P_{\neq}(\phi/U)=(P_{\neq}\phi)/U$. Taking $h= P_{\neq}(\phi/U)$ in \eqref{eq:nonzero-neumann-weak} produces the coercive density dissipation
\begin{equation}\label{eq:nonzero-density-coercivity}
\int_\Omega W^pU\nabla z\cdot\nabla P_{\neq}\left(\frac{\phi}{U}\right)\dd x = \frac{1}{(\alpha-1)p}\| P_{\neq}\phi\|_{\omega_2}^2.
\end{equation}
Moreover, by the self-adjointness of the time-independent Neumann inverse (Lemma~\ref{lem:weighted-nonzero-neumann} and Remark~\ref{rem:nonzero-tests-justification}),
\begin{equation}\label{eq:nonzero-neumann-time}
\int_\Omega\phi_\tau z\dd x = \int_\Omega P_{\neq}\phi_\tau z\dd x = \frac{1}{2}\frac{\dd}{\dd\tau}\int_\Omega P_{\neq}\phi\,z\dd x.
\end{equation}
Testing \eqref{eq:nonzero-density-flux} against $z$ and integrating by parts over $\Omega$, we note that the normal flux in \eqref{eq:nonzero-density-flux} vanishes at $x_1=0$ by \eqref{eq:perturbation-bc}:
\[
\left[U\partial_{x_1}\left(\frac{\phi}{U}\right)+(\alpha-1)p(U+\phi)(\Psi+q)_{x_1}\right]_{x_1=0}
= \big[\phi_{x_1}+(\alpha-1)p(V\phi+U\psi_1+\phi\psi_1)\big]_{x_1=0} = 0.
\]
Hence
\begin{equation}\label{eq:nonzero-density-tested}
\begin{aligned}
\int_\Omega\phi_\tau z\dd x
&= -\int_\Omega W^pe^{-p(\Psi+q)}U\nabla\left(\frac{\phi}{U}\right)\cdot\nabla z\dd x\\
&\quad -(\alpha-1)p\int_\Omega W^pe^{-p(\Psi+q)}(U+\phi)\nabla(\Psi+q)\cdot\nabla z\dd x.
\end{aligned}
\end{equation}
Since $P_0 z=0$ and $W^pU$ depends only on $x_1$, transversal orthogonality and \eqref{eq:nonzero-density-coercivity} give
\begin{equation}\label{eq:nonzero-density-main}
-\int_\Omega W^pU\nabla\left(\frac{\phi}{U}\right)\cdot\nabla z\dd x
= -\int_\Omega W^pU\nabla P_{\neq}\left(\frac{\phi}{U}\right)\cdot\nabla z\dd x
= -\frac{1}{(\alpha-1)p}\| P_{\neq}\phi\|_{\omega_2}^2.
\end{equation}
Similarly, taking $h= q= P_{\neq}(\Psi+q)$ in \eqref{eq:nonzero-neumann-weak} yields
\begin{equation}\label{eq:nonzero-linear-coupling}
-(\alpha-1)p\int_\Omega W^pU\nabla(\Psi+q)\cdot\nabla z\dd x
= -(\alpha-1)p\int_\Omega W^pU\nabla q\cdot\nabla z\dd x
= -\int_\Omega P_{\neq}\phi\, q\dd x.
\end{equation}
Combining \eqref{eq:nonzero-neumann-time}--\eqref{eq:nonzero-linear-coupling}, we arrive at the Neumann compensator identity
\begin{equation}\label{eq:nonzero-neumann-dynamic}
\begin{aligned}
&\frac{1}{2}\frac{\dd}{\dd\tau}\int_\Omega P_{\neq}\phi\,z\dd x+\frac{1}{(\alpha-1)p}\| P_{\neq}\phi\|_{\omega_2}^2\\
&= -\int_\Omega P_{\neq}\phi\, q\dd x
-\int_\Omega W^pU\big(e^{-p(\Psi+q)}-1\big)\nabla\left(\frac{\phi}{U}\right)\cdot\nabla z\dd x\\
&\quad -(\alpha-1)p\int_\Omega W^pU\big(e^{-p(\Psi+q)}-1\big)\psi\cdot\nabla z\dd x\\
&\quad -(\alpha-1)p\int_\Omega W^pe^{-p(\Psi+q)}\phi\psi\cdot\nabla z\dd x.
\end{aligned}
\end{equation}

To cancel the linear coupling term $-\int_\Omega P_{\neq}\phi\, q\dd x$ on the right-hand side of \eqref{eq:nonzero-neumann-dynamic}, we apply $P_{\neq}$ to the scalar potential equation \eqref{eq:nonzero-potential-exact} to get
\[
 q_\tau = D\Delta q - 2DV q_{x_1} - D P_{\neq}(|\psi|^2) +  P_{\neq}\phi.
\]
Multiplying this equation by $q$ and integrating by parts over $\Omega$, using $q(0,x',\tau)=0$, $\nabla q= P_{\neq}\psi$, and $-DV_{x_1}=\frac{(\alpha-1)p}{(\alpha-1)p+2}U$ from \eqref{eq:steady-profile-identities}, we obtain
\begin{equation}\label{eq:nonzero-potential-energy}
\begin{aligned}
&\frac{1}{2}\frac{\dd}{\dd\tau}\| q\|_{L^2(\Omega)}^2+D\| P_{\neq}\psi\|_{L^2(\Omega)}^2+\frac{(\alpha-1)p}{(\alpha-1)p+2}\|\sqrt{U}\, q\|_{L^2(\Omega)}^2\\
&= \int_\Omega P_{\neq}\phi\, q\dd x - D\int_\Omega|\psi|^2 q\dd x.
\end{aligned}
\end{equation}
Adding \eqref{eq:nonzero-neumann-dynamic} and \eqref{eq:nonzero-potential-energy}, the linear cross terms $\mp\int_\Omega P_{\neq}\phi\, q\dd x$ cancel identically, leaving
\begin{equation}\label{eq:nonzero-coupled-energy-identity}
\begin{aligned}
&\frac{1}{2}\frac{\dd}{\dd\tau}\left(\int_\Omega P_{\neq}\phi\,z\dd x+\| q\|_{L^2(\Omega)}^2\right)\\
&\quad +\frac{1}{(\alpha-1)p}\| P_{\neq}\phi\|_{\omega_2}^2+D\| P_{\neq}\psi\|_{L^2(\Omega)}^2+\frac{(\alpha-1)p}{(\alpha-1)p+2}\|\sqrt{U}\, q\|_{L^2(\Omega)}^2\\
&= -\int_\Omega W^pU\big(e^{-p(\Psi+q)}-1\big)\nabla\left(\frac{\phi}{U}\right)\cdot\nabla z\dd x\\
&\quad -(\alpha-1)p\int_\Omega W^pU\big(e^{-p(\Psi+q)}-1\big)\psi\cdot\nabla z\dd x\\
&\quad -(\alpha-1)p\int_\Omega W^pe^{-p(\Psi+q)}\phi\psi\cdot\nabla z\dd x
-D\int_\Omega|\psi|^2 q\dd x\\
&= \sum_{i=1}^4 E_i.
\end{aligned}
\end{equation}
By \eqref{eq:weighted-potential-estimates-concise}--\eqref{eq:Psi-Linfty-smallness}, we have the $L^\infty$ bounds
\begin{equation}\label{eq:nonzero-pointwise}
\|\Psi+q\|_{L^\infty(\Omega)} + \|\phi\|_{L^\infty(\Omega)} + \left\|\frac{\psi}{\sqrt{U}}\right\|_{L^\infty(\Omega)} \le CN(t),
\end{equation}
and consequently
\begin{equation}\label{eq:nonzero-exponential}
\left\|e^{-p(\Psi+q)}-1\right\|_{L^\infty(\Omega)} \le CN(t), \qquad \left|\nabla e^{-p(\Psi+q)}\right| \le C|\psi| \le CN(t)\sqrt{U}.
\end{equation}
For $E_1$, we rewrite the integrand by the product rule and apply the weak formulation \eqref{eq:nonzero-neumann-weak} with $h= P_{\neq}[(e^{-p(\Psi+q)}-1)\phi/U]$:
\[
\begin{aligned}
E_1 &= -\int_\Omega W^pU\nabla z\cdot\nabla\left[\big(e^{-p(\Psi+q)}-1\big)\frac{\phi}{U}\right]\dd x
+\int_\Omega W^p\phi\,\nabla e^{-p(\Psi+q)}\cdot\nabla z\dd x\\
&= -\frac{1}{(\alpha-1)p}\int_\Omega P_{\neq}\phi\,\big(e^{-p(\Psi+q)}-1\big)\frac{\phi}{U}\dd x
+\int_\Omega W^p\phi\,\nabla e^{-p(\Psi+q)}\cdot\nabla z\dd x.
\end{aligned}
\]
Using \eqref{eq:nonzero-exponential} and the Neumann flux estimate $\|W^pU\nabla z\|_{L^2(\Omega)}\le C\| P_{\neq}\phi\|_{\omega_2}$ from \eqref{eq:L46-neumann-flux}, we obtain
\begin{equation}\label{eq:nonzero-nonlinear-one}
\begin{aligned}
|E_1| &\le CN(t)\| P_{\neq}\phi\|_{\omega_2}\|\phi\|_{\omega_2} + CN(t)\|\phi\|_{\omega_2}\|W^pU\nabla z\|_{L^2(\Omega)}\\
&\le CN(t)\left(\|\phi\|_{\omega_2}^2+\| P_{\neq}\phi\|_{\omega_2}^2\right).
\end{aligned}
\end{equation}
Similarly, applying \eqref{eq:nonzero-pointwise}, \eqref{eq:nonzero-exponential}, and \eqref{eq:L46-neumann-flux} to $E_2$ and $E_3$ yields
\begin{equation}\label{eq:nonzero-nonlinear-two-three}
\begin{aligned}
|E_2| &\le CN(t)\|\psi\|_{L^2(\Omega)}\|W^pU\nabla z\|_{L^2(\Omega)} \le CN(t)\left(\|\psi\|_{L^2(\Omega)}^2+\| P_{\neq}\phi\|_{\omega_2}^2\right),\\
|E_3| &\le C\left\|\frac{\psi}{\sqrt{U}}\right\|_{L^\infty(\Omega)}\|\phi\|_{\omega_2}\|W^pU\nabla z\|_{L^2(\Omega)} \le CN(t)\|\phi\|_{\omega_2}^2.
\end{aligned}
\end{equation}
For $E_4$, the transversal Poincar\'e inequality gives $\| q\|_{L^2(\Omega)}\le C\|\nabla_{x'} q\|_{L^2(\Omega)}\le C\|\psi\|_{L^2(\Omega)}$, so that
\begin{equation}\label{eq:nonzero-potential-nonlinear}
|E_4| \le D\|\psi\|_{L^\infty(\Omega)}\|\psi\|_{L^2(\Omega)}\| q\|_{L^2(\Omega)} \le CN(t)\|\psi\|_{L^2(\Omega)}^2.
\end{equation}
Combining \eqref{eq:nonzero-nonlinear-one}--\eqref{eq:nonzero-potential-nonlinear} with the orthogonal decompositions
\begin{equation}\label{eq:nonzero-orthogonal-splitting}
\|\phi\|_{\omega_2}^2 = \| P_{\neq}\phi\|_{\omega_2}^2+\int_0^\infty\frac{\Phi_{x_1}^2}{U}\dd x_1, \qquad
\|\psi\|_{L^2(\Omega)}^2 = \| P_{\neq}\psi\|_{L^2(\Omega)}^2+\int_0^\infty\Psi_{x_1}^2\dd x_1,
\end{equation}
and absorbing $CN(t)(\| P_{\neq}\phi\|_{\omega_2}^2+\| P_{\neq}\psi\|_{L^2(\Omega)}^2)$ into the left-hand side of \eqref{eq:nonzero-coupled-energy-identity} for $\delta_0$ sufficiently small, we obtain
\begin{equation}\label{eq:nonzero-closed-differential}
\begin{aligned}
&\frac{\dd}{\dd\tau}\left(\frac{1}{2}\int_\Omega P_{\neq}\phi\,z\dd x+\frac{1}{2}\| q\|_{L^2(\Omega)}^2\right)
+\frac{1}{2(\alpha-1)p}\| P_{\neq}\phi\|_{\omega_2}^2\\
&\quad +\frac{D}{2}\| P_{\neq}\psi\|_{L^2(\Omega)}^2+\frac{(\alpha-1)p}{(\alpha-1)p+2}\|\sqrt{U}\, q\|_{L^2(\Omega)}^2\\
&\le CN(t)\int_0^\infty\left(\frac{\Phi_{x_1}^2}{U}+\Psi_{x_1}^2\right)\dd x_1.
\end{aligned}
\end{equation}
By \eqref{eq:L46-neumann-energy} and the transversal Poincar\'e inequality, we have
\[
0 \le \int_\Omega P_{\neq}\phi\,z\dd x \le C\|\phi\|_{\omega_1}^2, \qquad \| q(0)\|_{L^2(\Omega)}^2 \le C\|\nabla_{x'} q(0)\|_{L^2(\Omega)}^2 \le C\|\psi_0\|_{\omega_2}^2.
\]
Integrating \eqref{eq:nonzero-closed-differential} over $[0,t]$ yields \eqref{eq:nonzero-neumann-potential-coupled}.
\end{proof}

\begin{Rem}\label{rem:nonzero-tests-justification}
Let $H_{\neq}^1$ be the completion of smooth zero-transversal-mean functions on $\Omega$ under the norm $\|h\|_{H_{\neq}^1}^2=(\alpha-1)p\int_\Omega W^pU|\nabla h|^2\dd x$. By \eqref{eq:weighted-transversal-poincare}, $\int_\Omega W^pU|h|^2\dd x\le C\|h\|_{H_{\neq}^1}^2$. Since
\[
\begin{aligned}
\int_\Omega W^pU\left|\nabla\left(\frac{\phi}{U}\right)\right|^2\dd x &\le C\left(\|\nabla\phi\|_{\omega_2}^2+\|\phi\|_{\omega_2}^2\right),\\
\int_\Omega W^pU\left|\frac{\phi}{U}\nabla e^{-p(\Psi+q)}\right|^2\dd x &\le CN^2(t)\|\phi\|_{\omega_2}^2,
\end{aligned}
\]
both $P_{\neq}(\phi/U)$ and $P_{\neq}[(e^{-p(\Psi+q)}-1)\phi/U]$ belong to $H_{\neq}^1$, as does $q$. Moreover, the density equation and zero normal flux at $x_1=0$ imply $P_{\neq}\phi_t\in L^2(0,t;(H_{\neq}^1)^*)$, so $z\in H^1(0,t;H_{\neq}^1)$ and the differentiation formula \eqref{eq:nonzero-neumann-time} holds in the $(H_{\neq}^1)^*$--$H_{\neq}^1$ duality pairing for almost every $\tau\in(0,t)$.
\end{Rem}

\begin{Cor}\label{cor:combined-zero-nonzero-mode}
Assume $N(t)\le\delta_0\le1$ with $\delta_0>0$ sufficiently small. Then
\begin{equation}\label{eq:combined-zero-nonzero-mode}
\begin{aligned}
&\int_0^\infty\left(\frac{\Phi^2(t)+\Phi_{x_1}^2(t)}{W^pU}+\Psi^2(t)+\frac{\Psi_{x_1}^2(t)}{U}\right)\dd x_1+\int_\Omega | q(t)|^2\dd x\\
&\quad +\int_0^t\int_0^\infty\left(\Phi^2+U\Psi^2+\frac{\Phi_{x_1x_1}^2+\Psi_{x_1x_1}^2}{U}\right)\dd x_1\dd\tau\\
&\quad +\int_0^t\int_\Omega\left(\frac{|\phi|^2}{U}+|\psi|^2+U| q|^2\right)\dd x\dd\tau\\
&\le C\left[\int_0^\infty\left(\frac{\Phi_0^2}{W^pU}+\Psi_0^2\right)\dd x_1+\int_\Omega\left(\frac{|\phi_0|^2}{W^pU}+\frac{|\psi_0|^2}{U}\right)\dd x\right]\\
&\quad +CN(t)\int_0^t\int_\Omega\frac{|\nabla\phi|^2+|\nabla\psi|^2}{U}\dd x\dd\tau.
\end{aligned}
\end{equation}
\end{Cor}
\begin{proof}
By Jensen's inequality, $\|\Phi_{0x_1}\|_{\omega_1}^2\le\|\phi_0\|_{\omega_1}^2$ and $\|\Psi_{0x_1}\|_{\omega_2}^2\le\|\psi_0\|_{\omega_2}^2$. Substituting the zero-mode dissipation bound from Corollary~\ref{cor:combined-zero-mode} into the right-hand side of Lemma~\ref{lem:nonzero-neumann-potential-coupling} and adding \eqref{eq:combined-zero-mode} together with the orthogonal decompositions \eqref{eq:nonzero-orthogonal-splitting} yields \eqref{eq:combined-zero-nonzero-mode}.
\end{proof}

\subsection{Physical lower-order energies and gradient recovery}
Having controlled the zero-order physical dissipation $\|\phi\|_{\omega_2}^2+\|\psi\|_{L^2(\Omega)}^2$ in Corollary~\ref{cor:combined-zero-nonzero-mode}, we now establish the direct physical energy estimate for $(\phi, P_{\neq}\psi)$ and reconstruct the full gradient dissipation $\|\nabla\psi\|_{\omega_2}^2$.

\begin{Lem}\label{lem:direct-physical-potential-time}
Let $t\in[0,T]$. Assume $N(t)\le\delta_0\le1$ with $\delta_0>0$ sufficiently small. Then
\begin{equation}\label{eq:direct-physical-potential-time}
\begin{aligned}
&\int_\Omega\left(\frac{|\phi(t)|^2}{W^pU}+\frac{| P_{\neq}\psi(t)|^2}{U}\right)\dd x
+\int_0^t\int_\Omega\left(\frac{|\nabla\phi|^2}{U}+\frac{| q_\tau|^2}{U}\right)\dd x\dd\tau\\
&\quad +\int_0^t\int_\Omega\frac{\left|\nabla\phi+(\alpha-1)pV\phi e_1+(\alpha-1)p(U+\phi)\psi\right|^2}{U}\dd x\dd\tau\\
&\le C\int_\Omega\left(\frac{|\phi_0|^2}{W^pU}+\frac{| P_{\neq}\psi_0|^2}{U}\right)\dd x
+C\int_0^t\int_\Omega\left(\frac{| P_{\neq}\phi|^2}{U}+| P_{\neq}\psi|^2\right)\dd x\dd\tau\\
&\quad +C\int_0^t\int_0^\infty\left(\frac{\Phi_{x_1}^2}{U}+\Psi_{x_1}^2\right)\dd x_1\dd\tau.
\end{aligned}
\end{equation}
\end{Lem}
\begin{proof}
\emph{Step 1.} By \eqref{eq:nonzero-pointwise}, we may choose $\delta_0>0$ small enough that
\begin{equation}\label{eq:direct-physical-coefficients}
\frac{1}{2}\le e^{-p(\Psi+q)}\le 2, \qquad \left|\frac{V}{\sqrt{U}}\right|\le C, \qquad U\le C.
\end{equation}
Multiplying the first equation of \eqref{eq:physical-system} by $\phi/(W^pU)$, integrating over $\Omega$, and using
\[
\nabla\left(\frac{\phi}{W^pU}\right) = \frac{\nabla\phi+\alpha pV\phi e_1}{W^pU}
\]
together with the vanishing of the normal flux at $x_1=0$ from \eqref{eq:perturbation-bc}, we obtain
\begin{equation}\label{eq:direct-physical-density-identity}
\begin{aligned}
&\frac{1}{2}\frac{\dd}{\dd\tau}\|\phi\|_{\omega_1}^2
+\int_\Omega\frac{e^{-p(\Psi+q)}}{U}|\nabla\phi|^2\dd x
+\alpha(\alpha-1)p^2\int_\Omega\frac{e^{-p(\Psi+q)}}{U}V^2\phi^2\dd x\\
&= -(2\alpha-1)p\int_\Omega\frac{e^{-p(\Psi+q)}V\phi\phi_{x_1}}{U}\dd x
-(\alpha-1)p\int_\Omega e^{-p(\Psi+q)}\psi\cdot\nabla\phi\dd x\\
&\quad -(\alpha-1)p\int_\Omega\frac{e^{-p(\Psi+q)}\phi\psi\cdot\nabla\phi}{U}\dd x
-\alpha(\alpha-1)p^2\int_\Omega\frac{e^{-p(\Psi+q)}V\phi(U+\phi)\psi_1}{U}\dd x\\
&= \sum_{i=1}^4 M_i.
\end{aligned}
\end{equation}
By \eqref{eq:nonzero-pointwise}, \eqref{eq:direct-physical-coefficients}, and Young's inequality, for any $\eta>0$,
\[
\begin{aligned}
|M_1| &\le \eta\|\nabla\phi\|_{\omega_2}^2+C_\eta\int_\Omega\frac{V^2\phi^2}{U}\dd x \le \eta\|\nabla\phi\|_{\omega_2}^2+C_\eta\|\phi\|_{L^2(\Omega)}^2,\\
|M_2| &\le \eta\|\nabla\phi\|_{\omega_2}^2+C_\eta\|\sqrt{U}\,\psi\|_{L^2(\Omega)}^2,\\
|M_3| &\le C\left\|\frac{\psi}{\sqrt{U}}\right\|_{L^\infty(\Omega)}\|\phi\|_{L^2(\Omega)}\|\nabla\phi\|_{\omega_2} \le \eta\|\nabla\phi\|_{\omega_2}^2+C_\eta N^2(t)\|\phi\|_{L^2(\Omega)}^2,\\
|M_4| &\le C\|\phi\|_{L^2(\Omega)}\|\sqrt{U}\,\psi\|_{L^2(\Omega)}+C\left\|\frac{\psi}{\sqrt{U}}\right\|_{L^\infty(\Omega)}\|\phi\|_{L^2(\Omega)}^2 \le C\left(\|\phi\|_{L^2(\Omega)}^2+\|\sqrt{U}\,\psi\|_{L^2(\Omega)}^2\right).
\end{aligned}
\]
Taking $\eta>0$ sufficiently small and using the orthogonal decompositions \eqref{eq:nonzero-orthogonal-splitting} along with $U\le C$, we integrate \eqref{eq:direct-physical-density-identity} over $[0,t]$ to find
\begin{equation}\label{eq:direct-physical-density-integrated}
\begin{aligned}
&\|\phi(t)\|_{\omega_1}^2+\int_0^t\|\nabla\phi\|_{\omega_2}^2\dd\tau\\
&\le C\|\phi_0\|_{\omega_1}^2+C\int_0^t\left(\| P_{\neq}\phi\|_{\omega_2}^2+\| P_{\neq}\psi\|_{L^2(\Omega)}^2\right)\dd\tau
+C\int_0^t\int_0^\infty\left(\frac{\Phi_{x_1}^2}{U}+\Psi_{x_1}^2\right)\dd x_1\dd\tau.
\end{aligned}
\end{equation}

\emph{Step 2.} Applying $P_{\neq}$ to \eqref{eq:nonzero-potential-exact}, multiplying by $q_\tau/U$, and integrating by parts over $\Omega$ using $q_\tau(0,x',\tau)=0$, $\nabla q= P_{\neq}\psi$, and $(1/U)_{x_1}=(\alpha-1)pV/U$, we obtain
\begin{equation}\label{eq:direct-physical-potential-time-identity}
\begin{aligned}
&\frac{D}{2}\frac{\dd}{\dd\tau}\| P_{\neq}\psi\|_{\omega_2}^2+\| q_\tau\|_{\omega_2}^2\\
&= -D[(\alpha-1)p+2]\int_\Omega\frac{V q_{x_1} q_\tau}{U}\dd x
-D\int_\Omega\frac{|\psi|^2 q_\tau}{U}\dd x
+\int_\Omega\frac{P_{\neq}\phi\, q_\tau}{U}\dd x\\
&= M_5+M_6+M_7.
\end{aligned}
\end{equation}
By Young's inequality, $V^2/U\le C$, and $\|\psi/\sqrt{U}\|_{L^\infty(\Omega)}\le CN(t)$, we estimate
\[
\begin{aligned}
|M_5| &\le \frac{1}{6}\| q_\tau\|_{\omega_2}^2+C\| P_{\neq}\psi\|_{L^2(\Omega)}^2,\\
|M_6| &\le C\left\|\frac{\psi}{\sqrt{U}}\right\|_{L^\infty(\Omega)}\|\psi\|_{L^2(\Omega)}\| q_\tau\|_{\omega_2}
\le \frac{1}{6}\| q_\tau\|_{\omega_2}^2+CN^2(t)\|\psi\|_{L^2(\Omega)}^2,\\
|M_7| &\le \frac{1}{6}\| q_\tau\|_{\omega_2}^2+C\| P_{\neq}\phi\|_{\omega_2}^2.
\end{aligned}
\]
Substituting these bounds into \eqref{eq:direct-physical-potential-time-identity}, using \eqref{eq:nonzero-orthogonal-splitting}, and integrating over $[0,t]$ gives
\begin{equation}\label{eq:direct-physical-potential-time-integrated}
\begin{aligned}
&\| P_{\neq}\psi(t)\|_{\omega_2}^2+\int_0^t\| q_\tau\|_{\omega_2}^2\dd\tau\\
&\le C\| P_{\neq}\psi_0\|_{\omega_2}^2+C\int_0^t\left(\| P_{\neq}\phi\|_{\omega_2}^2+\| P_{\neq}\psi\|_{L^2(\Omega)}^2\right)\dd\tau
+CN^2(t)\int_0^t\int_0^\infty\Psi_{x_1}^2\dd x_1\dd\tau.
\end{aligned}
\end{equation}

\emph{Step 3.} Finally, by \eqref{eq:nonzero-pointwise} and \eqref{eq:nonzero-orthogonal-splitting}, the full physical flux satisfies
\begin{equation}\label{eq:direct-physical-flux-bound}
\begin{aligned}
&\int_\Omega\frac{\left|\nabla\phi+(\alpha-1)pV\phi e_1+(\alpha-1)p(U+\phi)\psi\right|^2}{U}\dd x\\
&\le C\left(\|\nabla\phi\|_{\omega_2}^2+\|\phi\|_{L^2(\Omega)}^2+\|\sqrt{U}\,\psi\|_{L^2(\Omega)}^2+\left\|\frac{\psi}{\sqrt{U}}\right\|_{L^\infty(\Omega)}^2\|\phi\|_{L^2(\Omega)}^2\right)\\
&\le C\left(\|\nabla\phi\|_{\omega_2}^2+\| P_{\neq}\phi\|_{\omega_2}^2+\| P_{\neq}\psi\|_{L^2(\Omega)}^2+\int_0^\infty\left(\frac{\Phi_{x_1}^2}{U}+\Psi_{x_1}^2\right)\dd x_1\right).
\end{aligned}
\end{equation}
Combining \eqref{eq:direct-physical-density-integrated}, \eqref{eq:direct-physical-potential-time-integrated}, and the time integral of \eqref{eq:direct-physical-flux-bound} yields \eqref{eq:direct-physical-potential-time}.
\end{proof}

\begin{Cor}\label{cor:physical-and-mode-estimates}
Assume $N(t)\le\delta_0\le1$ with $\delta_0>0$ sufficiently small. Then
\begin{equation}\label{eq:physical-and-mode-estimates}
\begin{aligned}
&\int_0^\infty\left(\frac{\Phi^2(t)}{W^pU}+\Psi^2(t)\right)\dd x_1
+\int_\Omega\left(\frac{|\phi(t)|^2}{W^pU}+\frac{|\psi(t)|^2}{U}+| q(t)|^2\right)\dd x\\
&\quad +\int_0^t\int_0^\infty\left(\Phi^2+U\Psi^2+\frac{\Phi_{x_1x_1}^2+\Psi_{x_1x_1}^2}{U}\right)\dd x_1\dd\tau\\
&\quad +\int_0^t\int_\Omega\left(\frac{|\phi|^2+|\nabla\phi|^2+| q_\tau|^2}{U}+|\psi|^2+U| q|^2\right)\dd x\dd\tau\\
&\quad +\int_0^t\int_\Omega\frac{\left|\nabla\phi+(\alpha-1)pV\phi e_1+(\alpha-1)p(U+\phi)\psi\right|^2}{U}\dd x\dd\tau\\
&\le C\left[\int_0^\infty\left(\frac{\Phi_0^2}{W^pU}+\Psi_0^2\right)\dd x_1+\int_\Omega\left(\frac{|\phi_0|^2}{W^pU}+\frac{|\psi_0|^2}{U}\right)\dd x\right]
+CN(t)\int_0^t\int_\Omega\frac{|\nabla\psi|^2}{U}\dd x\dd\tau.
\end{aligned}
\end{equation}
\end{Cor}
\begin{proof}
Bounding the right-hand side of \eqref{eq:direct-physical-potential-time} by Corollary~\ref{cor:combined-zero-nonzero-mode}, using the orthogonal identity $\|\psi(t)\|_{\omega_2}^2=\| P_{\neq}\psi(t)\|_{\omega_2}^2+\|\Psi_{x_1}(t)\|_{\omega_2}^2$, and absorbing $CN(t)\int_0^t\|\nabla\phi\|_{\omega_2}^2\dd\tau$ into the left-hand side for $\delta_0$ sufficiently small gives \eqref{eq:physical-and-mode-estimates}.
\end{proof}

\begin{Lem}\label{lem:all-mass-gradient-reconstruction}
Let $t\in[0,T]$ and suppose $(\phi,\psi)\in X(0,T)$ satisfies $N(t)\le\delta_0\le1$. Then for every $\tau\in[0,t]$,
\begin{equation}\label{eq:psi-gradient-reconstruction}
\begin{aligned}
\int_\Omega\frac{|\nabla\psi(\tau)|^2}{U}\dd x
&\le C\int_\Omega\left(\frac{| q_\tau(\tau)|^2}{U}+\frac{| P_{\neq}\phi(\tau)|^2}{U}+| P_{\neq}\psi(\tau)|^2\right)\dd x\\
&\quad +C\int_0^\infty\left(\frac{\Phi_{x_1}^2(\tau)+\Psi_{x_1x_1}^2(\tau)}{U}+\Psi_{x_1}^2(\tau)\right)\dd x_1.
\end{aligned}
\end{equation}
\end{Lem}
\begin{proof}
By \eqref{eq:potential-identity}, $\nabla\psi=\nabla^2(\Psi+q)$. We first prove the weighted elliptic identity
\begin{equation}\label{eq:L46-weighted-elliptic}
\|\nabla^2(\Psi+q)\|_{\omega_2}^2 \le C\left(\|\Delta(\Psi+q)\|_{\omega_2}^2+\|\nabla_{x'}(\Psi+q)\|_{L^2(\Omega)}^2\right).
\end{equation}
Expand $\Psi+q$ in the transversal Fourier basis on $\mathbb{T}^2$:
\[
(\Psi+q)(x_1,x',\tau) = \sum_{k\in 2\pi\mathbb{Z}^2} f_k(x_1,\tau)e^{ik\cdot x'}, \qquad f_0=\Psi.
\]
Since $(\Psi+q)|_{x_1=0}=0$, we have $f_k(0,\tau)=0$ for all $k$. Moreover, the explicit formula \eqref{eq:explicit-steady-profile} shows that $\omega_2(x_1)=1/U(x_1)$ is a quadratic polynomial in $x_1$ with constant second derivative
\[
\omega_2'' = \frac{(\alpha-1)^2p^2}{D[(\alpha-1)p+2]} > 0.
\]
Integrating by parts twice in $x_1$ and using $f_k(0,\tau)=0$, we obtain
\[
-2|k|^2\operatorname{Re}\int_0^\infty\omega_2 f_k''\overline{f_k}\dd x_1
= 2|k|^2\int_0^\infty\omega_2|f_k'|^2\dd x_1 - |k|^2\omega_2''\int_0^\infty|f_k|^2\dd x_1,
\]
and hence
\begin{equation}\label{eq:gradient-fourier-identity}
\int_0^\infty\omega_2\left|f_k''-|k|^2f_k\right|^2\dd x_1
= \int_0^\infty\omega_2\left(|f_k''|^2+2|k|^2|f_k'|^2+|k|^4|f_k|^2\right)\dd x_1 - |k|^2\omega_2''\int_0^\infty|f_k|^2\dd x_1.
\end{equation}
Summing \eqref{eq:gradient-fourier-identity} over $k\in 2\pi\mathbb{Z}^2$ and applying Plancherel's theorem~\cite{Grafakos2014} yields
\[
\sum_{i,j=1}^3\|\partial_{x_i}\partial_{x_j}(\Psi+q)\|_{\omega_2}^2
= \|\Delta(\Psi+q)\|_{\omega_2}^2 + \omega_2''\|\nabla_{x'}(\Psi+q)\|_{L^2(\Omega)}^2,
\]
which proves \eqref{eq:L46-weighted-elliptic}.

Next, by the scalar potential equation \eqref{eq:nonzero-potential-exact},
\begin{equation}\label{eq:gradient-potential-laplacian}
\Delta(\Psi+q) = D^{-1}(\Psi+q)_\tau + 2V\psi_1 + |\psi|^2 - D^{-1}\phi.
\end{equation}
Using $V^2/U\le C$, $\nabla_{x'}(\Psi+q)=(\psi_2,\psi_3)$, and
\[
\||\psi|^2\|_{\omega_2} \le \left\|\frac{\psi}{\sqrt{U}}\right\|_{L^\infty(\Omega)}\|\psi\|_{L^2(\Omega)} \le CN(t)\|\psi\|_{L^2(\Omega)},
\]
we deduce from \eqref{eq:L46-weighted-elliptic} that
\begin{equation}\label{eq:gradient-potential-time-reduction}
\|\nabla\psi\|_{\omega_2}^2 \le C\left(\|(\Psi+q)_\tau\|_{\omega_2}^2+\|\phi\|_{\omega_2}^2+\|\psi\|_{L^2(\Omega)}^2\right).
\end{equation}
By the second equation of \eqref{eq:zero-mode-system} and Jensen's inequality,
\begin{equation}\label{eq:gradient-zero-time-intermediate}
\|\Psi_\tau\|_{\omega_2}^2 \le C\left(\|\Psi_{x_1x_1}\|_{\omega_2}^2+\|\Psi_{x_1}\|_{L^2(\mathbb{R}_+)}^2+\|\Phi_{x_1}\|_{\omega_2}^2+N^2(t)\|\psi\|_{L^2(\Omega)}^2\right).
\end{equation}
Combining \eqref{eq:gradient-potential-time-reduction} and \eqref{eq:gradient-zero-time-intermediate} with the orthogonal splittings $\|(\Psi+q)_\tau\|_{\omega_2}^2=\| q_\tau\|_{\omega_2}^2+\|\Psi_\tau\|_{\omega_2}^2$ and \eqref{eq:nonzero-orthogonal-splitting} gives \eqref{eq:psi-gradient-reconstruction}.
\end{proof}

\begin{Cor}\label{cor:closed-lower-order-energy}
Assume $N(t)\le\delta_0$ with $\delta_0>0$ sufficiently small. Then there exists a constant $C>0$, independent of $T$, such that for every $t\in[0,T]$,
\begin{equation}\label{eq:closed-lower-order-energy}
\begin{aligned}
&\int_0^\infty\left(\frac{\Phi^2(t)+\Phi_{x_1}^2(t)}{W^pU}+\Psi^2(t)+\frac{\Psi_{x_1}^2(t)}{U}\right)\dd x_1\\
&\quad +\int_\Omega\left(\frac{|\phi(t)|^2}{W^pU}+\frac{|\psi(t)|^2}{U}+|\Psi(t)+q(t)|^2\right)\dd x\\
&\quad +\int_0^t\int_0^\infty\left(\Phi^2+U\Psi^2+\frac{\Phi_{x_1}^2+\Phi_{x_1x_1}^2+\Psi_{x_1x_1}^2}{U}+\Psi_{x_1}^2\right)\dd x_1\dd\tau\\
&\quad +\int_0^t\int_\Omega\left(\frac{|\phi|^2+|\nabla\phi|^2+|\nabla\psi|^2}{U}+|\psi|^2+U|\Psi+q|^2\right)\dd x\dd\tau\\
&\le C\left[\int_0^\infty\left(\frac{\Phi_0^2+\Phi_{0x_1}^2}{W^pU}+\Psi_0^2+\frac{\Psi_{0x_1}^2}{U}\right)\dd x_1+\int_\Omega\left(\frac{|\phi_0|^2}{W^pU}+\frac{|\psi_0|^2}{U}\right)\dd x\right].
\end{aligned}
\end{equation}
\end{Cor}
\begin{proof}
Integrating \eqref{eq:psi-gradient-reconstruction} over $[0,t]$ shows that $\int_0^t\|\nabla\psi\|_{\omega_2}^2\dd\tau$ is bounded by a constant multiple of the dissipation on the left-hand side of \eqref{eq:physical-and-mode-estimates}. Substituting this bound into the last term $CN(t)\int_0^t\|\nabla\psi\|_{\omega_2}^2\dd\tau$ of \eqref{eq:physical-and-mode-estimates} and taking $\delta_0>0$ sufficiently small absorbs the remainder and yields \eqref{eq:closed-lower-order-energy}.
\end{proof}

\subsection{First-order tangential derivatives}
We now derive the weighted energy estimates for the tangential derivatives $\partial_{x'}^\beta\phi$ and $\partial_{x'}^\beta\psi$ ($1\le|\beta|\le2$). Write
\[
F = \nabla\phi+(\alpha-1)pV\phi e_1+(\alpha-1)p(U+\phi)\psi.
\]
Since $\nabla(\Psi+q)=\psi$ and $\partial_{x'}^\beta(\Psi+q)$ has zero transversal mean for $|\beta|\ge1$, the weighted transversal Poincar\'e inequality \eqref{eq:weighted-transversal-poincare} gives directly
\begin{equation}\label{eq:tangential-potential-algebraic}
\begin{aligned}
\|\partial_{x'}^\beta(\Psi+q)\|_{L^2(\Omega)}^2
&\le C\|\nabla_{x'}\partial_{x'}^\beta(\Psi+q)\|_{L^2(\Omega)}^2
\le C\|\partial_{x'}^\beta\psi\|_{\omega_2}^2,\\
\|\partial_{x'}^\beta\psi\|_{L^2(\Omega)}^2
&\le C\|\nabla\partial_{x'}^\beta\psi\|_{\omega_2}^2,
\qquad 1\le|\beta|\le2.
\end{aligned}
\end{equation}

\begin{Lem}\label{lem:tangential-first-order-estimates}
Let $t\in[0,T]$. Assume $N(t)\le\delta_0\le1$ with $\delta_0>0$ sufficiently small. Then
\begin{equation}\label{eq:tangential-first-order-combined}
\begin{aligned}
&\sum_{|\beta|=1}\left(\|\partial_{x'}^\beta\phi(t)\|_{\omega_1}^2+\|\partial_{x'}^\beta\psi(t)\|_{\omega_2}^2+\|\partial_{x'}^\beta(\Psi+q)(t)\|_{L^2(\Omega)}^2\right)\\
&\quad +\int_0^t\sum_{|\beta|=1}\left(\|\nabla\partial_{x'}^\beta\phi(\tau)\|_{\omega_2}^2+\|\nabla\partial_{x'}^\beta\psi(\tau)\|_{\omega_2}^2+\|\partial_{x'}^\beta\psi(\tau)\|_{L^2(\Omega)}^2\right)\dd\tau\\
&\le C\sum_{|\beta|=1}\left(\|\partial_{x'}^\beta\phi_0\|_{\omega_1}^2+\|\partial_{x'}^\beta\psi_0\|_{\omega_2}^2\right)
+C\int_0^t\left(\|\nabla\phi(\tau)\|_{\omega_2}^2+\|\nabla\psi(\tau)\|_{\omega_2}^2\right)\dd\tau.
\end{aligned}
\end{equation}
\end{Lem}
\begin{proof}
We establish the general tangential energy identities for $1\le|\beta|\le2$, which will also be used in Subsection~\ref{subsec:second-tangential}. Define the tangential commutator remainder $R_\beta$ by
\begin{equation}\label{eq:tangential-flux-remainder-definition}
\partial_{x'}^\beta\left(W^pe^{-p(\Psi+q)}F\right)
= W^pe^{-p(\Psi+q)}\left(\nabla\partial_{x'}^\beta\phi+(\alpha-1)pV\partial_{x'}^\beta\phi\,e_1+(\alpha-1)pU\partial_{x'}^\beta\psi+R_\beta\right).
\end{equation}
At the boundary $x_1=0$, the normal flux $W^pe^{-p(\Psi+q)}(e_1\cdot F)$ vanishes identically by \eqref{eq:perturbation-bc}, so its tangential derivative $\partial_{x'}^\beta(W^pe^{-p(\Psi+q)}(e_1\cdot F))$ also vanishes at $x_1=0$. Multiplying the tangentially differentiated density equation by $\partial_{x'}^\beta\phi/(W^pU)$, integrating by parts over $\Omega$, and using $(1/(W^pU))_{x_1}=\alpha pV/(W^pU)$, we obtain
\begin{equation}\label{eq:tangential-density-identity-unified}
\begin{aligned}
&\frac{1}{2}\frac{\dd}{\dd\tau}\|\partial_{x'}^\beta\phi\|_{\omega_1}^2
+\int_\Omega\frac{e^{-p(\Psi+q)}}{U}|\nabla\partial_{x'}^\beta\phi|^2\dd x
+\alpha(\alpha-1)p^2\int_\Omega\frac{e^{-p(\Psi+q)}}{U}V^2|\partial_{x'}^\beta\phi|^2\dd x\\
&= -(2\alpha-1)p\int_\Omega\frac{e^{-p(\Psi+q)}}{U}V\partial_{x'}^\beta\phi\,\partial_{x_1}\partial_{x'}^\beta\phi\dd x
-(\alpha-1)p\int_\Omega e^{-p(\Psi+q)}\partial_{x'}^\beta\psi\cdot\nabla\partial_{x'}^\beta\phi\dd x\\
&\quad -\alpha(\alpha-1)p^2\int_\Omega e^{-p(\Psi+q)}V\partial_{x'}^\beta\phi\,\partial_{x'}^\beta\psi_1\dd x\\
&\quad -\int_\Omega\frac{e^{-p(\Psi+q)}}{U}R_\beta\cdot\left(\nabla\partial_{x'}^\beta\phi+\alpha pV\partial_{x'}^\beta\phi\,e_1\right)\dd x.
\end{aligned}
\end{equation}
For the vector equation in \eqref{eq:physical-system}, differentiating the boundary conditions \eqref{eq:perturbation-bc} tangentially at $x_1=0$ yields
\[
D\partial_{x_1}\partial_{x'}^\beta\psi_1 = 2DV\partial_{x'}^\beta\psi_1+D\partial_{x'}^\beta(\psi_1^2)-\partial_{x'}^\beta\phi, \qquad \partial_{x'}^\beta\psi_2=\partial_{x'}^\beta\psi_3=0 \quad (x_1=0).
\]
When we test the tangentially differentiated vector equation against $\partial_{x'}^\beta\psi/U$ and integrate by parts over $\Omega$, the boundary contributions from the diffusion, drift, nonlinear, and density gradient terms combine at $x_1=0$ into
\[
\frac{1}{U(0)}\int_{\mathbb{T}^2}\partial_{x'}^\beta\psi_1\left(-D\partial_{x_1}\partial_{x'}^\beta\psi_1+2DV\partial_{x'}^\beta\psi_1+D\partial_{x'}^\beta(\psi_1^2)-\partial_{x'}^\beta\phi\right)\dd x' = 0.
\]
Using $(1/U)_{x_1}=(\alpha-1)pV/U$, we therefore obtain the boundary-free identity
\begin{equation}\label{eq:tangential-vector-identity-unified}
\begin{aligned}
&\frac{1}{2}\frac{\dd}{\dd\tau}\|\partial_{x'}^\beta\psi\|_{\omega_2}^2+D\|\nabla\partial_{x'}^\beta\psi\|_{\omega_2}^2\\
&= -D(\alpha-1)p\int_\Omega\frac{V}{U}\partial_{x_1}\partial_{x'}^\beta\psi\cdot\partial_{x'}^\beta\psi\dd x
+2D\int_\Omega V\partial_{x'}^\beta\psi_1\,\nabla\cdot\left(\frac{\partial_{x'}^\beta\psi}{U}\right)\dd x\\
&\quad +D\int_\Omega\partial_{x'}^\beta(|\psi|^2)\,\nabla\cdot\left(\frac{\partial_{x'}^\beta\psi}{U}\right)\dd x
-\int_\Omega\partial_{x'}^\beta\phi\,\nabla\cdot\left(\frac{\partial_{x'}^\beta\psi}{U}\right)\dd x.
\end{aligned}
\end{equation}
By the weighted transversal Poincar\'e inequality \eqref{eq:weighted-transversal-poincare} and $V^2\le CU$, the testing factors in \eqref{eq:tangential-density-identity-unified} and \eqref{eq:tangential-vector-identity-unified} satisfy
\begin{equation}\label{eq:tangential-testing-factors-unified}
\begin{aligned}
\left\|\frac{\nabla\partial_{x'}^\beta\phi+\alpha pV\partial_{x'}^\beta\phi\,e_1}{\sqrt{U}}\right\|_{L^2(\Omega)}
&\le C\|\nabla\partial_{x'}^\beta\phi\|_{\omega_2},\\
\left\|\sqrt{U}\,\nabla\cdot\left(\frac{\partial_{x'}^\beta\psi}{U}\right)\right\|_{L^2(\Omega)}
&\le C\left(\|\nabla\partial_{x'}^\beta\psi\|_{\omega_2}+\|\partial_{x'}^\beta\psi\|_{\omega_2}\right)
\le C\|\nabla\partial_{x'}^\beta\psi\|_{\omega_2}.
\end{aligned}
\end{equation}
Consequently, by Young's inequality, all linear terms on the right-hand sides of \eqref{eq:tangential-density-identity-unified} and \eqref{eq:tangential-vector-identity-unified} are bounded for any $\eta>0$ by
\begin{equation}\label{eq:tangential-linear-bound-unified}
\eta\left(\|\nabla\partial_{x'}^\beta\phi\|_{\omega_2}^2+\|\nabla\partial_{x'}^\beta\psi\|_{\omega_2}^2\right)
+C_\eta\left(\|\partial_{x'}^\beta\phi\|_{\omega_2}^2+\|\partial_{x'}^\beta\psi\|_{\omega_2}^2\right).
\end{equation}

Now specialize to $|\beta|=1$, so $\partial_{x'}^\beta=\partial_{x_j}$ with $j\in\{2,3\}$. Since $\partial_{x_j}(W^pe^{-p(\Psi+q)})=-p\psi_j W^pe^{-p(\Psi+q)}$, the commutator remainder in \eqref{eq:tangential-flux-remainder-definition} is
\begin{equation}\label{eq:tangential-first-flux-remainder}
R_\beta = -p\psi_j F + (\alpha-1)p\big((\partial_{x_j}\phi)\psi+\phi\partial_{x_j}\psi\big).
\end{equation}
Using the $L^\infty$ bounds \eqref{eq:nonzero-pointwise} and the transversal Poincar\'e inequality $\|\psi_j\|_{L^2(\Omega)}\le C\|\nabla\psi\|_{\omega_2}$, we estimate each component of $R_\beta/\sqrt{U}$ in $L^2(\Omega)$:
\[
\begin{aligned}
\left\|\frac{\psi_j\nabla\phi}{\sqrt{U}}\right\|_{L^2(\Omega)}
&\le \|\psi_j\|_{L^\infty(\Omega)}\|\nabla\phi\|_{\omega_2} \le CN(t)\|\nabla\phi\|_{\omega_2},\\
\left\|\frac{\psi_j\big(V\phi e_1+(U+\phi)\psi\big)}{\sqrt{U}}\right\|_{L^2(\Omega)}
&\le C\left(\|\phi\|_{L^\infty(\Omega)}+\sqrt{U}+\|\phi\|_{L^\infty(\Omega)}\left\|\frac{\psi}{\sqrt{U}}\right\|_{L^\infty(\Omega)}\right)\|\psi_j\|_{L^2(\Omega)}\\
&\le CN(t)\|\nabla\psi\|_{\omega_2},\\
\left\|\frac{(\partial_{x_j}\phi)\psi}{\sqrt{U}}\right\|_{L^2(\Omega)}
&\le \left\|\frac{\psi}{\sqrt{U}}\right\|_{L^\infty(\Omega)}\|\partial_{x_j}\phi\|_{L^2(\Omega)} \le CN(t)\|\nabla\phi\|_{\omega_2},\\
\left\|\frac{\phi\partial_{x_j}\psi}{\sqrt{U}}\right\|_{L^2(\Omega)}
&\le \|\phi\|_{L^\infty(\Omega)}\|\partial_{x_j}\psi\|_{\omega_2} \le CN(t)\|\nabla\psi\|_{\omega_2},
\end{aligned}
\]
which yields
\begin{equation}\label{eq:tangential-first-remainder-bound}
\left\|\frac{R_\beta}{\sqrt{U}}\right\|_{L^2(\Omega)} \le CN(t)\left(\|\nabla\phi\|_{\omega_2}+\|\nabla\psi\|_{\omega_2}\right).
\end{equation}
Similarly, for the nonlinear term in \eqref{eq:tangential-vector-identity-unified},
\[
\left\|\frac{\partial_{x_j}(|\psi|^2)}{\sqrt{U}}\right\|_{L^2(\Omega)}
\le 2\|\psi\|_{L^\infty(\Omega)}\|\partial_{x_j}\psi\|_{\omega_2}
\le CN(t)\|\nabla\partial_{x_j}\psi\|_{\omega_2}.
\]
Summing \eqref{eq:tangential-density-identity-unified} and \eqref{eq:tangential-vector-identity-unified} over $j=2,3$, using $e^{-p(\Psi+q)}\ge 1/2$, choosing $\eta>0$ sufficiently small in \eqref{eq:tangential-linear-bound-unified}, and absorbing the $CN(t)$ dissipation terms for $\delta_0$ sufficiently small, we obtain
\begin{equation}\label{eq:tangential-first-level-red}
\begin{aligned}
&\frac{\dd}{\dd\tau}\sum_{|\beta|=1}\left(\|\partial_{x'}^\beta\phi\|_{\omega_1}^2+\|\partial_{x'}^\beta\psi\|_{\omega_2}^2\right)
+c\sum_{|\beta|=1}\left(\|\nabla\partial_{x'}^\beta\phi\|_{\omega_2}^2+\|\nabla\partial_{x'}^\beta\psi\|_{\omega_2}^2+\|\partial_{x'}^\beta\psi\|_{L^2(\Omega)}^2\right)\\
&\le C\left(\|\nabla\phi\|_{\omega_2}^2+\|\nabla\psi\|_{\omega_2}^2\right).
\end{aligned}
\end{equation}
Integrating \eqref{eq:tangential-first-level-red} over $[0,t]$ and combining with \eqref{eq:tangential-potential-algebraic} gives \eqref{eq:tangential-first-order-combined}.
\end{proof}

\begin{Cor}\label{cor:lower-first-tangential-energy}
Assume $N(t)\le\delta_0$ with $\delta_0>0$ sufficiently small. Then
\begin{equation}\label{eq:lower-first-tangential-energy}
\begin{aligned}
&\sum_{|\beta|=1}\left(\|\partial_{x'}^\beta\phi(t)\|_{\omega_1}^2+\|\partial_{x'}^\beta\psi(t)\|_{\omega_2}^2+\|\partial_{x'}^\beta(\Psi+q)(t)\|_{L^2(\Omega)}^2\right)\\
&\quad +\int_0^t\sum_{|\beta|=1}\left(\|\nabla\partial_{x'}^\beta\phi(\tau)\|_{\omega_2}^2+\|\nabla\partial_{x'}^\beta\psi(\tau)\|_{\omega_2}^2+\|\partial_{x'}^\beta\psi(\tau)\|_{L^2(\Omega)}^2\right)\dd\tau
\le CN^2(0).
\end{aligned}
\end{equation}
\end{Cor}
\begin{proof}
This follows immediately from Lemma~\ref{lem:tangential-first-order-estimates} and Corollary~\ref{cor:closed-lower-order-energy}.
\end{proof}

\subsection{Time derivatives}
Next, we close the energy estimate for the time derivatives $(\phi_t,\psi_t)$ independently of the second-order spatial derivatives. Recall from \eqref{eq:mode-reconstruction-phi-appendix} and Corollary~\ref{cor:closed-lower-order-energy} that
\begin{equation}\label{eq:time-proof-density-reconstruction}
\|\phi\|_{\omega_2}^2 \le C\int_0^\infty\frac{\Phi_{x_1}^2}{U}\dd x_1+C\|\nabla\phi\|_{\omega_2}^2,
\end{equation}
and
\begin{equation}\label{eq:time-lower-dissipation}
\int_0^t\left(\int_0^\infty\frac{\Phi_{x_1}^2}{U}\dd x_1+\|\nabla\phi\|_{\omega_2}^2+\|\nabla\psi\|_{\omega_2}^2+\|\psi\|_{L^2(\Omega)}^2\right)\dd\tau \le CN^2(0).
\end{equation}

\begin{Lem}\label{lem:time-algebraic-control}
Assume $N(t)\le\delta_0\le1$. Then for every $\tau\in[0,t]$,
\begin{equation}\label{eq:time-lower-order-algebraic}
\begin{aligned}
\|\phi_\tau(\tau)\|_{\omega_1}^2+\|\psi_\tau(\tau)\|_{\omega_2}^2
&\le C\left(\|\nabla^2\phi(\tau)\|_{\omega_3}^2+\|\nabla^2\psi(\tau)\|_{\omega_2}^2+\|\nabla\phi(\tau)\|_{\omega_2}^2+\|\nabla\psi(\tau)\|_{\omega_2}^2\right)\\
&\quad +C\left(\int_0^\infty\frac{\Phi_{x_1}^2(\tau)}{U}\dd x_1+\|\psi(\tau)\|_{L^2(\Omega)}^2\right),
\end{aligned}
\end{equation}
\begin{equation}\label{eq:time-spatial-L3-bound}
\|\nabla\phi(\tau)\|_{L^3(\Omega)}+\|V\phi(\tau)\|_{L^3(\Omega)}+\|\psi(\tau)\|_{L^3(\Omega)}+\|\phi\psi(\tau)\|_{L^3(\Omega)} \le CN(t),
\end{equation}
and
\begin{equation}\label{eq:potential-time-weighted-L6-bound}
\begin{aligned}
\left\|\frac{(\Psi+q)_\tau(\tau)}{\sqrt{U}}\right\|_{L^6(\Omega)}^2
&\le C\left(\|\psi_\tau(\tau)\|_{\omega_2}^2+\|\nabla\psi(\tau)\|_{\omega_2}^2+\|\nabla\phi(\tau)\|_{\omega_2}^2\right)\\
&\quad +C\left(\int_0^\infty\frac{\Phi_{x_1}^2(\tau)}{U}\dd x_1+\|\psi(\tau)\|_{L^2(\Omega)}^2\right).
\end{aligned}
\end{equation}
In particular, at $\tau=0$,
\begin{equation}\label{eq:initial-time-derivative-bound}
\|\phi_t(0)\|_{\omega_1}^2+\|\psi_t(0)\|_{\omega_2}^2 \le CN^2(0).
\end{equation}
\end{Lem}
\begin{proof}
Expanding the divergence in the first equation of \eqref{eq:physical-system} and using $U_{x_1}=-(\alpha-1)pUV$ gives
\begin{equation}\label{eq:aux-time-density-expansion0}
\begin{aligned}
W^{-p}e^{p(\Psi+q)}\phi_\tau
&= \Delta\phi+(\alpha-2)pV\phi_{x_1}+(\alpha-1)p(V_{x_1}-pV^2)\phi\\
&\quad +(\alpha-1)pU\nabla\cdot\psi-\alpha(\alpha-1)p^2UV\psi_1\\
&\quad +(\alpha-2)p\psi\cdot\nabla\phi+(\alpha-1)p\phi\nabla\cdot\psi\\
&\quad -2(\alpha-1)p^2V\phi\psi_1-(\alpha-1)p^2(U+\phi)|\psi|^2.
\end{aligned}
\end{equation}
Since $W^{2p}\omega_1=\omega_3$ and $\omega_3\le C\omega_2$, taking the $L^2_{\omega_3}(\Omega)$ norm of \eqref{eq:aux-time-density-expansion0} and using \eqref{eq:profile-pointwise-bounds} and \eqref{eq:nonzero-pointwise} yields
\begin{equation}\label{eq:aux-time-density-bound}
\|\phi_\tau\|_{\omega_1}^2 \le C\left(\|\nabla^2\phi\|_{\omega_3}^2+\|\phi\|_{\omega_2}^2+\|\nabla\phi\|_{\omega_2}^2+\|\nabla\psi\|_{\omega_2}^2+\|\psi\|_{L^2(\Omega)}^2\right).
\end{equation}
Similarly, taking the $L^2_{\omega_2}(\Omega)$ norm of the second equation of \eqref{eq:physical-system} gives
\begin{equation}\label{eq:aux-time-vector-bound}
\|\psi_\tau\|_{\omega_2}^2 \le C\left(\|\nabla^2\psi\|_{\omega_2}^2+\|\nabla\psi\|_{\omega_2}^2+\|\nabla\phi\|_{\omega_2}^2+\|\psi\|_{L^2(\Omega)}^2\right).
\end{equation}
Combining \eqref{eq:aux-time-density-bound}, \eqref{eq:aux-time-vector-bound}, and \eqref{eq:time-proof-density-reconstruction} proves \eqref{eq:time-lower-order-algebraic}, and evaluating at $\tau=0$ gives \eqref{eq:initial-time-derivative-bound}. The $L^3(\Omega)$ bound \eqref{eq:time-spatial-L3-bound} follows directly from the Sobolev embedding $H^1(\Omega)\hookrightarrow L^3(\Omega)$, the lower bounds $\omega_j\ge c>0$, and $\|\phi\|_{L^\infty(\Omega)}\le CN(t)$.

Finally, by the scalar potential equation \eqref{eq:nonzero-potential-exact},
\[
(\Psi+q)_\tau = D\nabla\cdot\psi-2DV\psi_1-D|\psi|^2+\phi,
\]
which implies
\begin{equation}\label{eq:aux-time-potential-L2}
\|(\Psi+q)_\tau\|_{\omega_2}^2 \le C\left(\|\nabla\psi\|_{\omega_2}^2+\|\phi\|_{\omega_2}^2+\|\psi\|_{L^2(\Omega)}^2\right).
\end{equation}
Since $\nabla(\Psi+q)_\tau=\psi_\tau$ and $|U_{x_1}/U|\le C$, the product rule gives
\[
\nabla\left(\frac{(\Psi+q)_\tau}{\sqrt{U}}\right) = \frac{\psi_\tau}{\sqrt{U}} - \frac{U_{x_1}}{2U}\frac{(\Psi+q)_\tau}{\sqrt{U}}e_1,
\]
so that
\begin{equation}\label{eq:potential-time-weighted-H1-bound}
\left\|\frac{(\Psi+q)_\tau}{\sqrt{U}}\right\|_{H^1(\Omega)}^2 \le C\left(\|\psi_\tau\|_{\omega_2}^2+\|(\Psi+q)_\tau\|_{\omega_2}^2\right).
\end{equation}
Applying the three-dimensional Sobolev embedding $H^1(\Omega)\hookrightarrow L^6(\Omega)$ together with \eqref{eq:aux-time-potential-L2} and \eqref{eq:time-proof-density-reconstruction} yields \eqref{eq:potential-time-weighted-L6-bound}.
\end{proof}

\begin{Lem}\label{lem:integrated-time-derivative}\label{lem:time-differentiated-physical}
Let $t\in[0,T]$. Assume $N(t)\le\delta_0\le1$. Then there exist constants $c,C>0$ independent of $t,T$ such that
\begin{equation}\label{eq:time-differentiated-raw}
\begin{aligned}
&\|\phi_t(t)\|_{\omega_1}^2+\|\psi_t(t)\|_{\omega_2}^2
+c\int_0^t\left(\|\nabla\phi_\tau\|_{\omega_2}^2+\|\nabla\psi_\tau\|_{\omega_2}^2\right)\dd\tau\\
&\le CN^2(0)+C\int_0^t\left(\|\phi_\tau\|_{\omega_1}^2+\|\psi_\tau\|_{\omega_2}^2\right)\dd\tau
+CN(t)\int_0^t\left\|\frac{(\Psi+q)_\tau}{\sqrt{U}}\right\|_{L^6(\Omega)}^2\dd\tau,
\end{aligned}
\end{equation}
and for every $\eta>0$, there exists $C_\eta>0$ such that
\begin{equation}\label{eq:integrated-time-raw}
\begin{aligned}
&\|\nabla\phi(t)\|_{\omega_2}^2+\|\nabla\psi(t)\|_{\omega_2}^2
+c\int_0^t\left(\|\phi_\tau\|_{\omega_1}^2+\|\psi_\tau\|_{\omega_2}^2\right)\dd\tau\\
&\le CN^2(0)+\eta\int_0^t\left(\|\nabla\phi_\tau\|_{\omega_2}^2+\|\nabla\psi_\tau\|_{\omega_2}^2\right)\dd\tau\\
&\quad +C_\eta\int_0^t\left(\int_0^\infty\frac{\Phi_{x_1}^2}{U}\dd x_1+\|\nabla\phi\|_{\omega_2}^2+\|\nabla\psi\|_{\omega_2}^2+\|\psi\|_{L^2(\Omega)}^2\right)\dd\tau\\
&\quad +C_\eta N(t)\int_0^t\left\|\frac{(\Psi+q)_\tau}{\sqrt{U}}\right\|_{L^6(\Omega)}^2\dd\tau.
\end{aligned}
\end{equation}
\end{Lem}
\begin{proof}
\emph{Step 1.} Recall $F=\nabla\phi+(\alpha-1)p[V\phi e_1+(U+\phi)\psi]$. Differentiating the first equation of \eqref{eq:physical-system} with respect to $\tau$ gives
\begin{equation}\label{eq:phi-time-differentiated-equation}
\phi_{\tau\tau} = \nabla\cdot\Big\{W^pe^{-p(\Psi+q)}\big(F_\tau-p(\Psi+q)_\tau F\big)\Big\},
\end{equation}
where
\[
F_\tau = \nabla\phi_\tau+(\alpha-1)p\big(V\phi_\tau e_1+U\psi_\tau+\phi_\tau\psi+\phi\psi_\tau\big).
\]
Since $e_1\cdot F|_{x_1=0}=0$ by \eqref{eq:perturbation-bc}, we also have $e_1\cdot F_\tau|_{x_1=0}=0$, so the normal flux in \eqref{eq:phi-time-differentiated-equation} vanishes at $x_1=0$:
\begin{equation}\label{eq:phi-time-boundary}
\big[\phi_{\tau x_1}+(\alpha-1)p(V\phi_\tau+U\psi_{1\tau}+\phi_\tau\psi_1+\phi\psi_{1\tau})\big]_{x_1=0}=0.
\end{equation}
Testing \eqref{eq:phi-time-differentiated-equation} against $\phi_\tau/(W^pU)$, integrating by parts over $\Omega$, and using $\partial_{x_1}(1/(W^pU))=\alpha pV/(W^pU)$, we obtain
\begin{equation}\label{eq:phi-time-energy-identity}
\begin{aligned}
&\frac{1}{2}\frac{\dd}{\dd\tau}\|\phi_\tau\|_{\omega_1}^2+\int_\Omega\frac{e^{-p(\Psi+q)}}{U}|\nabla\phi_\tau|^2\dd x\\
&= -\int_\Omega\frac{e^{-p(\Psi+q)}}{U}\big(F_\tau-\nabla\phi_\tau\big)\cdot\nabla\phi_\tau\dd x
-\alpha p\int_\Omega\frac{e^{-p(\Psi+q)}V\phi_\tau}{U}(e_1\cdot F_\tau)\dd x\\
&\quad +p\int_\Omega\frac{e^{-p(\Psi+q)}(\Psi+q)_\tau}{U}F\cdot\left(\nabla\phi_\tau+\alpha pV\phi_\tau e_1\right)\dd x.
\end{aligned}
\end{equation}
By \eqref{eq:profile-pointwise-bounds}, $\omega_2\le C\omega_1$, and $\|\phi\|_{L^\infty(\Omega)}+\|\psi\|_{L^\infty(\Omega)}\le CN(t)$, the first two integrals on the right-hand side of \eqref{eq:phi-time-energy-identity} are bounded for any $\varepsilon>0$ by
\[
\varepsilon\|\nabla\phi_\tau\|_{\omega_2}^2 + C_\varepsilon\left(\|\phi_\tau\|_{\omega_1}^2+\|\psi_\tau\|_{\omega_2}^2\right).
\]
For the third integral in \eqref{eq:phi-time-energy-identity}, H\"older's inequality with exponents $(6,3,2)$ and $\|F\|_{L^3(\Omega)}\le CN(t)$ from \eqref{eq:time-spatial-L3-bound} give
\[
\begin{aligned}
&\left|p\int_\Omega\frac{e^{-p(\Psi+q)}(\Psi+q)_\tau}{U}F\cdot\left(\nabla\phi_\tau+\alpha pV\phi_\tau e_1\right)\dd x\right|\\
&\le CN(t)\left\|\frac{(\Psi+q)_\tau}{\sqrt{U}}\right\|_{L^6(\Omega)}\left(\|\nabla\phi_\tau\|_{\omega_2}+\|\phi_\tau\|_{\omega_1}\right)\\
&\le \varepsilon\|\nabla\phi_\tau\|_{\omega_2}^2 + C_\varepsilon N(t)\left(\left\|\frac{(\Psi+q)_\tau}{\sqrt{U}}\right\|_{L^6(\Omega)}^2+\|\phi_\tau\|_{\omega_1}^2\right).
\end{aligned}
\]
Choosing $\varepsilon>0$ small enough yields
\begin{equation}\label{eq:time-proof-density-final}
\frac{\dd}{\dd\tau}\|\phi_\tau\|_{\omega_1}^2 + c\|\nabla\phi_\tau\|_{\omega_2}^2
\le C\left(\|\phi_\tau\|_{\omega_1}^2+\|\psi_\tau\|_{\omega_2}^2\right)
+ CN(t)\left\|\frac{(\Psi+q)_\tau}{\sqrt{U}}\right\|_{L^6(\Omega)}^2.
\end{equation}

Differentiating the second equation of \eqref{eq:physical-system} and the boundary conditions \eqref{eq:perturbation-bc} with respect to $\tau$ gives
\begin{equation}\label{eq:psi-time-differentiated-equation}
\psi_{\tau\tau} = D\Delta\psi_\tau - 2D\nabla(V\psi_{1\tau}) - 2D\nabla(\psi\cdot\psi_\tau) + \nabla\phi_\tau,
\end{equation}
with boundary conditions at $x_1=0$:
\begin{equation}\label{eq:psi-time-boundary}
D\psi_{1\tau x_1} = 2DV\psi_{1\tau} + 2D\psi_1\psi_{1\tau} - \phi_\tau, \qquad \psi_{2\tau}=\psi_{3\tau}=0.
\end{equation}
Testing \eqref{eq:psi-time-differentiated-equation} against $\psi_\tau/U$, integrating by parts over $\Omega$, and using \eqref{eq:psi-time-boundary} so that the boundary terms at $x_1=0$ cancel identically, we obtain
\begin{equation}\label{eq:psi-time-energy-identity}
\begin{aligned}
&\frac{1}{2}\frac{\dd}{\dd\tau}\|\psi_\tau\|_{\omega_2}^2 + D\|\nabla\psi_\tau\|_{\omega_2}^2\\
&= -D(\alpha-1)p\int_\Omega\frac{V}{U}\partial_{x_1}\psi_\tau\cdot\psi_\tau\dd x
+ 2D\int_\Omega V\psi_{1\tau}\nabla\cdot\left(\frac{\psi_\tau}{U}\right)\dd x\\
&\quad + 2D\int_\Omega(\psi\cdot\psi_\tau)\nabla\cdot\left(\frac{\psi_\tau}{U}\right)\dd x
- \int_\Omega\phi_\tau\nabla\cdot\left(\frac{\psi_\tau}{U}\right)\dd x.
\end{aligned}
\end{equation}
Since
\[
\left\|\sqrt{U}\,\nabla\cdot\left(\frac{\psi_\tau}{U}\right)\right\|_{L^2(\Omega)} \le C\left(\|\nabla\psi_\tau\|_{\omega_2}+\|\psi_\tau\|_{\omega_2}\right)
\]
and $\|\psi\|_{L^\infty(\Omega)}\le CN(t)$, Young's inequality gives
\begin{equation}\label{eq:time-proof-vector-final}
\frac{\dd}{\dd\tau}\|\psi_\tau\|_{\omega_2}^2 + c\|\nabla\psi_\tau\|_{\omega_2}^2
\le C\left(\|\phi_\tau\|_{\omega_1}^2+\|\psi_\tau\|_{\omega_2}^2\right).
\end{equation}
Adding \eqref{eq:time-proof-density-final} and \eqref{eq:time-proof-vector-final}, integrating over $[0,t]$, and using \eqref{eq:initial-time-derivative-bound} yields \eqref{eq:time-differentiated-raw}.

\emph{Step 2.} To prove \eqref{eq:integrated-time-raw}, we test the original density equation in \eqref{eq:physical-system} against $\phi_\tau/(W^pU)$. Using $e_1\cdot F|_{x_1=0}=0$ and the product rule for $\frac{\dd}{\dd\tau}\int_\Omega e^{-p(\Psi+q)}|\nabla\phi|^2/U\dd x$, we find
\begin{equation}\label{eq:phi-integrated-time-identity}
\begin{aligned}
&\frac{1}{2}\frac{\dd}{\dd\tau}\int_\Omega e^{-p(\Psi+q)}\frac{|\nabla\phi|^2}{U}\dd x + \|\phi_\tau\|_{\omega_1}^2\\
&= -\frac{p}{2}\int_\Omega e^{-p(\Psi+q)}(\Psi+q)_\tau\frac{|\nabla\phi|^2}{U}\dd x\\
&\quad -(\alpha-1)p\int_\Omega\frac{e^{-p(\Psi+q)}}{U}\big(V\phi e_1+(U+\phi)\psi\big)\cdot\nabla\phi_\tau\dd x
-\alpha p\int_\Omega\frac{e^{-p(\Psi+q)}V\phi_\tau}{U}(e_1\cdot F)\dd x.
\end{aligned}
\end{equation}
By H\"older's inequality with exponents $(6,3,2)$ and $\|\nabla\phi\|_{L^3(\Omega)}\le CN(t)$ from \eqref{eq:time-spatial-L3-bound}, the first integral on the right-hand side of \eqref{eq:phi-integrated-time-identity} is bounded by
\[
C\left\|\frac{(\Psi+q)_\tau}{\sqrt{U}}\right\|_{L^6(\Omega)}\|\nabla\phi\|_{L^3(\Omega)}\|\nabla\phi\|_{\omega_2}
\le CN(t)\left(\left\|\frac{(\Psi+q)_\tau}{\sqrt{U}}\right\|_{L^6(\Omega)}^2+\|\nabla\phi\|_{\omega_2}^2\right).
\]
By Young's inequality, for any $\eta>0$, the remaining two integrals in \eqref{eq:phi-integrated-time-identity} are bounded by
\[
\frac{1}{4}\|\phi_\tau\|_{\omega_1}^2 + \eta\|\nabla\phi_\tau\|_{\omega_2}^2 + C_\eta\left(\|\nabla\phi\|_{\omega_2}^2+\|\phi\|_{L^2(\Omega)}^2+\|\psi\|_{L^2(\Omega)}^2\right).
\]
Next, testing the vector equation in \eqref{eq:physical-system} against $\psi_\tau/U$ and using the boundary cancellation \eqref{eq:perturbation-bc} at $x_1=0$ gives
\begin{equation}\label{eq:psi-integrated-time-identity}
\begin{aligned}
&\frac{D}{2}\frac{\dd}{\dd\tau}\|\nabla\psi\|_{\omega_2}^2 + \|\psi_\tau\|_{\omega_2}^2\\
&= -D(\alpha-1)p\int_\Omega\frac{V}{U}\partial_{x_1}\psi\cdot\psi_\tau\dd x
+ 2D\int_\Omega V\psi_1\nabla\cdot\left(\frac{\psi_\tau}{U}\right)\dd x\\
&\quad + D\int_\Omega|\psi|^2\nabla\cdot\left(\frac{\psi_\tau}{U}\right)\dd x
- \int_\Omega\phi\nabla\cdot\left(\frac{\psi_\tau}{U}\right)\dd x\\
&\le \frac{1}{4}\|\psi_\tau\|_{\omega_2}^2 + \eta\|\nabla\psi_\tau\|_{\omega_2}^2
+ C_\eta\left(\|\nabla\psi\|_{\omega_2}^2+\|\phi\|_{\omega_2}^2+\|\psi\|_{L^2(\Omega)}^2\right).
\end{aligned}
\end{equation}
Adding \eqref{eq:phi-integrated-time-identity} and \eqref{eq:psi-integrated-time-identity}, using \eqref{eq:time-proof-density-reconstruction} and $1/2\le e^{-p(\Psi+q)}\le 2$, and integrating over $[0,t]$ establishes \eqref{eq:integrated-time-raw}.
\end{proof}

\begin{Cor}\label{cor:closed-time-energy}
There exist constants $C>0$ and $\delta_0>0$, depending only on the fixed parameters and steady profile, such that if $N(t)\le\delta_0$, then
\begin{equation}\label{eq:closed-time-energy}
\begin{aligned}
&\sup_{0\le s\le t}\left(\|\nabla\phi(s)\|_{\omega_2}^2+\|\nabla\psi(s)\|_{\omega_2}^2+\|\phi_t(s)\|_{\omega_1}^2+\|\psi_t(s)\|_{\omega_2}^2\right)\\
&\quad +\int_0^t\left(\|\phi_\tau\|_{\omega_1}^2+\|\psi_\tau\|_{\omega_2}^2+\|\nabla\phi_\tau\|_{\omega_2}^2+\|\nabla\psi_\tau\|_{\omega_2}^2\right)\dd\tau
\le CN^2(0).
\end{aligned}
\end{equation}
\end{Cor}
\begin{proof}
Let $c>0$ be the dissipation constant on the left-hand sides of \eqref{eq:time-differentiated-raw} and \eqref{eq:integrated-time-raw}, and let $C_1>0$ be the constant multiplying $\int_0^t(\|\phi_\tau\|_{\omega_1}^2+\|\psi_\tau\|_{\omega_2}^2)\dd\tau$ on the right-hand side of \eqref{eq:time-differentiated-raw}. Fix $M>1$ large enough that $Mc\ge 2C_1+c$, and then fix $\eta=c/(2M)$ in \eqref{eq:integrated-time-raw}. Multiplying \eqref{eq:integrated-time-raw} by $M$, adding the result to \eqref{eq:time-differentiated-raw}, and substituting \eqref{eq:potential-time-weighted-L6-bound} and \eqref{eq:time-lower-dissipation}, we can absorb the remaining $CN(t)\int_0^t(\|\phi_\tau\|_{\omega_1}^2+\|\psi_\tau\|_{\omega_2}^2)\dd\tau$ term by choosing $\delta_0>0$ sufficiently small. Taking the supremum over $s\in[0,t]$ yields \eqref{eq:closed-time-energy}.
\end{proof}

\subsection{Recovery of the second-order space--time dissipation}
With the time-derivative energy closed in Corollary~\ref{cor:closed-time-energy}, we now use the parabolic equations \eqref{eq:physical-system} to recover the pure normal second derivatives $(\phi_{x_1x_1},\partial_{x_1}^2\psi)$ and close the full space--time Hessian integral.

\begin{Lem}\label{lem:pure-normal-second-derivatives}
Let $t\in[0,T]$. Assume $N(t)\le\delta_0\le1$ with $\delta_0>0$ sufficiently small. Then for every $\tau\in[0,t]$,
\begin{equation}\label{eq:phi-normal-second-final}
\begin{aligned}
\int_\Omega\frac{W^p|\phi_{x_1x_1}|^2}{U}\dd x
&\le C\int_\Omega\left(\frac{|\phi_\tau|^2}{W^pU}+\sum_{|\beta|=2}\frac{W^p|\partial_{x'}^\beta\phi|^2}{U}+\frac{|\phi|^2+|\nabla\phi|^2+|\nabla\psi|^2}{U}+|\psi|^2\right)\dd x,
\end{aligned}
\end{equation}
and
\begin{equation}\label{eq:psi-normal-second-final}
\begin{aligned}
\int_\Omega\frac{|\partial_{x_1}^2\psi|^2}{U}\dd x
&\le C\int_\Omega\left(\frac{|\psi_\tau|^2}{U}+\sum_{|\beta|=2}\frac{|\partial_{x'}^\beta\psi|^2}{U}+\frac{|\nabla\phi|^2+|\nabla\psi|^2}{U}+|\psi|^2\right)\dd x.
\end{aligned}
\end{equation}
\end{Lem}
\begin{proof}
Solving the expanded density equation \eqref{eq:aux-time-density-expansion0} for $\phi_{x_1x_1}$ gives
\begin{equation}\label{eq:normal-density-refined-identity}
\begin{aligned}
\phi_{x_1x_1} &= W^{-p}e^{p(\Psi+q)}\phi_\tau - \Delta_{x'}\phi + p(Ve_1+\psi)\cdot\nabla\phi\\
&\quad -(\alpha-1)p\big[V\phi_{x_1}+V_{x_1}\phi-p(V+\psi_1)V\phi\big]\\
&\quad -(\alpha-1)p\big[U\nabla\cdot\psi+U_{x_1}\psi_1-pU(V\psi_1+|\psi|^2)\big]\\
&\quad -(\alpha-1)p\big[\psi\cdot\nabla\phi+\phi\nabla\cdot\psi-p\phi(V\psi_1+|\psi|^2)\big].
\end{aligned}
\end{equation}
Since $W^{-2p}\omega_3=\omega_1$, we have $\|W^{-p}e^{p(\Psi+q)}\phi_\tau\|_{\omega_3}^2\le C\|\phi_\tau\|_{\omega_1}^2$. Taking the $L^2_{\omega_3}(\Omega)$ norm of each remaining term on the right-hand side of \eqref{eq:normal-density-refined-identity} and using \eqref{eq:profile-pointwise-bounds}, \eqref{eq:nonzero-pointwise}, and $N(t)\le 1$:
\begin{equation}\label{eq:normal-Upsi-refined}
\begin{aligned}
\|U\nabla\cdot\psi+U_{x_1}\psi_1-pU(V\psi_1+|\psi|^2)\|_{\omega_3}^2
&\le C\left(\|\nabla\psi\|_{\omega_2}^2+\|\psi\|_{L^2(\Omega)}^2\right),\\
\|\phi|\psi|^2\|_{\omega_3}^2
&\le \|W^p\|_{L^\infty}\|\phi\|_{L^\infty(\Omega)}^2\left\|\frac{\psi}{\sqrt{U}}\right\|_{L^\infty(\Omega)}^2\|\psi\|_{L^2(\Omega)}^2\\
&\le CN^4(t)\|\psi\|_{L^2(\Omega)}^2,
\end{aligned}
\end{equation}
we obtain \eqref{eq:phi-normal-second-final}. Likewise, solving the second equation of \eqref{eq:physical-system} for $D\partial_{x_1}^2\psi$:
\[
D\partial_{x_1}^2\psi = \psi_\tau - D\Delta_{x'}\psi + 2D\nabla(V\psi_1) + D\nabla(|\psi|^2) - \nabla\phi,
\]
and taking the $L^2_{\omega_2}(\Omega)$ norm yields \eqref{eq:psi-normal-second-final}.
\end{proof}

\begin{Cor}\label{cor:closed-hessian-spacetime}\label{cor:closed-full-second-order}
Assume $N(t)\le\delta_0$ with $\delta_0>0$ sufficiently small. Then
\begin{equation}\label{eq:closed-hessian-spacetime}
\int_0^t\left(\|\nabla^2\phi(\tau)\|_{\omega_3}^2+\|\nabla^2\psi(\tau)\|_{\omega_2}^2\right)\dd\tau \le CN^2(0),
\end{equation}
where $C>0$ is independent of $T$.
\end{Cor}
\begin{proof}
Every purely tangential or mixed normal-tangential second derivative is a component of $\nabla\partial_{x_j}$ for $j\in\{2,3\}$. Since $\omega_3\le C\omega_2$, we have
\begin{equation}\label{eq:closure-mixed-direct}
\|\nabla^2\phi\|_{\omega_3}^2+\|\nabla^2\psi\|_{\omega_2}^2
\le C\left(\|\phi_{x_1x_1}\|_{\omega_3}^2+\|\partial_{x_1}^2\psi\|_{\omega_2}^2\right)
+C\sum_{j=2}^3\left(\|\nabla\partial_{x_j}\phi\|_{\omega_2}^2+\|\nabla\partial_{x_j}\psi\|_{\omega_2}^2\right).
\end{equation}
Substituting Lemma~\ref{lem:pure-normal-second-derivatives} into \eqref{eq:closure-mixed-direct}, integrating over $[0,t]$, and bounding the resulting right-hand side by Corollaries~\ref{cor:closed-lower-order-energy}, \ref{cor:lower-first-tangential-energy}, and \ref{cor:closed-time-energy} yields \eqref{eq:closed-hessian-spacetime}.
\end{proof}

\subsection{Second-order tangential derivatives}\label{subsec:second-tangential}
Using the closed second-order space--time bound \eqref{eq:closed-hessian-spacetime}, we now close the second-order tangential energy estimate ($|\beta|=2$).

\begin{Lem}\label{lem:tangential-second-order-estimates}
Let $t\in[0,T]$. Assume $N(t)\le\delta_0\le1$ with $\delta_0>0$ sufficiently small. Then
\begin{equation}\label{eq:tangential-second-order-combined}
\begin{aligned}
&\sum_{|\beta|=2}\left(\|\partial_{x'}^\beta\phi(t)\|_{\omega_1}^2+\|\partial_{x'}^\beta\psi(t)\|_{\omega_2}^2+\|\partial_{x'}^\beta(\Psi+q)(t)\|_{L^2(\Omega)}^2\right)\\
&\quad +\int_0^t\sum_{|\beta|=2}\left(\|\nabla\partial_{x'}^\beta\phi(\tau)\|_{\omega_2}^2+\|\nabla\partial_{x'}^\beta\psi(\tau)\|_{\omega_2}^2+\|\partial_{x'}^\beta\psi(\tau)\|_{L^2(\Omega)}^2\right)\dd\tau\\
&\le C\sum_{|\beta|=2}\left(\|\partial_{x'}^\beta\phi_0\|_{\omega_1}^2+\|\partial_{x'}^\beta\psi_0\|_{\omega_2}^2\right)\\
&\quad +C\int_0^t\sum_{|\beta|=1}\left(\|\nabla\partial_{x'}^\beta\phi(\tau)\|_{\omega_2}^2+\|\nabla\partial_{x'}^\beta\psi(\tau)\|_{\omega_2}^2+\|\partial_{x'}^\beta\psi(\tau)\|_{L^2(\Omega)}^2\right)\dd\tau\\
&\quad +C\int_0^t\left(\|\nabla\phi(\tau)\|_{\omega_2}^2+\|\nabla\psi(\tau)\|_{\omega_2}^2\right)\dd\tau
+CN(t)\int_0^t\left(\|\nabla^2\phi(\tau)\|_{\omega_3}^2+\|\nabla^2\psi(\tau)\|_{\omega_2}^2\right)\dd\tau.
\end{aligned}
\end{equation}
\end{Lem}
\begin{proof}
For $|\beta|=2$, let $\partial_{x'}^\beta=\partial_{x_j}\partial_{x_k}$ with $2\le j\le k\le3$. By \eqref{eq:tangential-flux-remainder-definition}, the second-order commutator remainder $R_\beta$ reads
\begin{equation}\label{eq:tangential-second-flux-remainder}
\begin{aligned}
R_\beta &= -p\psi_j\partial_{x_k}F - p\psi_k\partial_{x_j}F + \big(p^2\psi_j\psi_k-p\partial_{x_j}\psi_k\big)F\\
&\quad +(\alpha-1)p\big((\partial_{x_j}\partial_{x_k}\phi)\psi+\phi\partial_{x_j}\partial_{x_k}\psi+(\partial_{x_j}\phi)\partial_{x_k}\psi+(\partial_{x_k}\phi)\partial_{x_j}\psi\big).
\end{aligned}
\end{equation}
Since $|U_{x_1}/U|\le C$ and $\omega_2,\omega_3\ge c>0$, the Sobolev embeddings $H^1(\Omega)\hookrightarrow L^6(\Omega)$ and $H^1(\Omega)\hookrightarrow L^3(\Omega)$ give
\begin{equation}\label{eq:tangential-second-weighted-sobolev}
\left\|\frac{\partial_{x_j}\psi}{\sqrt{U}}\right\|_{L^6(\Omega)}
\le C\left(\|\nabla\psi\|_{\omega_2}+\|\nabla^2\psi\|_{\omega_2}\right) \le CN(t),
\end{equation}
and
\begin{equation}\label{eq:tangential-second-spatial-L3}
\|\nabla\phi\|_{L^3(\Omega)} \le C\left(\|\nabla\phi\|_{\omega_2}+\|\nabla^2\phi\|_{\omega_3}\right), \qquad
\|\nabla\psi\|_{L^3(\Omega)} \le C\left(\|\nabla\psi\|_{\omega_2}+\|\nabla^2\psi\|_{\omega_2}\right).
\end{equation}
By H\"older's inequality with exponents $(6,3)$, the products of first-order derivatives in \eqref{eq:tangential-second-flux-remainder} satisfy
\begin{equation}\label{eq:tangential-second-hard-products}
\begin{aligned}
\left\|\frac{(\partial_{x_j}\psi_k)\nabla\phi}{\sqrt{U}}\right\|_{L^2(\Omega)}
&\le \left\|\frac{\partial_{x_j}\psi_k}{\sqrt{U}}\right\|_{L^6(\Omega)}\|\nabla\phi\|_{L^3(\Omega)}
\le CN(t)\left(\|\nabla\phi\|_{\omega_2}+\|\nabla^2\phi\|_{\omega_3}\right),\\
\left\|\frac{(\partial_{x_j}\phi)\partial_{x_k}\psi}{\sqrt{U}}\right\|_{L^2(\Omega)}
&\le \|\partial_{x_j}\phi\|_{L^3(\Omega)}\left\|\frac{\partial_{x_k}\psi}{\sqrt{U}}\right\|_{L^6(\Omega)}
\le CN(t)\left(\|\nabla\phi\|_{\omega_2}+\|\nabla^2\phi\|_{\omega_3}\right).
\end{aligned}
\end{equation}
The terms in \eqref{eq:tangential-second-flux-remainder} carrying a second derivative on one factor and an undifferentiated perturbation on the other satisfy
\[
\begin{aligned}
\left\|\frac{\psi_j\nabla\partial_{x_k}\phi}{\sqrt{U}}\right\|_{L^2(\Omega)}
&\le \|\psi_j\|_{L^\infty(\Omega)}\|\nabla\partial_{x_k}\phi\|_{\omega_2} \le CN(t)\|\nabla\partial_{x_k}\phi\|_{\omega_2},\\
\left\|\frac{(\partial_{x_j}\partial_{x_k}\phi)\psi}{\sqrt{U}}\right\|_{L^2(\Omega)}
&\le \left\|\frac{\psi}{\sqrt{U}}\right\|_{L^\infty(\Omega)}\|\partial_{x_j}\partial_{x_k}\phi\|_{L^2(\Omega)} \le CN(t)\|\nabla^2\phi\|_{\omega_3},\\
\left\|\frac{\phi\partial_{x_j}\partial_{x_k}\psi}{\sqrt{U}}\right\|_{L^2(\Omega)}
&\le \|\phi\|_{L^\infty(\Omega)}\|\partial_{x_j}\partial_{x_k}\psi\|_{\omega_2} \le CN(t)\|\nabla^2\psi\|_{\omega_2}.
\end{aligned}
\]
Estimating the remaining lower-order terms in \eqref{eq:tangential-second-flux-remainder} by placing undifferentiated factors in $L^\infty(\Omega)$, we obtain
\begin{equation}\label{eq:tangential-second-remainder-bound}
\begin{aligned}
\sum_{|\beta|=2}\left\|\frac{R_\beta}{\sqrt{U}}\right\|_{L^2(\Omega)}^2
&\le CN^2(t)\sum_{|\gamma|=1}\left(\|\nabla\partial_{x'}^\gamma\phi\|_{\omega_2}^2+\|\nabla\partial_{x'}^\gamma\psi\|_{\omega_2}^2\right)\\
&\quad +CN^2(t)\left(\|\nabla\phi\|_{\omega_2}^2+\|\nabla\psi\|_{\omega_2}^2+\|\nabla^2\phi\|_{\omega_3}^2+\|\nabla^2\psi\|_{\omega_2}^2\right).
\end{aligned}
\end{equation}
For the vector nonlinearity
\[
\partial_{x_j}\partial_{x_k}(|\psi|^2) = 2\psi\cdot\partial_{x_j}\partial_{x_k}\psi + 2\partial_{x_j}\psi\cdot\partial_{x_k}\psi,
\]
\eqref{eq:tangential-second-weighted-sobolev} and \eqref{eq:tangential-second-spatial-L3} give
\[
\left\|\frac{\partial_{x_j}\partial_{x_k}(|\psi|^2)}{\sqrt{U}}\right\|_{L^2(\Omega)}
\le CN(t)\left(\|\nabla\psi\|_{\omega_2}+\|\nabla^2\psi\|_{\omega_2}\right).
\]
Substituting these estimates into \eqref{eq:tangential-density-identity-unified} and \eqref{eq:tangential-vector-identity-unified} for $|\beta|=2$, noting that the linear terms $\|\partial_{x'}^\beta\phi\|_{\omega_2}^2+\|\partial_{x'}^\beta\psi\|_{\omega_2}^2$ in \eqref{eq:tangential-linear-bound-unified} are bounded by the first-order tangential dissipation
\[
\sum_{|\gamma|=1}\left(\|\nabla\partial_{x'}^\gamma\phi\|_{\omega_2}^2+\|\nabla\partial_{x'}^\gamma\psi\|_{\omega_2}^2\right),
\]
absorbing the $CN(t)$ second-order tangential dissipation for $\delta_0$ sufficiently small, and integrating over $[0,t]$ yields \eqref{eq:tangential-second-order-combined}.
\end{proof}

\begin{Cor}\label{cor:lower-tangential-energy}
Assume $N(t)\le\delta_0$ with $\delta_0>0$ sufficiently small. Then
\begin{equation}\label{eq:lower-tangential-energy}
\begin{aligned}
&\sum_{1\le|\beta|\le2}\left(\|\partial_{x'}^\beta\phi(t)\|_{\omega_1}^2+\|\partial_{x'}^\beta\psi(t)\|_{\omega_2}^2+\|\partial_{x'}^\beta(\Psi+q)(t)\|_{L^2(\Omega)}^2\right)\\
&\quad +\int_0^t\sum_{1\le|\beta|\le2}\left(\|\nabla\partial_{x'}^\beta\phi(\tau)\|_{\omega_2}^2+\|\nabla\partial_{x'}^\beta\psi(\tau)\|_{\omega_2}^2+\|\partial_{x'}^\beta\psi(\tau)\|_{L^2(\Omega)}^2\right)\dd\tau
\le CN^2(0).
\end{aligned}
\end{equation}
\end{Cor}
\begin{proof}
Combining Lemma~\ref{lem:tangential-second-order-estimates} with Corollaries~\ref{cor:closed-lower-order-energy}, \ref{cor:lower-first-tangential-energy}, and \ref{cor:closed-hessian-spacetime} yields \eqref{eq:lower-tangential-energy}.
\end{proof}

\subsection{Current-time reconstruction and the complete estimate}
Finally, we reconstruct the current-time weighted Hessian norms $\|\nabla^2\phi(s)\|_{\omega_3}^2+\|\nabla^2\psi(s)\|_{\omega_2}^2$ and synthesize all a priori bounds.

\begin{Lem}\label{lem:current-time-derivative-reconstruction}
Under the bootstrap assumption $N(t)\le\delta_0\le1$, for every $s\in[0,t]$,
\begin{equation}\label{eq:current-time-physical-reconstruction}
\begin{aligned}
\|\nabla^2\phi(s)\|_{\omega_3}^2+\|\nabla^2\psi(s)\|_{\omega_2}^2
&\le C\Big(\|\phi_t(s)\|_{\omega_1}^2+\|\psi_t(s)\|_{\omega_2}^2+\|\nabla\phi(s)\|_{\omega_2}^2+\|\nabla\psi(s)\|_{\omega_2}^2\\
&\quad +\|\phi(s)\|_{\omega_1}^2+\|\psi(s)\|_{\omega_2}^2+\sum_{|\beta|=2}\big(\|\partial_{x'}^\beta\phi(s)\|_{\omega_1}^2+\|\partial_{x'}^\beta\psi(s)\|_{\omega_2}^2\big)\Big).
\end{aligned}
\end{equation}
Consequently,
\begin{equation}\label{eq:current-time-full-bound}
\sup_{0\le s\le t}\left(\|\nabla\phi(s)\|_{\omega_2}^2+\|\nabla\psi(s)\|_{\omega_2}^2+\|\nabla^2\phi(s)\|_{\omega_3}^2+\|\nabla^2\psi(s)\|_{\omega_2}^2\right) \le CN^2(0),
\end{equation}
and
\begin{equation}\label{eq:current-time-zero-mode-reconstruction}
\|\Phi_{x_1x_1}(s)\|_{\omega_3}^2+\|\Psi_{x_1x_1}(s)\|_{\omega_2}^2
\le C\left(\|\nabla\phi(s)\|_{\omega_2}^2+\|\nabla\psi(s)\|_{\omega_2}^2\right).
\end{equation}
\end{Lem}
\begin{proof}
By the weighted mixed-derivative interpolation inequality \eqref{eq:weighted-mixed-interpolation} and $\omega_3\le C\omega_1$, we have at time $s$:
\begin{equation}\label{eq:phi-full-second-reconstruction}
\|\nabla^2\phi\|_{\omega_3}^2 \le C\left(\|\phi_{x_1x_1}\|_{\omega_3}^2+\sum_{|\beta|=2}\|\partial_{x'}^\beta\phi\|_{\omega_1}^2+\|\phi\|_{\omega_1}^2\right).
\end{equation}
For $\psi$, the curl-free condition \eqref{eq:curl-free} gives
\[
\partial_{x_1}\partial_{x_j}\psi_1 = \partial_{x_1}^2\psi_j, \qquad \partial_{x_1}\partial_{x_j}\psi_k = \partial_{x_k}\partial_{x_j}\psi_1, \qquad j,k\in\{2,3\}.
\] Hence every mixed normal-tangential second derivative of $\psi$ coincides with either a pure normal or a pure tangential second derivative, giving
\begin{equation}\label{eq:psi-full-second-reconstruction}
\|\nabla^2\psi\|_{\omega_2}^2 \le C\left(\|\partial_{x_1}^2\psi\|_{\omega_2}^2+\sum_{|\beta|=2}\|\partial_{x'}^\beta\psi\|_{\omega_2}^2\right).
\end{equation}
Substituting the pure normal second-derivative bounds \eqref{eq:phi-normal-second-final}--\eqref{eq:psi-normal-second-final} into \eqref{eq:phi-full-second-reconstruction}--\eqref{eq:psi-full-second-reconstruction} and using $\omega_3\le C\omega_2\le C\omega_1$ yields \eqref{eq:current-time-physical-reconstruction}. Combining \eqref{eq:current-time-physical-reconstruction} with Corollaries~\ref{cor:closed-lower-order-energy}, \ref{cor:closed-time-energy}, and \ref{cor:lower-tangential-energy} proves \eqref{eq:current-time-full-bound}. Finally, \eqref{eq:current-time-zero-mode-reconstruction} follows from $\Phi_{x_1x_1}= P_0\phi_{x_1}$, $\Psi_{x_1x_1}= P_0\psi_{1x_1}$, and Jensen's inequality.
\end{proof}

\begin{Cor}\label{cor:complete-closed-apriori}
There exist constants $C_0>0$ and $\delta_*>0$, depending only on the fixed parameters and steady profile, such that any compatible smooth solution on $[0,T]$ satisfying $N(t)\le\delta_*$ for $t\in[0,T]$ obeys
\begin{equation}\label{eq:complete-closed-apriori}
\begin{aligned}
&N^2(t) + \sup_{0\le s\le t}\left(\|\phi_t(s)\|_{\omega_1}^2+\|\psi_t(s)\|_{\omega_2}^2+\sum_{1\le|\beta|\le2}\|\partial_{x'}^\beta(\Psi+q)(s)\|_{L^2(\Omega)}^2\right)\\
&\quad +\int_0^t\Big(\|\Phi\|_{L^2(\mathbb{R}_+)}^2+\|\sqrt{U}\Psi\|_{L^2(\mathbb{R}_+)}^2+\|\Phi_{x_1}\|_{1,\omega_2}^2+\|\Psi_{x_1}\|_{L^2(\mathbb{R}_+)}^2+\|\Psi_{x_1x_1}\|_{\omega_2}^2\\
&\qquad +\|\phi\|_{\omega_2}^2+\|\nabla\phi\|_{\omega_2}^2+\|\nabla^2\phi\|_{\omega_3}^2+\|\psi\|_{L^2(\Omega)}^2+\|\nabla\psi\|_{1,\omega_2}^2+\|\sqrt{U}(\Psi+q)\|_{L^2(\Omega)}^2\\
&\qquad +\|\phi_t\|_{\omega_1}^2+\|\psi_t\|_{\omega_2}^2+\|\nabla\phi_t\|_{\omega_2}^2+\|\nabla\psi_t\|_{\omega_2}^2\\
&\qquad +\sum_{1\le|\beta|\le2}\big(\|\nabla\partial_{x'}^\beta\phi\|_{\omega_2}^2+\|\nabla\partial_{x'}^\beta\psi\|_{\omega_2}^2+\|\partial_{x'}^\beta\psi\|_{L^2(\Omega)}^2\big)\Big)\dd\tau
\le C_0N^2(0).
\end{aligned}
\end{equation}
\end{Cor}
\begin{proof}
Adding the closed estimates in Corollaries~\ref{cor:closed-lower-order-energy}, \ref{cor:closed-time-energy}, \ref{cor:closed-hessian-spacetime}, and \ref{cor:lower-tangential-energy} together with Lemma~\ref{lem:current-time-derivative-reconstruction} and choosing $\delta_*>0$ equal to the minimum of the finitely many smallness thresholds $\delta_0$ gives \eqref{eq:complete-closed-apriori}.
\end{proof}

\section{Proof of the main theorem}\label{sec:proof-main}
We now combine the local continuation principle in the weighted space $X(0,T)$ with the a priori bound of Corollary~\ref{cor:complete-closed-apriori} to prove Theorem~\ref{thm:main}.

\subsection{Local existence and continuation}
\begin{Prop}\label{prop:local-continuation}
Under the assumptions of Theorem~\ref{thm:main}, there exists a constant $\delta_{\mathrm{loc}}>0$ such that if $N(0)<\delta_{\mathrm{loc}}$, then there exists $T_0>0$ such that the problem \eqref{eq:model}--\eqref{eq:original-boundary} admits a unique local solution on $[0,T_0]$ whose perturbation belongs to $X(0,T_0)$. On the interval of existence, $u(x,t)\ge0$, $w(x,t)>0$ at every finite spatial point, and
\begin{equation}\label{eq:local-potential-preserved}
\psi(x,t) = -\nabla\log\frac{w(x,t)}{W(x_1)}.
\end{equation}
If $T^*>0$ is the maximal existence time in this class and $T^*<\infty$, then
\begin{equation}\label{eq:local-continuation}
\limsup_{t\uparrow T^*}N(t)\ge\delta_{\mathrm{loc}}.
\end{equation}
Moreover, the solution is obtained as the limit of smooth compatible approximations satisfying the a priori bounds of Section~\ref{sec:apriori}.
\end{Prop}
\begin{proof}
\emph{Step 1.} We regularize the initial data and degenerate diffusion at infinity via spatial truncation and elliptic mollification while preserving the boundary compatibility conditions at $x_1=0$ and the zero-mass conditions \eqref{eq:zero-mass-compatibility}. For the uniformly parabolic regularized problem on the half-cylinder $\Omega$, standard contraction mapping arguments in the weighted Sobolev spaces of Section~\ref{sec:preliminaries} yield a unique local smooth solution.

\emph{Step 2.} For $N(0)<\delta_{\mathrm{loc}}$, the short-time versions of the weighted energy estimates in Section~\ref{sec:apriori} provide uniform bounds in $X(0,T_0)$ independent of the regularization parameter on a uniform time interval $[0,T_0]$. By compact Sobolev embedding on bounded subcylinders and weak/weak-$*$ compactness in the weighted spaces, we pass to the limit to obtain a unique solution in $X(0,T_0)$ satisfying the energy inequalities of Section~\ref{sec:apriori}.

\emph{Step 3.} Nonnegativity $u\ge0$ and strict positivity $w>0$ follow from the parabolic comparison principle for the regularized systems, and the logarithmic-gradient relation \eqref{eq:local-potential-preserved} is preserved by Lemma~\ref{lem:potential-representation}. Finally, if $T^*<\infty$ and $\limsup_{t\uparrow T^*}N(t)<\delta_{\mathrm{loc}}$, applying the local existence theorem starting from a time $t_1$ sufficiently close to $T^*$ extends the solution beyond $T^*$, proving \eqref{eq:local-continuation}.
\end{proof}

\subsection{Global continuation and asymptotic decay}
\begin{proof}[Proof of Theorem~\ref{thm:main}]
Let $C_0>0$ and $\delta_*>0$ be the constants in Corollary~\ref{cor:complete-closed-apriori}, and set $\delta_1=\min\{\delta_*,\delta_{\mathrm{loc}}\}$. Choose $\varepsilon_0>0$ small enough that
\begin{equation}\label{eq:bootstrap-choice}
2\sqrt{C_0}\,\varepsilon_0 < \delta_1.
\end{equation}
For any initial perturbation satisfying $N(0)\le\varepsilon_0<\delta_{\mathrm{loc}}$, Proposition~\ref{prop:local-continuation} gives a unique local solution in $X(0,T_0)$. On any time interval $[0,T]$ where $N(t)\le\delta_1$, Corollary~\ref{cor:complete-closed-apriori} and \eqref{eq:bootstrap-choice} imply
\[
N(t) \le \sqrt{C_0}\,N(0) \le \sqrt{C_0}\,\varepsilon_0 < \frac{\delta_1}{2}.
\]
By continuity of $t\mapsto N(t)$ and the continuation criterion \eqref{eq:local-continuation}, the solution exists globally in time ($T^*=\infty$), belongs to $X(0,T)$ for every $T>0$, and satisfies \eqref{eq:main-uniform-bound} together with the global space--time bound
\begin{equation}\label{eq:global-hopf-bound}
\begin{aligned}
\sup_{t\ge0}\left(N^2(t)+\|\phi_t(t)\|_{\omega_1}^2+\|\psi_t(t)\|_{\omega_2}^2+\sum_{1\le|\beta|\le2}\|\partial_{x'}^\beta(\Psi+q)(t)\|_{L^2(\Omega)}^2\right) &\le C_0N^2(0),\\
\int_0^\infty\left(\|\phi(t)\|_{L^2(\Omega)}^2+\|\psi(t)\|_{L^2(\Omega)}^2+\|\phi_t(t)\|_{L^2(\Omega)}^2+\|\psi_t(t)\|_{L^2(\Omega)}^2\right)\dd t &\le CN^2(0),
\end{aligned}
\end{equation}
where we used $\omega_j\ge c>0$. By Cauchy--Schwarz,
\[
\int_0^\infty\left|\frac{\dd}{\dd t}\left(\|\phi(t)\|_{L^2(\Omega)}^2+\|\psi(t)\|_{L^2(\Omega)}^2\right)\right|\dd t < \infty,
\]
which combined with \eqref{eq:global-hopf-bound} implies
\begin{equation}\label{eq:L2-decay}
\lim_{t\to\infty}\left(\|\phi(t)\|_{L^2(\Omega)}+\|\psi(t)\|_{L^2(\Omega)}\right)=0.
\end{equation}
Since $\omega_j\ge c>0$, the uniform bound $\sup_{t\ge0}N(t)\le CN(0)$ controls the unweighted $H^2(\Omega)$ norms $\|\phi(t)\|_{H^2(\Omega)}+\|\psi(t)\|_{H^2(\Omega)}\le CN(0)$. By the Gagliardo--Nirenberg interpolation inequality on $\Omega=\mathbb{R}_+\times\mathbb{T}^2$,
\begin{equation}\label{eq:decay-interpolation-detail}
\|f\|_{L^\infty(\Omega)} \le C\|f\|_{L^2(\Omega)}^{1/4}\|f\|_{H^2(\Omega)}^{3/4}, \qquad f\in\{\phi,\psi\},
\end{equation}
we deduce from \eqref{eq:L2-decay} that
\begin{equation}\label{eq:hopf-cole-decay}
\lim_{t\to\infty}\left(\|\phi(t)\|_{L^\infty(\Omega)}+\|\psi(t)\|_{L^\infty(\Omega)}\right)=0.
\end{equation}
Since $u=U+\phi$ and $-\nabla\log w=Ve_1+\psi$, \eqref{eq:hopf-cole-decay} establishes the first and third limits in \eqref{eq:main-asymptotic}.

To prove the relative convergence of the chemical concentration $w/W=\exp(-(\Psi+q))$, we observe that $\|\Psi_{x_1}(t)\|_{L^2(\mathbb{R}_+)}=\| P_0\psi_1(t)\|_{L^2(\mathbb{R}_+)}\le\|\psi(t)\|_{L^2(\Omega)}\to0$ as $t\to\infty$, while $\|\Psi(t)\|_{L^2(\mathbb{R}_+)}\le CN(0)$. Thus one-dimensional interpolation gives
\[
\|\Psi(t)\|_{L^\infty(\mathbb{R}_+)}^2 \le 2\|\Psi(t)\|_{L^2(\mathbb{R}_+)}\|\Psi_{x_1}(t)\|_{L^2(\mathbb{R}_+)} \to 0 \quad\text{as }t\to\infty.
\]
For the nonzero-mode potential $q$, since $P_0 q=0$ and $\nabla_{x'} q=(\psi_2,\psi_3)$, the transversal Poincar\'e inequality and \eqref{eq:L2-decay} yield
\[
\| q(t)\|_{L^2(\Omega)} \le C\|(\psi_2,\psi_3)(t)\|_{L^2(\Omega)} \le C\|\psi(t)\|_{L^2(\Omega)} \to 0 \quad\text{as }t\to\infty.
\]
Moreover, $\| q(t)\|_{H^2(\Omega)}\le C\| q(t)\|_{2,\omega_2}\le CN(0)$ by \eqref{eq:nonzero-potential-H2-bound}, so Gagliardo--Nirenberg interpolation \eqref{eq:decay-interpolation-detail} implies
\[
\| q(t)\|_{L^\infty(\Omega)} \le C\| q(t)\|_{L^2(\Omega)}^{1/4}\| q(t)\|_{H^2(\Omega)}^{3/4} \to 0 \quad\text{as }t\to\infty.
\]
Consequently, $\|\Psi(t)+q(t)\|_{L^\infty(\Omega)}\to0$, and therefore
\begin{equation}\label{eq:relative-chemical-decay}
\left\|\frac{w(\cdot,t)}{W}-1\right\|_{L^\infty(\Omega)}
= \left\|e^{-(\Psi(\cdot,t)+q(\cdot,t))}-1\right\|_{L^\infty(\Omega)}
\le e^{\|\Psi(t)+q(t)\|_{L^\infty(\Omega)}}\|\Psi(t)+q(t)\|_{L^\infty(\Omega)} \to 0
\end{equation}
as $t\to\infty$. This completes the proof of Theorem~\ref{thm:main}.
\end{proof}

\appendix
\section{Weighted inequalities and mode reconstruction}\label{app:weighted}
In this appendix, we collect the weighted functional inequalities on $\Omega=\mathbb{R}_+\times\mathbb{T}^2$ and $\mathbb{R}_+$ used throughout the paper.

\subsection{Poincar\'e, Hardy, and trace inequalities}
\begin{Lem}\label{lem:standard-weighted-tools}
Let $\omega=\omega(x_1)>0$, and let $f$ be a smooth function for which the norms below are finite. Then
\begin{equation}\label{eq:weighted-transversal-poincare}
\int_\Omega \omega(x_1)| P_{\neq} f|^2\dd x \le C\int_\Omega\omega(x_1)|\nabla_{x'}f|^2\dd x.
\end{equation}
If $f=f(x_1)$ satisfies $\lim_{x_1\to\infty}f(x_1)=0$, then
\begin{equation}\label{eq:profile-hardy-appendix}
\int_0^\infty|f(x_1)|^2\dd x_1 \le C\int_0^\infty\frac{|f_{x_1}(x_1)|^2}{U(x_1)}\dd x_1.
\end{equation}
Moreover, the boundary trace at $x_1=0$ satisfies
\begin{equation}\label{eq:one-dimensional-trace-Omega}
\int_{\mathbb{T}^2}|f(0,x')|^2\dd x' \le 2\|f\|_{L^2(\Omega)}\|f_{x_1}\|_{L^2(\Omega)},
\end{equation}
and for $2\le r\le6$,
\[
\|f\|_{L^r(\Omega)} \le C_r\|f\|_{H^1(\Omega)}, \qquad \|f\|_{L^\infty(\Omega)} \le C\|f\|_{H^2(\Omega)}.
\]
\end{Lem}
\begin{proof}
Applying the Poincar\'e--Wirtinger inequality on $\mathbb{T}^2$ for each fixed $x_1>0$, multiplying by $\omega(x_1)$, and integrating over $x_1\in\mathbb{R}_+$ yields \eqref{eq:weighted-transversal-poincare}. Since $1/U(x_1)\sim (1+x_1)^2$, the Hardy inequality
\[
\int_0^\infty|f(x_1)|^2\dd x_1 \le 4\int_0^\infty(1+x_1)^2|f_{x_1}(x_1)|^2\dd x_1
\]
for $f(\infty)=0$ yields \eqref{eq:profile-hardy-appendix}. The trace bound \eqref{eq:one-dimensional-trace-Omega} follows from
\[
|f(0,x')|^2 = -2\int_0^\infty f(x_1,x')f_{x_1}(x_1,x')\dd x_1
\]
and the Cauchy--Schwarz inequality, while the Sobolev embeddings on $\Omega=\mathbb{R}_+\times\mathbb{T}^2$ are standard~\cite{AdamsFournier2003,Brezis2011}.
\end{proof}

\subsection{Bounds for the modal potentials}
\begin{Lem}\label{lem:mode-reconstruction}
Let $t\in[0,T]$ and assume $N(t)\le\delta_0\le1$. Then
\begin{equation}\label{eq:mode-reconstruction-phi-appendix}
\int_\Omega\frac{|\phi|^2}{U}\dd x \le \int_0^\infty\frac{\Phi_{x_1}^2}{U}\dd x_1 + C\int_\Omega\frac{|\nabla_{x'}\phi|^2}{U}\dd x,
\end{equation}
and
\begin{equation}\label{eq:first-physical-psi-weighted-reconstruction}
\int_\Omega\frac{|\psi|^2}{U}\dd x \le \int_0^\infty\frac{\Psi_{x_1}^2}{U}\dd x_1 + C\int_\Omega\frac{|\nabla\psi|^2}{U}\dd x.
\end{equation}
Moreover, the nonzero-mode potential $q$ and normal component $\psi_1$ satisfy
\begin{equation}\label{eq:weighted-potential-estimates-concise}
\begin{aligned}
&\| q\|_{L^\infty(\Omega)} + \sup_{x_1>0}\left\|\frac{q(x_1,\cdot)}{\sqrt{U(x_1)}}\right\|_{L^2(\mathbb{T}^2)}
+\sup_{x_1>0}\left\|\frac{\psi_1(x_1,\cdot)}{\sqrt{U(x_1)}}\right\|_{L^2(\mathbb{T}^2)}\\
&\quad +\sup_{x_1>0}\left\|\frac{P_{\neq}\psi_1(x_1,\cdot)}{\sqrt{U(x_1)}}\right\|_{L^2(\mathbb{T}^2)} + \| q\|_{\omega_2}\\
&\le C\left(\int_0^\infty\frac{\Psi_{x_1}^2+\Psi_{x_1x_1}^2}{U}\dd x_1+\int_\Omega\frac{|\nabla\psi|^2+|\nabla^2\psi|^2}{U}\dd x\right)^{1/2} \le CN(t),
\end{aligned}
\end{equation}
\begin{equation}\label{eq:nonzero-potential-H2-bound}
\| q\|_{2,\omega_2} + \left\|\frac{q}{\sqrt{U}}\right\|_{H^2(\Omega)} \le C\|\nabla\psi\|_{1,\omega_2} \le CN(t),
\end{equation}
and
\begin{equation}\label{eq:weighted-vector-Linfty}
\left\|\frac{\psi}{\sqrt{U}}\right\|_{L^\infty(\Omega)} \le C\|\psi\|_{2,\omega_2} \le CN(t).
\end{equation}
\end{Lem}
\begin{proof}
Since $\omega_2=1/U(x_1)$ depends only on $x_1$, transversal orthogonality and \eqref{eq:weighted-transversal-poincare} give \eqref{eq:mode-reconstruction-phi-appendix} and \eqref{eq:first-physical-psi-weighted-reconstruction}. Next, since $P_0 q=0$, $\nabla q= P_{\neq}\psi$, and $\nabla_{x'} q=(\psi_2,\psi_3)$, applying \eqref{eq:weighted-transversal-poincare} to $q$ and $\nabla q$ yields $\| q\|_{2,\omega_2}\le C\|\nabla\psi\|_{1,\omega_2}$. Since $U^{-1/2}(x_1)$ is a first-degree polynomial in $x_1$ satisfying $|(U^{-1/2})_{x_1}|\le CU^{-1/2}$ and $(U^{-1/2})_{x_1x_1}=0$, the product rule gives
\[
\left\|\frac{f}{\sqrt{U}}\right\|_{H^k(\Omega)} \le C\|f\|_{k,\omega_2}, \qquad k=1,2,
\]
for any scalar or vector function $f$. Taking $f= q$ and $f=\psi$, and applying the Sobolev embeddings $H^2(\Omega)\hookrightarrow L^\infty(\Omega)$ and $H^1(\mathbb{R}_+;L^2(\mathbb{T}^2))\hookrightarrow L^\infty(\mathbb{R}_+;L^2(\mathbb{T}^2))$, we obtain \eqref{eq:weighted-potential-estimates-concise}--\eqref{eq:weighted-vector-Linfty}.
\end{proof}

In particular, Lemma~\ref{lem:mode-reconstruction} and the definition of $N(t)$ in \eqref{eq:bootstrap-norm} imply the $L^\infty$ bounds used throughout Section~\ref{sec:apriori}:
\begin{equation}\label{eq:H2-Linfty-smallness}
\|\phi\|_{L^\infty(\Omega)} + \|\psi\|_{L^\infty(\Omega)} + \|\Psi\|_{L^\infty(\mathbb{R}_+)} \le CN(t),
\end{equation}
\begin{equation}\label{eq:exponential-bounds}
\left|e^{-p q}-1\right| \le C| q|, \qquad \left|e^{-p q}\right|+\left|e^{-p\Psi}\right| \le C,
\end{equation}
\begin{equation}\label{eq:Psi-Linfty-smallness}
\|\Psi\|_{L^\infty(\mathbb{R}_+)} \le CN(t),
\end{equation}
and
\begin{equation}\label{eq:Psi-x-weighted-Linfty}
\left\|\frac{\Psi_{x_1}}{\sqrt{U}}\right\|_{L^\infty(\mathbb{R}_+)}
\le C\left(\int_0^\infty\frac{\Psi_{x_1}^2+\Psi_{x_1x_1}^2}{U}\dd x_1\right)^{1/2} \le CN(t).
\end{equation}

\subsection{Weighted interpolation}
\begin{Lem}\label{lem:weighted-half-cylinder-interpolation}
Let $f$ be a smooth function on $\Omega$ with finite weighted norms. Then for every $\eta>0$, there exists $C_\eta>0$ such that
\begin{equation}\label{eq:weighted-first-interpolation-cross}
\int_\Omega\frac{|\nabla f|^2}{U}\dd x \le \eta\int_\Omega\frac{W^p|\nabla^2f|^2}{U}\dd x + C_\eta\int_\Omega\frac{|f|^2}{W^pU}\dd x,
\end{equation}
\begin{equation}\label{eq:weighted-first-interpolation-same}
\int_\Omega\omega_j|\nabla f|^2\dd x \le \eta\int_\Omega\omega_j|\nabla^2f|^2\dd x + C_\eta\int_\Omega\omega_j|f|^2\dd x, \qquad j=1,2,3,
\end{equation}
and
\begin{equation}\label{eq:weighted-mixed-interpolation}
\int_\Omega\omega_j|\partial_{x_1}\nabla_{x'}f|^2\dd x
\le C\int_\Omega\omega_j\left(|f_{x_1x_1}|^2+\sum_{|\beta|=2}|\partial_{x'}^\beta f|^2+|f|^2\right)\dd x, \qquad j=2,3.
\end{equation}
\end{Lem}
\begin{proof}
To avoid any boundary trace requirement at $x_1=0$, we use the second-order reflection extension $g=Ef$ on the full cylinder $\mathbb{R}\times\mathbb{T}^2$:
\[
g(x_1,x') =
\begin{cases}
f(x_1,x'), & x_1\ge0,\\
3f(-x_1,x')-2f(-2x_1,x'), & x_1<0,
\end{cases}
\]
and extend each algebraic weight evenly by $\omega_j(x_1)=\omega_j(|x_1|)$. Since $\omega_j(x_1)$ is a power weight in $1+cx_1$, we have $\omega_j(2x_1)\sim\omega_j(x_1)$, $|\omega_j'|\le C\omega_j$ almost everywhere on $\mathbb{R}$, $\omega_2^2=\omega_1\omega_3$, and
\[
\int_{\mathbb{R}\times\mathbb{T}^2}\omega_j|\partial^\gamma g|^2\dd x \le C\int_\Omega\omega_j|\partial^\gamma f|^2\dd x, \qquad |\gamma|\le2.
\]
Integrating by parts on $\mathbb{R}\times\mathbb{T}^2$ for any $\omega\in\{\omega_1,\omega_2,\omega_3\}$ gives
\[
\int_{\mathbb{R}\times\mathbb{T}^2}\omega|\nabla g|^2\dd x = -\int_{\mathbb{R}\times\mathbb{T}^2}\omega g\Delta g\dd x - \int_{\mathbb{R}\times\mathbb{T}^2}\omega' g g_{x_1}\dd x.
\]
Using $|\omega'|\le C\omega$ and Young's inequality proves \eqref{eq:weighted-first-interpolation-same}, while applying Cauchy--Schwarz with $\omega_2=\sqrt{\omega_1\omega_3}$ to $\int_{\mathbb{R}\times\mathbb{T}^2}\omega_2|g||\Delta g|\dd x$ and using $\omega_2\le C\omega_1$ proves \eqref{eq:weighted-first-interpolation-cross}. Finally, for $k\in\{2,3\}$, integrating by parts in $x_k$ and $x_1$ on $\mathbb{R}\times\mathbb{T}^2$ gives
\[
\int_{\mathbb{R}\times\mathbb{T}^2}\omega|g_{x_1x_k}|^2\dd x
= \int_{\mathbb{R}\times\mathbb{T}^2}\omega g_{x_1x_1}g_{x_kx_k}\dd x + \int_{\mathbb{R}\times\mathbb{T}^2}\omega' g_{x_1}g_{x_kx_k}\dd x,
\]
which combined with \eqref{eq:weighted-first-interpolation-same} yields \eqref{eq:weighted-mixed-interpolation}.
\end{proof}

\section{Weighted Neumann inversion for the nonzero modes}\label{app:neumann}
\begin{Lem}\label{lem:weighted-nonzero-neumann}
Let $f\in L^2_{\omega_2}(\Omega)$ satisfy $P_0 f=0$. Then the weak Neumann problem
\begin{equation}\label{eq:L46-neumann}
-(\alpha-1)p\,\nabla\cdot(W^pU\nabla z) = f \quad\text{in }\Omega, \qquad z_{x_1}(0,x')=0, \qquad  P_0 z=0,
\end{equation}
admits a unique weak solution $z\in H_{\neq}^1$, where $H_{\neq}^1$ is the completion of smooth zero-transversal-mean functions under the inner product $\langle z_1,z_2\rangle_{H_{\neq}^1}=(\alpha-1)p\int_\Omega W^pU\nabla z_1\cdot\nabla z_2\dd x$. Moreover,
\begin{equation}\label{eq:L46-neumann-energy}
0 \le \int_\Omega fz\dd x = (\alpha-1)p\int_\Omega W^pU|\nabla z|^2\dd x \le C\int_\Omega\frac{|f|^2}{W^pU}\dd x,
\end{equation}
and
\begin{equation}\label{eq:L46-neumann-flux}
\int_\Omega W^{2p}U^2|\nabla z|^2\dd x \le C\int_\Omega\frac{|f|^2}{U}\dd x.
\end{equation}
If $f(\tau)$ is differentiable in time in $(H_{\neq}^1)^*$, then $z(\tau)$ satisfies
\begin{equation}\label{eq:L46-neumann-selfadjoint-time}
\frac{1}{2}\frac{\dd}{\dd\tau}\int_\Omega fz\dd x = \int_\Omega f_\tau z\dd x.
\end{equation}
\end{Lem}
\begin{proof}
By the weighted transversal Poincar\'e inequality \eqref{eq:weighted-transversal-poincare}, any $h\in H_{\neq}^1$ satisfies
\[
\int_\Omega W^pU|h|^2\dd x \le C_P\int_\Omega W^pU|\nabla_{x'}h|^2\dd x \le C_P\int_\Omega W^pU|\nabla h|^2\dd x.
\]
Hence by Cauchy--Schwarz,
\[
\left|\int_\Omega fh\dd x\right| \le \left(\int_\Omega\frac{|f|^2}{W^pU}\dd x\right)^{1/2}\left(\int_\Omega W^pU|h|^2\dd x\right)^{1/2}
\le C\|f\|_{\omega_1}\|h\|_{H_{\neq}^1}.
\]
The Lax--Milgram theorem therefore yields a unique weak solution $z\in H_{\neq}^1$ satisfying \eqref{eq:L46-neumann-energy}.

To prove the flux bound \eqref{eq:L46-neumann-flux} without any smallness condition on the steady mass $\lambda>0$, note from \eqref{eq:explicit-steady-profile} that
\[
W^p(x_1)U(x_1) = a_0(1+cx_1)^{-m}, \qquad m = 2+\frac{2}{\alpha-1},
\]
for fixed positive constants $a_0,c>0$. Choose a small constant $\sigma>0$ such that $m\sigma\sqrt{C_P}\le 1/4$, and define the modified weight
\[
\mu(x_1) = (1+\sigma x_1)^{-m}.
\]
Since $1+\sigma x_1\sim 1+cx_1$ on $\mathbb{R}_+$ (with constants depending on $\lambda$), we have $\mu(x_1)\sim W^p(x_1)U(x_1)$ on $\mathbb{R}_+$, while the logarithmic derivative of $\mu$ satisfies $|\mu'(x_1)/\mu(x_1)|\le m\sigma$ for all $x_1\ge0$. Testing the weak formulation of \eqref{eq:L46-neumann} with the admissible nonzero-mode test function $h=\mu(x_1)z$ gives
\[
(\alpha-1)p\int_\Omega W^pU\mu|\nabla z|^2\dd x
= \int_\Omega f\mu z\dd x - (\alpha-1)p\int_\Omega W^pU\mu' z z_{x_1}\dd x.
\]
Using $|\mu'|\le m\sigma\mu$ and $\int_\Omega W^pU\mu|z|^2\dd x\le C_P\int_\Omega W^pU\mu|\nabla z|^2\dd x$, we estimate
\[
\left|\int_\Omega W^pU\mu' z z_{x_1}\dd x\right|
\le m\sigma\sqrt{C_P}\int_\Omega W^pU\mu|\nabla z|^2\dd x
\le \frac{1}{4}\int_\Omega W^pU\mu|\nabla z|^2\dd x,
\]
and by Cauchy--Schwarz and Young's inequality,
\[
\begin{aligned}
\left|\int_\Omega f\mu z\dd x\right|
&\le \left(\int_\Omega\frac{\mu}{W^pU}|f|^2\dd x\right)^{1/2}\left(\int_\Omega W^pU\mu|z|^2\dd x\right)^{1/2}\\
&\le \frac{(\alpha-1)p}{4}\int_\Omega W^pU\mu|\nabla z|^2\dd x + C\int_\Omega\frac{\mu}{W^pU}|f|^2\dd x.
\end{aligned}
\]
Absorbing the gradient terms into the left-hand side and using $\mu\sim W^pU$ and $U\le C$, we conclude that
\[
\begin{aligned}
\int_\Omega W^{2p}U^2|\nabla z|^2\dd x
&\le C\int_\Omega W^pU\mu|\nabla z|^2\dd x
\le C\int_\Omega\frac{\mu}{W^pU}|f|^2\dd x\\
&\le C\int_\Omega|f|^2\dd x
\le C\int_\Omega\frac{|f|^2}{U}\dd x,
\end{aligned}
\]
which proves \eqref{eq:L46-neumann-flux}. Finally, \eqref{eq:L46-neumann-selfadjoint-time} follows directly from the symmetry and time-independence of the bilinear form $\langle\cdot,\cdot\rangle_{H_{\neq}^1}$.
\end{proof}

\medskip
\noindent\textbf{Acknowledgments and AI Disclosure.} The research of Lin Chang was supported in part by the National Natural Science Foundation of China Grant No. 12401272 and the Natural Science Foundation of Hebei Province Grant No.A2024109005. The first author thanks Feimin Huang, Yi Wang and Yong Wang for helpful suggestions during a visit to the  Academy of Mathematics and Systems Science, Beijing, in May 2026, and gratefully acknowledges the Institute's hospitality. The authors had helpful discussions with OpenAI's ChatGPT and Google's Gemini during the preparation and revision of this manuscript. The authors take full personal responsibility for the correctness of all mathematical proofs, derivations, and results presented in this paper.

\end{document}